\documentclass[mnsc,nonblindrev,copyedit]{informs3a}

\usepackage[utf8]{inputenc}
\usepackage[letterpaper,margin=1in]{geometry}
\let\OriginalManuscriptNormalSize\normalsizeXIyesCE
\renewcommand\normalsizeXIyesCE{\OriginalManuscriptNormalSize\fontsize{11}{19.71}\selectfont\abovedisplayskip=5pt plus 2pt minus 1pt\belowdisplayskip=5pt plus 2pt minus 1pt\abovedisplayshortskip=2pt plus 1pt minus 1pt\belowdisplayshortskip=4pt plus 1pt minus 1pt}
\let\OriginalManuscriptSmallSize\smallXyesCE
\renewcommand\smallXyesCE{\OriginalManuscriptSmallSize\fontsize{10}{15}\selectfont}
\let\OriginalManuscriptFootnoteSize\footnotesizeIXyesCE
\renewcommand\footnotesizeIXyesCE{\OriginalManuscriptFootnoteSize\fontsize{9}{13.5}\selectfont}
\renewcommand\ABSfont{\fontsize{9}{18.75}\selectfont}
\AtBeginDocument{\newgeometry{margin=1in,hoffset=0pt,voffset=0pt}\normalsize\raggedbottom\predisplaypenalty=10000}

\usepackage{amsmath,amssymb,mathtools}
\usepackage{booktabs}
\usepackage{array}
\usepackage{enumerate}
\usepackage[dvipsnames]{xcolor}
\usepackage{graphicx}
\usepackage{placeins}
\usepackage{tikz}
\usetikzlibrary{arrows.meta,shapes.geometric}
\usepackage{pgfplots}
\pgfplotsset{compat=1.13}
\usepackage{natbib}
\usepackage[backref=false,bookmarks,breaklinks=true,colorlinks=true,hypertexnames=false,plainpages=false,citecolor=DarkBlue,urlcolor=DarkBlue,filecolor=DarkBlue,linkcolor=DarkBlue,pdftitle={A 5/2 Guarantee for Echelon Stock (s,S) Policies in One-Warehouse Multi-Retailer Systems},pdfauthor={Ming Hu},pdfkeywords={multi-echelon inventory; periodic review; fixed ordering costs; stochastic demand; supplier lead time; approximation guarantees}]{hyperref}
\urldef{\provemeurl}\url{https://prove2.me/missions/A_%245%2F2%24_Guarantee_for_Echelon_Stock_%24(s%2CS)%24_Policies_in_One-Warehouse_Multi-Retailer_Systems}
\definecolor{DarkBlue}{rgb}{0,0.08,0.45}
\definecolor{lsTeal}{rgb}{0.1567,0.4717,0.4234}
\definecolor{lsOrng}{rgb}{0.8436,0.4582,0.1704}
\definecolor{lsPale}{rgb}{0.9470,0.7987,0.6011}
\definecolor{lsBand}{rgb}{0.9550,0.8924,0.7724}

\bibpunct[, ]{(}{)}{,}{a}{}{,}
\def\bibfont{\fontsize{11}{14.5}\selectfont}

\setlength{\floatsep}{4pt plus 2pt minus 2pt}
\setlength{\textfloatsep}{4pt plus 2pt minus 2pt}
\setlength{\intextsep}{4pt plus 2pt minus 2pt}
\TheoremsNumberedThrough
\EquationsNumberedThrough
\MANUSCRIPTNO{}
\renewcommand{\Halmos}{\unskip\nobreak\mbox{\quad$\square$}}

\newcommand{\E}{\mathbb E}

\newcommand{\R}{\mathbb R}
\newcommand{\OPT}{\operatorname{OPT}}

\newcommand{\cL}{\mathcal L}
\newcommand{\Prob}{\mathbb P}
\newcommand{\J}{\mathcal J}

\newcommand{\nt}[1]{{#1}}

\allowdisplaybreaks
\newcommand{\ind}{\mathbf1}
\newcommand{\osc}{\operatorname{osc}}
\newcommand{\NV}{\operatorname{NV}}

\DeclareUnicodeCharacter{03BC}{\ensuremath{\mu}}
\DeclareUnicodeCharacter{03C6}{\ensuremath{\phi}}
\DeclareUnicodeCharacter{03A6}{\ensuremath{\Phi}}
\DeclareUnicodeCharacter{03C0}{\ensuremath{\pi}}
\DeclareUnicodeCharacter{03BA}{\ensuremath{\kappa}}
\DeclareUnicodeCharacter{03B7}{\ensuremath{\eta}}
\DeclareUnicodeCharacter{03B5}{\ensuremath{\epsilon}}
\DeclareUnicodeCharacter{03B4}{\ensuremath{\delta}}
\DeclareUnicodeCharacter{03B8}{\ensuremath{\theta}}
\DeclareUnicodeCharacter{03B1}{\ensuremath{\alpha}}
\DeclareUnicodeCharacter{03B2}{\ensuremath{\beta}}
\DeclareUnicodeCharacter{03C1}{\ensuremath{\rho}}
\DeclareUnicodeCharacter{03BB}{\ensuremath{\lambda}}
\DeclareUnicodeCharacter{0393}{\ensuremath{\Gamma}}
\DeclareUnicodeCharacter{03B3}{\ensuremath{\gamma}}
\DeclareUnicodeCharacter{2113}{\ensuremath{\ell}}
\DeclareUnicodeCharacter{2264}{\ensuremath{\le}}
\DeclareUnicodeCharacter{2265}{\ensuremath{\ge}}
\DeclareUnicodeCharacter{221E}{\ensuremath{\infty}}
\DeclareUnicodeCharacter{2192}{\ensuremath{\to}}
\DeclareUnicodeCharacter{21D2}{\ensuremath{\Rightarrow}}
\DeclareUnicodeCharacter{2208}{\ensuremath{\in}}
\DeclareUnicodeCharacter{221A}{\ensuremath{\surd}}
\DeclareUnicodeCharacter{2211}{\ensuremath{\sum}}
\DeclareUnicodeCharacter{222B}{\ensuremath{\int}}
\DeclareUnicodeCharacter{2112}{\ensuremath{\mathcal L}}
\DeclareUnicodeCharacter{1D49C}{\ensuremath{\alpha}}
\DeclareUnicodeCharacter{2193}{\ensuremath{\downarrow}}
\DeclareUnicodeCharacter{2202}{\ensuremath{\partial}}
\DeclareUnicodeCharacter{2260}{\ensuremath{\ne}}
\DeclareUnicodeCharacter{21A6}{\ensuremath{\mapsto}}
\DeclareUnicodeCharacter{00B2}{\ensuremath{{}^{2}}}
\DeclareUnicodeCharacter{00B3}{\ensuremath{{}^{3}}}
\DeclareUnicodeCharacter{2074}{\ensuremath{{}^{4}}}
\DeclareUnicodeCharacter{2075}{\ensuremath{{}^{5}}}
\DeclareUnicodeCharacter{2076}{\ensuremath{{}^{6}}}
\DeclareUnicodeCharacter{2077}{\ensuremath{{}^{7}}}
\DeclareUnicodeCharacter{2078}{\ensuremath{{}^{8}}}
\DeclareUnicodeCharacter{2080}{\ensuremath{{}_{0}}}
\DeclareUnicodeCharacter{2081}{\ensuremath{{}_{1}}}
\DeclareUnicodeCharacter{2082}{\ensuremath{{}_{2}}}
\DeclareUnicodeCharacter{208B}{\ensuremath{{}_{-}}}
\usepackage{mathrsfs}

\allowdisplaybreaks
\begin{document}

\RUNAUTHOR{Hu}

\RUNTITLE{A $5/2$ Guarantee for Echelon Stock $(s,S)$ Policies in OWMR Systems}

\TITLE{\Large \textcolor{black}{A $5/2$ Guarantee for Echelon Stock $(s,S)$ Policies\\in One-Warehouse Multi-Retailer Systems}}

\ARTICLEAUTHORS{\AUTHOR{Ming Hu$^{*}$} \AFF{Rotman School of Management, University of Toronto, Toronto, Ontario M5S 3E6, Canada\\\EMAIL{ming.hu@rotman.utoronto.ca}}}

\ABSTRACT{We study a periodic-review warehouse serving multiple heterogeneous retailers with stochastic demand and fixed costs and lead times at both echelons. We construct a deterministic echelon stock $(s,S)$ policy whose long-run expected average cost is at most $5/2$ times the optimum over all nonanticipative policies. The guarantee holds when, as in a retail chain with common service penalties, the uniform backlog penalty and heterogeneous retailer holding rates are at least the warehouse holding rate and heterogeneous retailer setup costs are at most $3/2$ times the supplier's. The latter setup-cost condition extends beyond supplier setup-cost dominance: a retailer dispatch may incur a setup cost up to 50 percent higher than a supplier order's. Demand vectors are i.i.d.\ over periods with finite positive means, and demands of different retailers may be correlated. Each retailer selects the full cycle of its single-location plan or a shifted cycle of half its width, under an inflated setup cost, and the warehouse serves requests in arrival order, allowing partial shipments. The upper bound on the policy's cost includes each retailer's mirror cost, which prices the warehouse imbalance created by that retailer's batching. The upper and lower bounds share a pooled supplier problem and add one problem per retailer, giving a factor independent of all lead times and the number of retailers. For arbitrary nonnegative setup costs and associated retailer demands, a consolidation policy holds allocated stock until a dispatch threshold is met and gives a factor below $4.9$. On 864 instances built from Walmart store sales, the simulated cost of an integer-rounded implementation of the FIFO policy averages 1.45 times the lower bound and reaches 1.75, and our analytical upper bound on its cost is below 1.94 times the lower bound; both ratios to the lower bound have nearly the same mean at every supplier lead time.}

\def\HOOKa{\vspace*{-6pt}}\def\HOOKb{\vspace*{-8pt}}%
\maketitle
\begingroup
\renewcommand{\thefootnote}{\fnsymbol{footnote}}%
\footnotetext[1]{The author used AI tools (such as ChatGPT and Claude) to assist with derivations, numerical verification, and editing. The author developed and verified the substantive research questions, theoretical framework, analysis, and results. The author takes sole responsibility for the paper's content. The main results are formally verified in Lean~4; the formal proofs are public at \provemeurl. Their formal statements, including the definitions of the model and of the policies, match those of the paper, and these proofs rely only on Lean's standard axioms.}%
\endgroup
\setlength{\abovedisplayskip}{5pt plus 2pt minus 1pt}
\setlength{\belowdisplayskip}{5pt plus 2pt minus 1pt}
\setlength{\abovedisplayshortskip}{2pt plus 1pt minus 1pt}
\setlength{\belowdisplayshortskip}{4pt plus 1pt minus 1pt}
\vspace{-2em}
\section{Introduction}\label{sec:intro}

A warehouse serving several retailers must coordinate supplier purchases, retailer dispatches, and inventory at both echelons. Fixed costs favor fewer, larger orders and holding costs favor smaller, more frequent ones, while uncertain demand and lead times keep the echelons from timing orders together. An optimal policy decides when to purchase, which retailers to serve, and how much to send, using the joint inventory and pipeline history. These interactions make stochastic one-warehouse multi-retailer (OWMR) systems hard even with independent demands and divisible quantities. Existing results optimize or bound a policy within a prescribed policy class, such as fixed reorder intervals \citep{ChuShen2010,ShangTaoZhou2015}, or give bounds that depend on problem primitives and become tight in a limit \citep{GallegoEtAl2007,ZhuChenHuYang2021}.

Two challenges stand in the way of a constant-factor guarantee. The first comes from the supplier lead time: purchases cannot immediately clear pending requests, so the warehouse must position stock against uncertain aggregate demand, and a retailer request may span several warehouse receipts and incur a setup charge for each portion dispatched separately. The second is warehouse imbalance, a batching-induced mismatch between stock on hand and retailer orders: independently chosen retailer replenishment cycles may let stock wait for an order, or an order wait for stock. A uniform guarantee must address both challenges at all finite lead times and network sizes.

We provide a uniform $5/2$ guarantee for a deterministic policy with first-in, first-out (FIFO) warehouse fulfillment in periodic-review OWMR systems with fixed costs at both echelons, relative to the optimum over all nonanticipative policies, including history-dependent and randomized ones. The guarantee holds at every finite supplier and retailer lead time, for any number of retailers and nonnegative demand vectors that are i.i.d.\ over periods with finite positive component means, including demand correlated across retailers within a period. Demand may have point masses and unbounded support. We focus on a retail chain with a common charge for delayed service and with retailer holding rates and the backlog penalty at least as high as the warehouse holding rate. These assumptions represent shared service standards and relatively inexpensive warehouse storage. Each retailer setup cost is at most $3/2$ times the supplier's. This condition includes supplier setup-cost dominance, a plausible regime when each inbound order is a long-haul shipment, by truck from a national distribution center or by container from overseas, with order processing, carrier booking, and receiving, while each store dispatch is a short haul from the regional warehouse to a local store. It also allows a retailer dispatch to cost up to 50 percent more than a supplier order. The restriction applies to each retailer separately, so the combined setup cost of dispatching to all retailers may exceed the setup cost of a supplier order by any amount in a large network.

Our FIFO policy separates retailer planning from physical replenishment. Each retailer first solves a single-location problem at an inflated planning setup cost, $5/3$ times its actual setup cost. It then follows the resulting optimal $(\hat s_i,\hat S_i)$ cycle (its full window $[\hat s_i,\hat S_i]$) or the best shifted cycle of half that width, whichever has the lower augmented cost, which adds the mirror cost defined below to the retailer's own cost. Here $\hat S_i$ is the order-up-to level, or target. The selected policy generates orders from observed demand, ignoring warehouse stock; each positive order enters the warehouse queue as a request. The warehouse purchases according to one $(s,S)$ rule on an aggregate inventory position, its echelon inventory position less a constant, and serves requests in FIFO order, allowing a partial final request when stock runs out and continuing it at later receipts until the request is filled.\footnote{Modified echelon $(r,Q)$ policies \citep{HuYang2014,ZhuChenHuYang2021} also ship part of an order when stock is short, but a shipment lifting the downstream position above its reorder point ends the order and cancels the rest.} This $(s,S)$ rule is optimal for a pooled supplier problem, a single-location problem facing aggregate demand across retailers over the supplier lead time at the warehouse holding rate and the backlog penalty, with the supplier setup cost raised by the largest retailer setup cost to pay for these continuations. Thus, our FIFO policy is an echelon stock $(s,S)$ policy. All choices are deterministic functions of the model, including the fulfillment sequence of requests placed in the same period, which serves lower-indexed retailers first, so no initial randomization is required.

\begin{figure}[t]
\SingleSpaced
\centering
\includegraphics[width=.99\linewidth]{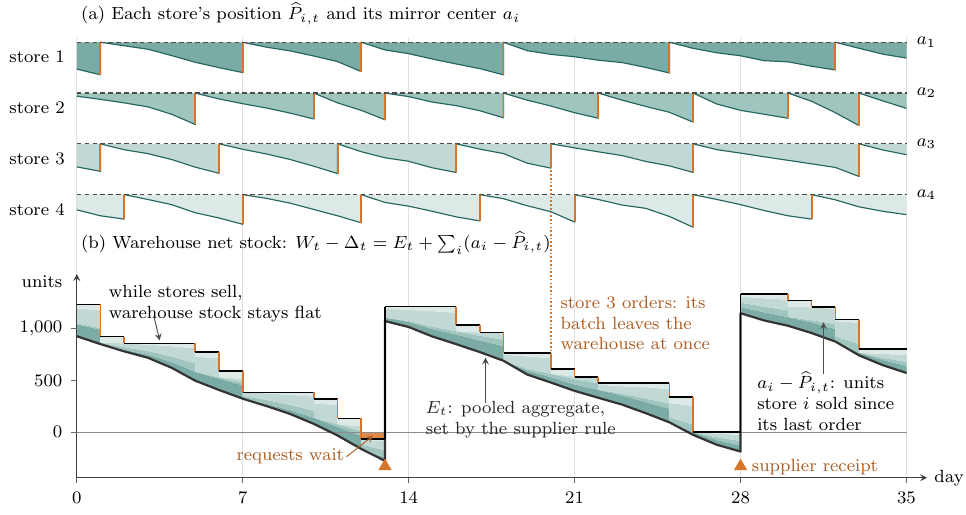}
\caption{The mirror view of the warehouse on Walmart sales.}
\label{fig:intro-mirror}
\begin{minipage}{0.98\linewidth}
{\SingleSpaced\footnotesize\emph{Notes.} Our policy on 35 days of actual sales of \texttt{FOODS\_3\_587} at the four Walmart California stores, with a supplier lead time of three days, retailer lead times of one day, and setup costs corresponding to retailer and warehouse EOQ cycles of seven and fourteen days. Since the backlog penalty is about twelve times the warehouse holding rate, each mirror center is the store's order-up-to level, so the shaded height $a_i-\widehat P_{i,t}$ is what store $i$ has sold since its last order. The position $\widehat P_{i,t}$ includes net inventory, shipments in transit, and requests awaiting dispatch. Panel (b) stacks these heights on the pooled aggregate $E_t$, and their top is the warehouse net stock in each period.}
\end{minipage}
\end{figure}

Our analysis rests on a mirror view of the warehouse. A retailer's order reduces warehouse net stock, defined as stock less requests awaiting dispatch, by its whole batch at once,\footnote{If stock runs out, the unfilled part of a request waits in the FIFO queue until later supplier receipts, so net stock still falls by the full batch size and becomes negative; dispatching the waiting part later leaves net stock unchanged.} while its later sales leave that net stock unchanged. Relative to a pooled aggregate that the supplier rule controls, warehouse net stock is the mirror image of the retailer's sawtooth position about a fixed level, its mirror center (Figure~\ref{fig:intro-mirror}): lowest just after an order and highest late in the cycle. We charge each retailer for its reflection at the warehouse rate, the backlog penalty $p$ while the retailer is above its center and the holding rate $h$ while it is below, and call the smallest long-run average charge its mirror cost. We use the best center, which solves a newsvendor problem with the retailer's inventory position as demand, so it is a $p/(h+p)$ quantile of the position's long-run distribution.\footnote{In Figure~\ref{fig:intro-mirror}, each center is the store's order-up-to level because cycles average fewer than $1+p/h\approx13$ days; with longer cycles or a smaller penalty, the center may lie inside the window. Its location also depends on how much time the post-order position spends at the target, where it stays through every period without demand.} Allocation relaxations for a depot without stock \citep{EppenSchrage1981,FedergruenZipkin1984} remove this imbalance by rebalancing retailer inventories; in contrast, the mirror cost prices the batching-induced imbalance between warehouse net stock and the pooled position. Base stock never swings and has zero mirror cost, so the charge falls only on retailers that order in batches.

The guarantee $5/2$ is $1+3/2$. The lower bound splits $5/2$ copies of any policy's cost: one copy goes to a pooled supplier problem over the supplier lead time at setup cost $5K_0/2$, and the other $3/2$ copies to each retailer's single-location problem at setup cost $(5/2)K_i/(3/2)=5K_i/3$, so retailers plan at that cost. Here $K_0$ is the supplier setup cost, $K_i$ is retailer $i$'s setup cost, and $K_{\max}=\max_iK_i$. Under FIFO service, at most one previously started request, the partly filled request at the head of the queue, can continue at a supplier receipt, and we assign its extra setup charge, at most $K_{\max}$, to that receipt's supplier order. An inventory identity splits the warehouse's net stock into a pooled aggregate and the mirror images of the retailer positions, so the policy's cost is at most the same supplier problem at setup cost $K_0+K_{\max}\le5K_0/2$ plus each retailer's augmented cost. Since both bounds contain the same pooled supplier problem, the supplier lead time does not enlarge the factor. The bounds also use only each retailer's demand distribution and the distribution of aggregate demand, so the guarantee allows demand correlated across stores.

On the retailer side, the upper bound must match the lower bound's $3/2$ copies of each retailer's single-location problem at setup cost $5K_i/3$. Therefore, we construct the policy, and bound its cost, with each retailer's planning setup cost raised by two thirds, to $5K_i/3$. The larger planning setup cost can produce larger batches and hence larger swings in warehouse net stock, which raise the mirror cost charged to the retailer. Each retailer selects the full cycle of its plan or the best shifted half-width window of that cycle, whichever has the smaller augmented cost.\footnote{In the integer-rounded implementation used in our Walmart experiments (Section~\ref{sec:numerics}), every retailer selects the half window. Its savings in retailer inventory costs and mirror cost more than offset its additional setup cost: the full window's augmented cost, evaluated at actual setup costs, is 5 to 17 percent higher than the half window's.} The selection balances the setup costs that larger batches save against the larger swings and mirror cost they create. Two observations explain why the full cycle and its shifted half suffice. First, shortening a cycle can raise its order frequency while reducing holding and backlog costs. Some shifted half of an optimal cycle preserves at least half of its expected cost savings. The shift matters because the favorable part of the cycle need not begin at its original target. Second, a wider window lets the inventory position swing more, which raises its mirror-cost bound. Halving the window lowers this bound but may give up savings in holding and backlog costs; these lost savings and the half window's mirror cost cannot both be large. Together, the two observations bound the better candidate's cost, mirror cost included, by $3/2$ times the retailer's optimal single-location cost at its inflated planning setup cost, which matches the retailer's share of the lower bound at every retailer lead time.

The same FIFO policy has a stronger guarantee when retailer holding rates $H_i$ are higher. Under $H_i\ge mh$, $m\ge1$, and $p\ge h$, its guarantee factor is $\max\{\gamma_m,1+K_{\max}/K_0\}$ with $\gamma_m=9/4+1/(4m)$; it is $\gamma_m$ if $K_{\max}\le(\gamma_m-1)K_0$, which includes supplier setup-cost dominance $K_i\le K_0$. In particular, $H_i\ge5h/3$ gives $12/5$ under the same demand and lead-time scope as the $5/2$ result; the price of the smaller factor is a slightly tighter setup-cost condition, $K_{\max}\le7K_0/5$ instead of $K_{\max}\le3K_0/2$.

For arbitrary nonnegative setup costs and associated retailer demands, a deterministic consolidation policy tuned to $m$ achieves a guarantee factor of $13/3+32/(57m)$ under $H_i\ge mh$, $m\ge1$, and $p\ge h$. At $m=1$, this is $93/19<4.9$. The warehouse retains each request's allocated stock until a dispatch threshold is reached, so each request ships in at most three portions, trading warehouse holding and slower retailer service for fewer dispatches. With $m=\min_iH_i/h$, the consolidation and FIFO factors depend only on cost parameters, so we can use the policy with the smaller factor.

On 864 instances built from Walmart store sales, the simulated cost of an integer-rounded FIFO implementation averages $1.45$ times the best evaluated lower bound and reaches $1.75$, while our upper bound on its cost is below $1.94$ times that bound. The policy requires no simulation search, yet its estimated cost comes within $8.2$ and $12.2$ percent, on average over 216 instances, of echelon $(R,nQ)$ and integer-ratio benchmarks, both tuned by simulation for each instance. Its mean simulated-cost ratio to the lower bound changes only from $1.46$ to $1.44$ as the supplier lead time rises from $1$ to $14$ days. The tuned benchmarks' simulated costs average $1.33$ and $1.29$ times the same lower bound, so, up to simulation error, they average at most these multiples of the optimum. For a more direct comparison with optimized decisions, we also compute bounded-state dynamic-programming values on 144 two-store instances. The ratio of the FIFO policy's simulated cost to these values averages $1.22$, while that of the lower bound to these values averages $0.82$ and is at least $0.71$. The policy's ratio rises as store setup costs grow relative to supplier setup costs.

\section{Literature Review}\label{sec:literature}

\textbf{Stochastic distribution systems.} \citet{ClarkScarf1960} establish the serial decomposition, and \citet{EppenSchrage1981} and \citet{FedergruenZipkin1984} develop allocation relaxations for distribution systems, which let retailer inventories be rebalanced freely; our mirror cost prices the imbalance that retailer batching creates between warehouse net stock and a pooled position. With fixed ordering costs, \citet{ChenZheng1997} and \citet{Marklund2011} analyze echelon batch-ordering and time-based shipment-consolidation policies. \citet{ChenSamroengraja2000} stagger identical retailers' orders so that one retailer orders in each period, a scheduling remedy for the imbalance that our mirror cost prices. \citet{ZhuChenHuYang2021} derive modified echelon $(r,Q)$ policies from independent single-stage problems for continuous-review OWMR systems with Poisson demand, with bounds that become asymptotically tight under specified cost scalings, as do the policies of \citet{HuangEtAl2025} for perishable goods as the number of retailers grows. The upper and lower bounds of \citet{GallegoEtAl2007} for base-stock systems differ by a warehouse newsvendor cost on pooled lead-time demand, whereas our pooled supplier cost enters both of our bounds. For distribution systems with expediting and supply disruptions, \citet{ZhaoBirgeDeValveInman2026} derive a stochastic-program lower bound on the optimal cost and set base-stock levels from that bound.

\indent\textbf{Lower bounds and guarantees.} Our lower bound builds on the cost-allocation framework of \citet{ChenZheng1994LB}, whose induced-penalty and parameter-allocation bounds assign part of each cost to retailer problems and the rest to a supplier problem; in joint replenishment, \citet{AtkinsIyogun1988} allocate the common setup cost among single-item problems to evaluate periodic and can-order policies. We allocate setup, holding, and backlog costs through bounded functions that satisfy single-location cost inequalities, so the lower bound contains the same pooled supplier problem as our upper bound. Since supplier purchases and retailer dispatches can occur in different periods, such an allocation needs a physical shipment rule: FIFO charges first portions to retailer orders and continuations to supplier receipts, which extends the upstream charge for extra downstream shipments in \citet{HuYang2014}. For serial systems, \citet{Chen1999}, \citet{HuYang2014}, and \citet{Hu2026} give guarantees under their own demand and cost assumptions, and \citet{LeviRoundyTruongWang2017} give factors between two and three against all nonanticipative policies without fixed costs. Asymptotic analyses instead show that simple policies become optimal in a limit, as the regular lead time grows in dual sourcing \citep{XinGoldberg2018} or as the product lifetime grows for perishable inventory \citep{BuGongChao2023}, whereas our factor holds at every lead time.

\indent\textbf{Approximation algorithms.} Constant-factor guarantees are known for deterministic OWMR and joint replenishment problems \citep{Roundy1985,LeviEtAl2008}, for single-stage stochastic inventory control through dual balancing \citep{LeviEtAl2007} and with capacities and setup costs \citep{ShiZhangChaoLevi2014}, and for base-stock policies on stochastic inventory networks without fixed costs or lead times \citep{DeValveMyles2025}. For a stochastic OWMR safety-stock approximation with independent normal demands, \citet{NajyDiabatElbassioni2025} develop an FPTAS and obtain power-of-two policies with a factor $1.021(1+\epsilon)$ when the base period is optimized. Their service-level formula ignores allocation issues, and the guarantee concerns that approximate optimization model. Our comparison includes physical allocation and dispatch costs and uses the optimum over all nonanticipative policies of the stochastic system. In e-commerce fulfillment, \citet{Ma2023} splits multi-item orders across fulfillment centers within a factor logarithmic in the number of items per order. Our analysis of the single-location problems rests on the $(s,S)$ optimality results of \citet{Iglehart1963}, \citet{Zheng1991}, and \citet{FeinbergLewis2018}, the last without requiring demand to be either discrete or continuous.

\indent\textbf{Synchronized replenishment.} With constant deterministic demand and no shortages, power-of-two policies with an optimized base period coordinate replenishment through nested cycles and come within about 2 percent of optimal cost \citep{Roundy1985}. With stochastic demand, our retailers order in response to their own demand, while the warehouse cannot predict demand exactly and, because of the supplier lead time, must position stock in advance. Retailer orders then need not match warehouse receipts, so stock waits for orders or orders wait for stock. Hence, coordination must handle uncertainty: the warehouse positions stock like a newsvendor against pooled retailer demand over the supplier lead time, at its holding rate, and dispatches it as retailer orders arrive, while the target of each retailer covers only the demand over its own lead time.

\section{Model and Main Results}\label{sec:model}

\begin{figure}[htbp]
\SingleSpaced
\centering
\includegraphics[width=.99\linewidth]{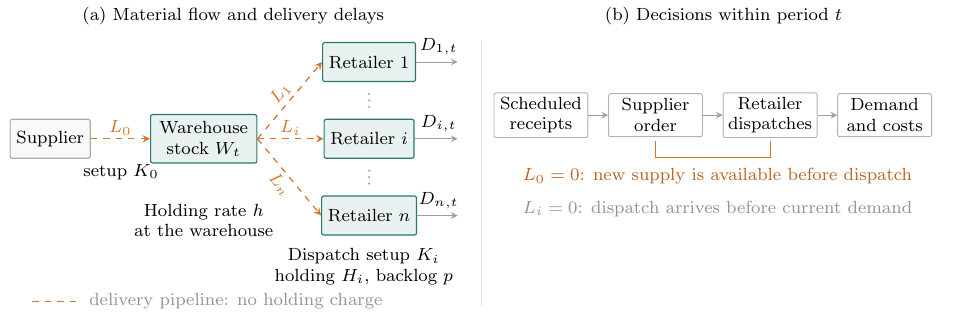}
\caption{The distribution system, its delivery delays, and the sequence of events within a period.}
\label{fig:system}
\vspace{12pt}
\end{figure}

Figure~\ref{fig:system} shows a warehouse that replenishes retailers $i=1,\ldots,n$ from an outside supplier in discrete periods. Each period $t$ has four steps. First, shipments sent in earlier periods arrive as scheduled. Second, the warehouse may place a supplier order. Third, the warehouse may dispatch stock to the retailers. Fourth, demand occurs at the retailers and is met from retailer stock or backlogged; holding and backlog costs are then charged on the resulting inventories. A supplier order placed in period $t$ arrives in period $t+L_0$, in time for that period's dispatches, and a dispatch to retailer $i$ in period $t$ arrives in period $t+L_i$, in time for that period's demand. The lead times $L_0$ and $L_i$ are finite nonnegative integers; with a zero lead time, a shipment can be used in the period it is sent. Quantities are divisible, capacities are unlimited, and unmet demand is backlogged. As in the classical literature, stock moves only downstream and cannot be disposed of.

The demand vectors $(D_{1,t},\ldots,D_{n,t})$ are independent and identically distributed (i.i.d.) over periods, with nonnegative components and $0<\E D_{i,t}<\infty$. Retailers' demands in a period may be dependent, as when stores share weather or promotions. Theorem~\ref{thm:arbitrary} assumes within-period association (Section~\ref{sec:arbitrary}). Let $D_i$ have retailer $i$'s marginal distribution and $D$ that of $D_t:=\sum_iD_{i,t}$.

Setup costs are $K_0,K_i\ge0$ per supplier order and per retailer dispatch. The warehouse holding cost is $h>0$ per unit per period, the retailer holding costs are $H_i>0$, and the common backlog penalty is $p>0$. We use a common penalty to represent a chain-wide charge for delayed service. The nearly uniform prices documented for U.S. chains by \citet{DellaVignaGentzkow2019} motivate treating the product's sales value similarly across stores; a common delay penalty additionally assumes that a delay costs comparable customer goodwill at every store. Holding costs include local rent, labor, and handling and may differ by store. Section~\ref{sec:conclusion} discusses extending the analysis to retailer-specific penalties. Units in transit incur no holding cost, which loses no generality.\footnote{Every unit spends exactly $L_0$ periods in the supplier pipeline and $L_i$ periods in the pipeline to retailer $i$, and a policy with finite average cost orders and dispatches at the demand rates in the long run, so all such policies have the same long-run average pipeline stock, $L_0\E D+\sum_iL_i\E D_i$. A linear charge on this stock adds the same constant to every policy's cost and to the optimal cost, which preserves every guarantee stated below.}

Let $q_{0,t}$ and $q_{i,t}$ be the supplier order and retailer dispatch quantities. Let $W_t\ge0$ be warehouse stock after dispatches and $I_{i,t}$ retailer net inventory after demand. The inventory balances are $W_t=W_{t-1}+q_{0,t-L_0}-\sum_iq_{i,t}$ and $I_{i,t}=I_{i,t-1}+q_{i,t-L_i}-D_{i,t}$.
Quantities indexed by negative periods describe a prescribed finite initial state with nonnegative warehouse stock, finite retailer inventories, and finite pipelines, which we call a valid state. All units dispatched to the same retailer in a period form one dispatch with one setup charge. With $x^-=(-x)^+$, the period cost is
\begin{equation}\label{eq:cost}
 c_t:=K_0\ind_{\{q_{0,t}>0\}}+\sum_iK_i\ind_{\{q_{i,t}>0\}}
       +hW_t+\sum_i\bigl(H_iI_{i,t}^++pI_{i,t}^-\bigr).
\end{equation}
A nonanticipative policy may use the observed history and auxiliary randomness independent of all demands. It is feasible if its order and dispatch quantities are nonnegative and its dispatches never exceed the warehouse stock available in that period. Define
\begin{equation}\label{eq:objective}
 C(\pi):=\limsup_{T\to\infty}\frac1T\E_\pi\sum_{t=0}^{T-1}c_t,
 \qquad \OPT:=\inf_{\pi\text{ nonanticipative and feasible}}C(\pi).
\end{equation}

We now define the deterministic FIFO policy; Section~\ref{sec:arbitrary} constructs a second policy for arbitrary setup costs. The FIFO policy serves retailer requests from warehouse stock in arrival order, and each retailer follows an inventory rule derived from a single-location problem. For $k\ge0$ and an integer $L\ge0$, the single-location problem of retailer $i$ with setup cost $k$ and lead time $L$ treats its supply as unlimited: each order costs $k$ and arrives after $L$ periods, and the retailer pays holding cost $H_i$ and backlog penalty $p$ on its net inventory after demand, as in~\eqref{eq:cost}. An $(s,S)$ policy raises the inventory position to $S$ whenever it falls to $s$ or below, and base stock is the case $s=S$; such a policy is optimal for every single-location problem \citep{Iglehart1963,Zheng1991,FeinbergLewis2018}. The supplier problem is the analogous single-location problem for warehouse replenishment against aggregate demand $D$: each supplier order costs $k$, and a position $y$ after ordering incurs $\E[h(y-D^{(L_0)})^++p(D^{(L_0)}-y)^+]$, where $D^{(L_0)}$ is aggregate demand over $L_0$ periods. Its optimal long-run average cost, including setup charges and inventory costs, is denoted by $J_0(k)$.

The FIFO policy also measures how each retailer's orders affect the warehouse. For a retailer policy whose inventory position after ordering in period $t$ is $\widehat P_{i,t}$, its mirror cost is the minimum over $a\in\R$ of $\lim_{T\to\infty}T^{-1}\E\sum_{t=0}^{T-1}[h{(a-\widehat P_{i,t})}^++p{(\widehat P_{i,t}-a)}^+]$, and a minimizing $a$ is its mirror center. As in the introduction, warehouse net stock mirrors the swings of the retailer's position about this center, and the mirror cost charges them at the warehouse rate. Base stock restores the same position every period and so has no swings; its mirror cost is zero even though demand is random.

For a retailer $i$ with $K_i>0$, let $(\hat s_i,\hat S_i)$ be an optimal $(s,S)$ policy for its single-location problem with planning setup cost $5K_i/3$ and lead time $L_i$, chosen as in Lemma~\ref{lem:scalar} to satisfy the canonical conditions~\eqref{eq:canonical}. Let $(s_i',S_i')$ have the smallest long-run average holding and backlog cost among $(s,S)$ policies with $\hat s_i\le s<S\le\hat S_i$ and $S-s=(\hat S_i-\hat s_i)/2$; all these policies order equally often. A retailer policy's augmented cost is its long-run average cost in the single-location problem with actual setup cost $K_i$ and lead time $L_i$, together with the mirror cost of its retailer position.

\begin{definition}[{\mdseries\scshape FIFO Policy}]\label{def:fifo}
Let $K_{\max}:=\max_iK_i$. After the initialization in Appendix~\ref{app:initialization}, the deterministic FIFO policy $\pi_{\rm FIFO}$ follows three rules and breaks ties by fixed rules.

(i)~\emph{Retailer policies.} Retailer $i$ uses whichever of $(\hat s_i,\hat S_i)$ and $(s_i',S_i')$ has the smaller augmented cost if $K_i>0$, and optimal base stock for its single-location problem with lead time $L_i$ if $K_i=0$.

(ii)~\emph{Requests and dispatch.} After the period's scheduled arrivals at the warehouse and the retailers, retailer $i$ evaluates its rule at its current inventory position including outstanding requests, before adding any new request. Write $\widehat P_{i,t}=P_{i,t}+\Delta_{i,t}$ for this position after the period's dispatches and before demand is realized: $P_{i,t}$ is net inventory plus shipments already dispatched and in transit, and $\Delta_{i,t}$ is the quantity of requests awaiting dispatch. The position, including outstanding requests, changes only through retailer $i$'s own requests and demand. Each positive request joins one FIFO queue. After the supplier order and any resulting immediate receipt, the warehouse fills requests in queue order from available stock, the last one possibly in part, and sends all the portions allocated to retailer $i$ in that period as one dispatch, with a single setup charge $K_i$.

(iii)~\emph{Supplier orders.} The warehouse applies an $(s,S)$ policy, attaining $J_0(K_0+K_{\max})$ in~\eqref{eq:source-value}, to the aggregate position $Y_t^{\rm pre}=W+O+\sum_i(P_i-a_i)$, the warehouse's echelon inventory position less the sum of the retailers' mirror centers, measured after these arrivals and immediately before the supplier order. Here $W$ is warehouse stock, $O$ is supplier stock in transit, $P_i$ is retailer $i$'s current inventory position excluding unfilled requests, and $a_i$ is its fixed mirror center. Retailer requests and dispatches leave this aggregate unchanged, so its post-order value is $Y_t=Y_t^{\rm pre}+q_{0,t}$.\Halmos
\end{definition}

Each retailer fixes its policy at the start and orders as if supply were unlimited, so its orders do not depend on warehouse stock. Since a retailer has no downstream stock, its echelon inventory position is $\widehat P_{i,t}$, and the warehouse's is $W_t+O_t+\sum_iP_{i,t}$. Hence, rule~(iii) is an $(s,S)$ rule on the warehouse's echelon inventory position with both thresholds raised by $\sum_ia_i$, and $\pi_{\rm FIFO}$ is an echelon stock $(s,S)$ policy with FIFO service. Section~\ref{sec:numerics} compares it with the echelon stock $(s,S)$ policy that plans each stage at its actual setup cost and uses only full windows.

The physical fulfillment rule keeps the inventory positions including outstanding requests and the request queue as records, and only its dispatches move stock to retailers. The warehouse ships every unit in the period it is allocated; the consolidation policy of Section~\ref{sec:arbitrary} keeps the FIFO queue but delays dispatch. The extra $K_{\max}$ in the supplier setup cost pays for the one physical dispatch that may continue a partly shipped request after each supplier delivery, as Lemma~\ref{lem:continuation} shows.

\begin{theorem}[{\mdseries\scshape General FIFO Guarantee}]\label{thm:main}
Let $m\ge1$, $\gamma_m:=9/4+1/(4m)$, and $\gamma\ge\gamma_m$. Suppose $H_i\ge mh$ for every retailer, $p\ge h$, and $K_{\max}\le(\gamma-1)K_0$. From every valid initial state, the deterministic policy $\pi_{\rm FIFO}$ of Definition~\ref{def:fifo} satisfies $C(\pi_{\rm FIFO})\le\gamma\OPT$ at all finite lead times.
\end{theorem}

Taking $m=1$ and $\gamma=\gamma_1=5/2$ gives $C(\pi_{\rm FIFO})\le(5/2)\OPT$ when $H_i\ge h$, $p\ge h$, and $K_i\le3K_0/2$ for every $i$. This condition includes supplier setup-cost dominance $K_i\le K_0$ and enters the proof only via $K_0+K_{\max}\le5K_0/2$. With $K_0>0$, a larger $\gamma$ admits larger retailer setup costs.

Two candidates suffice. The first, $(\hat s_i,\hat S_i)$, pays actual setup cost $K_i$, below the planning setup cost $5K_i/3$, but its position swings over the full window planned at that cost and creates a mirror cost. The second, $(s_i',S_i')$, halves the width, narrowing these swings, but can order more often and may miss part of the first candidate's savings in holding and backlog cost. Its window is placed where these savings are largest, which need not be at the top of $[\hat s_i,\hat S_i]$. Proposition~\ref{prop:scalar} bounds the smaller augmented cost by $3/2$ times the retailer's optimal single-location cost at setup cost $5K_i/3$.

Three holding and penalty assumptions play separate roles. First, $H_i\ge h$ ensures that the width $\hat S_i-\hat s_i$ times $hp/(h+p)$ is at most the optimal single-location cost at setup cost $5K_i/3$, which bounds the mirror cost, and lets the lower bound pool retailer costs at the warehouse rate. Second, $p\ge h$ enters only Proposition~\ref{prop:scalar}, which compares the candidates' augmented costs. Third, a common penalty lets Propositions~\ref{prop:upper} and~\ref{prop:lower} price each shortage at rate $p$, whichever retailer's request waits.

Retailer holding rates $H_i\ge mh$ raise the multiplier of $hp/(h+p)$ to $mhp/(mh+p)$, which shrinks the allowed width, and with it, the bound on the mirror cost. Theorem~\ref{thm:main} thus lowers the factor to $\gamma_m$ for the same FIFO policy. At $m=5/3$ and $\gamma=\gamma_m=12/5$, it gives $C(\pi_{\rm FIFO})\le(12/5)\OPT$ if $H_i\ge5h/3$, $p\ge h$, and $K_i\le7K_0/5$ for every $i$. This includes supplier setup-cost dominance $K_i\le K_0$ and retains the model's general demand distributions and arbitrary finite lead times.

\begin{figure}[htbp]
\SingleSpaced
\centering
\includegraphics[width=.99\linewidth]{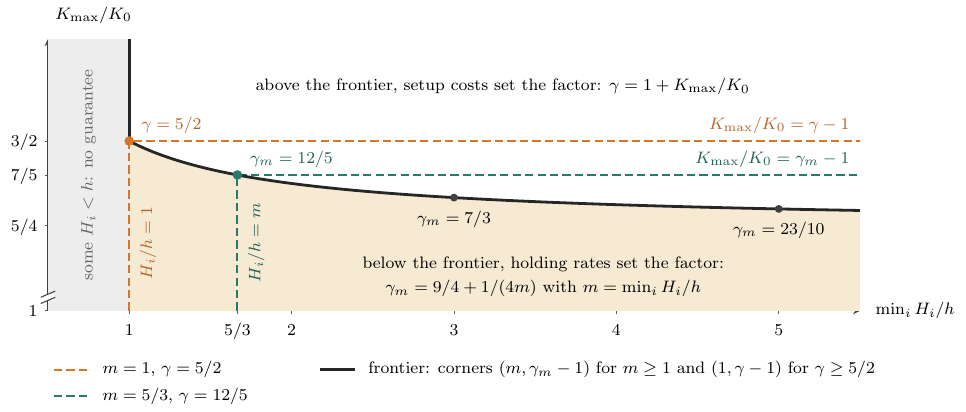}
\caption{Guarantee factors of Theorem~\ref{thm:main}.}
\label{fig:frontier}
\begin{minipage}{0.98\linewidth}
{\SingleSpaced\footnotesize\emph{Notes.} Here $K_0>0$ and $p\ge h$. Each dashed line bounds the region of a special case, $(m,\gamma)=(1,5/2)$ or $(5/3,12/5)$: to the right of its vertical part and below its horizontal part. The corners $(m,\gamma_m-1)$ of such regions for $m\ge1$, where $9/4<\gamma_m\le5/2$, trace the curved part of the frontier, and the corners $(1,\gamma-1)$ for $\gamma\ge5/2$ trace its vertical part. Dots mark the corners at $m=1$, $5/3$, $3$, and $5$, labeled with their factors. The vertical axis starts at $1$, and below $K_{\max}/K_0=5/4$ the factor depends only on $\min_iH_i/h$, because the frontier lies above $5/4$.}
\end{minipage}
\end{figure}

For $K_0>0$, set $m=\min_iH_i/h\ge1$. Choosing $\gamma=\max\{\gamma_m,1+K_{\max}/K_0\}$ in Theorem~\ref{thm:main} gives
\begin{equation}\label{eq:setup-ratio}
 C(\pi_{\rm FIFO})\le\max\!\left\{\frac94+\frac1{4m},\;1+\frac{K_{\max}}{K_0}\right\}\OPT.
\end{equation}
Figure~\ref{fig:frontier} shows this guarantee when $K_0>0$ and $p\ge h$. Below the frontier $K_{\max}/K_0=\gamma_m-1$, retailer holding rates set the factor $\gamma_m$; above it, setup costs set the factor $1+K_{\max}/K_0$. The guarantee approaches $9/4$ as retailer holding rates rise, provided $K_{\max}\le5K_0/4$. With associated retailer demands, the consolidation policy of Section~\ref{sec:arbitrary} gives the factor $\bar\gamma_m=13/3+32/(57m)<4.9$ for every setup-cost ratio, and selecting the policy with the smaller factor gives the smaller guarantee.

\section{The Single-Location Comparison}\label{sec:scalar}

This section proves Proposition~\ref{prop:scalar}: planning at setup cost $5K/3$ and selecting the cheaper of the full window and the best shifted half window gives augmented cost at most $3/2$ times the planning value. The proof uses three ingredients. First, Lemma~\ref{lem:reward} expresses the canonical cycle's savings in terms of its planning setup cost. Second, Lemma~\ref{lem:shift} finds a half window that retains at least half of those savings. Third, Lemma~\ref{lem:tail} links the savings in the upper half of the full window to the half-window mirror cost, allowing the two candidate costs to be compared.

\subsection{Single-Location Inventory Problems}

For nonnegative i.i.d.\ demand $D$ with $0<\E D<\infty$ and a nonnegative convex, coercive (${f(x)\to\infty}$ as ${|x|\to\infty}$), globally Lipschitz function $f$, consider the following single-location inventory problem with general cost $f$. In each period, the position may be raised; a positive increase costs $K$; the post-action position $x$ incurs $f(x)$; and demand then reduces it. Its single-location value $\J_D(f;K)$ is the infimum of the expected long-run average cost, defined by the limsup criterion in~\eqref{eq:objective}, over all nonanticipative policies; by Lemma~\ref{lem:scalar} and Appendix~\ref{app:scalar}, it does not depend on the initial position. Coercivity ensures an optimal policy with finite reorder and order-up-to levels $\hat s$ and $\hat S$.

The following lemma combines the classical optimality of $(s,S)$ policies \citep{Iglehart1963,Zheng1991,FeinbergLewis2018} with a bounded function satisfying a one-period cost inequality. Appendix~\ref{app:scalar} proves it for real positions and every nonnegative demand distribution with finite positive mean, without requiring demand to be discrete or to have a density.

\begin{lemma}[{\mdseries\scshape Single-Location Inventory Problem}]\label{lem:scalar}
Let $J=\J_D(f;K)$. If $K>0$, an optimal policy orders to $\hat S$ whenever its pre-action position is at most $\hat s$, where $\hat s<\hat S$ satisfy the canonical conditions
\begin{equation}\label{eq:canonical}
 f(\hat s)=J,\qquad f(\hat S)\le J,\qquad f(x)\le J\quad(\hat s\le x\le\hat S).
\end{equation}
The policy does not order on $(\hat s,\hat S]$. When initialized at $\hat S$, the policy has long-run average cost $J$, bounded post-action positions, and a cycle length with finite moments of every order. There is a bounded function $u$ with oscillation $\osc(u):=\sup u-\inf u\le K$ such that
\begin{equation}\label{eq:account}
 J+(f(x)-J)^+\le f(x)+\E[u(x-D)-u(x)]\qquad(x\in\R).
\end{equation}
For every nonanticipative policy whose information available before each period's demand is realized, including any observed network history and auxiliary randomness, is independent of that demand, the post-action positions $x_t$ and order quantities $q_t$ satisfy, from every initial position,
\begin{equation}\label{eq:finite-horizon}
 TJ-K\le\E\sum_{t=0}^{T-1}\bigl[K\ind_{\{q_t>0\}}+f(x_t)\bigr]\qquad(T\ge1).
\end{equation}
If $K=0$, base stock at a minimizer of $f$ gives $J=\min f$, and both~\eqref{eq:account} and~\eqref{eq:finite-horizon} hold with $u=0$.
\end{lemma}

The bounded function $u$ records costs transferred between periods, and telescoping these transfers gives the finite-horizon bound~\eqref{eq:finite-horizon}. Since $\osc(u)\le K$ caps the endpoint difference, the bound holds without a terminal inventory restriction. Appendix~\ref{app:scalar} proves it for policies that also use network history and auxiliary randomness, as Proposition~\ref{prop:lower} needs. Section~\ref{sec:lower} drops the nonnegative excess $(f(x)-J)^+$ in~\eqref{eq:account}; retaining it would strengthen the shared-cost lower bound.

We now connect this single-location problem to a retailer with lead time. Let $P_{i,t}$ be retailer $i$'s net inventory plus shipments already dispatched and in transit to it, measured immediately after dispatches and before demand is realized. Unfilled warehouse requests are excluded from $P_{i,t}$. The inventory balance equations imply
\begin{equation}\label{eq:retailer-lead}
 I_{i,t+L_i}=P_{i,t}-\sum_{\ell=t}^{t+L_i}D_{i,\ell},
\end{equation}
because later dispatches arrive after $t+L_i$. Let $Z_i$ be a sum of $L_i+1$ fresh demands and define
\begin{equation}\label{eq:retailer-value}
 F_i(x):=\E\bigl[H_i(x-Z_i)^++p(Z_i-x)^+\bigr],
 \qquad J_i(k):=\J_{D_i}(F_i;k).
\end{equation}
Given the history through the dispatches of period $t$, which fixes $P_{i,t}$ and is independent of later demands, the expected physical retailer loss at $t+L_i$ is $F_i(P_{i,t})$ by~\eqref{eq:retailer-lead}. The sum $Z_i$ has $L_i+1$ terms: a dispatch from period $t$ is in transit during $L_i$ periods of demand, and the demand of its arrival period $t+L_i$ also occurs before that period's cost is charged. Each $F_i$ satisfies the assumptions of Lemma~\ref{lem:scalar}.\footnote{The finite first moment of $Z_i$ makes $F_i$ finite and globally Lipschitz. Jensen's inequality for the convex loss $y\mapsto H_iy^++p(-y)^+$ gives $F_i(x)\ge H_i(x-\E Z_i)^++p(\E Z_i-x)^+$, which tends to infinity as $|x|\to\infty$, so $F_i$ is coercive.} The policy uses $F_i$ to choose each retailer's replenishment cycle, and the shared-cost lower bound of Section~\ref{sec:lower} uses the single-location values $J_i(k)$ for every retailer in the network.

\subsection{The Renewal Measure and Short Intervals}\label{sec:renewal}

To compare cycles of different widths, we express their lengths and costs through cumulative demand. For generic demand, write $D^{(0)}=0$ and $D^{(n)}=D_1+\cdots+D_n$. The renewal measure $U(A)$ is the expected number of partial sums in $A$, and $M(w)$ is the expected number of such sums below $w$:\footnote{Each $\Prob(D^{(n)}\in\cdot)$ is a probability measure on $[0,\infty)$. For disjoint sets $A_1,A_2,\ldots$, exchanging the order of summation of nonnegative terms gives $U(\bigcup_jA_j)=\sum_{n\ge0}\sum_j\Prob(D^{(n)}\in A_j)=\sum_jU(A_j)$, so $U$ is a measure. Let $\tau_{w}=\inf\{n\ge1:D^{(n)}\ge w\}$ be the number of periods needed for cumulative demand to reach or exceed $w$. Since demands are nonnegative, $D^{(n)}<w$ exactly when $n<\tau_{w}$, so $M(w)=U[0,w)=\sum_{n\ge0}\Prob(D^{(n)}<w)=\E\sum_{n\ge0}\ind_{\{D^{(n)}<w\}}=\E\sum_{n=0}^{\tau_{w}-1}1=\E\tau_{w}$. Choose $d>0$ with $\Prob(D\ge d)>0$. Cumulative demand reaches $w$ once $\lceil w/d\rceil$ demands of at least $d$ have occurred, and the wait for each has mean $1/\Prob(D\ge d)$, so $M(w)\le\lceil w/d\rceil/\Prob(D\ge d)<\infty$. Thus, $U$ is finite on bounded sets.}
\begin{equation}\label{eq:renewal}
 U(A):=\sum_{n\ge0}\Prob(D^{(n)}\in A),\qquad M(w):=U[0,w),\quad w>0.
\end{equation}
For every $w>0$, $M(w)<\infty$. A width-$w$ policy with target $\Gamma$ is the $(s,S)$ policy with $S=\Gamma$ and $s=\Gamma-w$; it orders back to $\Gamma$ when its depletion, the amount by which its position lies below the target, reaches $w$. Started at its target $\Gamma$, the policy has i.i.d.\ cycles. In the periods within a cycle, the post-action depletion equals $D^{(0)},\ldots,D^{(\tau_{w}-1)}$. A bounded measurable $v$ therefore has expected cycle cost $\sum_{n\ge0}\E[v(D^{(n)})\ind_{\{D^{(n)}<w\}}]=\int_{[0,w)}v\,dU$, and the expected cycle length is $M(w)$. The renewal-reward theorem for expected time averages \citep[Section~3.6]{Ross1996} gives order frequency $1/M(w)$ and long-run post-action position $\Gamma-X_{w}$, where the depletion $X_{w}$ is a random variable with distribution $U(\cdot\cap[0,w))/M(w)$, a probability measure because $U[0,w)=M(w)$, so that
\begin{equation}\label{eq:occupation}
 \E v(X_{w})=\frac1{M(w)}\int_{[0,w)}v(x)\,U(dx)
\end{equation}
for bounded measurable $v$. Taking $v(x)=F(\Gamma-x)$ in this renewal-reward identity gives the average holding and backlog cost. Here $X_{w}$ is distributed as the long-run fraction of periods spent at each depletion level, so~\eqref{eq:occupation} gives time averages, not limits: the depletion dynamics need not converge to a steady state; under deterministic demand, depletion cycles through fixed levels.

Two conventions keep~\eqref{eq:occupation} exact for every demand distribution, including discrete ones with point masses. First, intervals of depletion levels, such as $[0,w)$, are half open: a cumulative demand equal to the width triggers a new order and is excluded from that cycle. Second, repeated zero demands are retained, so $U\{0\}=1/\Prob(D>0)$. Thus, no point mass is split or smoothed. Lemma~\ref{lem:mirror} uses the mean bound $\E X_{w}\le w/2$ in the next lemma, and Lemma~\ref{lem:tail} and Proposition~\ref{prop:scalar} use its doubling bound $M(2w)\le2M(w)$: doubling the width at most doubles the expected cycle length.

\begin{lemma}[{\mdseries\scshape Renewal Interval Bounds}]\label{lem:renewal}
For $x\ge0$ and $w>0$, $U[x,x+w)\le M(w)$ and $M(2w)\le2M(w)$, and the long-run depletion $X_w$ of a width-$w$ policy satisfies $\E X_{w}\le w/2$.
\end{lemma}

\proof{Proof.} Fresh demands after the first passage $\tau_x$ give the interval-count bound $U[x,x+w)\le M(w)$. Setting $x=w$ gives the doubling bound, and integrating over $x$ gives $\E X_{w}\le w/2$.

\emph{Step 1. Bound the interval counts.} We first show $U[x,x+w)\le M(w)$: in expectation, no interval of length $w$ holds more partial sums than $[0,w)$. At $x=0$ equality holds. For $x>0$, partial sums before $\tau_x$ lie below $x$, so all partial sums in $[x,x+w)$ occur from $\tau_x$ on. If $D^{(\tau_x)}\ge x+w$, there are none. Otherwise, since $\tau_x$ is a stopping time, the subsequent demands are i.i.d.\ copies of $D$, independent of $D_1,\ldots,D_{\tau_x}$. The partial sums from $\tau_x$ onward, therefore, equal $D^{(\tau_x)}$ plus the cumulative sums of these fresh demands, starting with zero. They lie below $x+w$ exactly when the corresponding fresh cumulative sum lies below $x+w-D^{(\tau_x)}$. Given $D^{(\tau_x)}$, their expected number is, therefore, $M(x+w-D^{(\tau_x)})\le M(w)$, because $D^{(\tau_x)}\ge x$ and $M$ is nondecreasing. Taking expectations gives $U[x,x+w)\le M(w)$, and $x=w$ gives $M(2w)=U[0,w)+U[w,2w)\le2M(w)$.

\emph{Step 2. Integrate to bound mean depletion.} For the last inequality, $U[x,w)\le M(w-x)$ for ${0<x<w}$. Integrating $\ind_{\{x\le y\}}$ against $dx\,U(dy)$ over $y$ first gives $U[x,w)$, and over $x$ first gives $y$. Since Tonelli's theorem allows either order, $\int_0^{w}U[x,w)\,dx=\int_{[0,w)}y\,U(dy)=M(w)\E X_{w}$, where the last equality is~\eqref{eq:occupation} with $v(y)=y$. Integrating $U[x,w)\le M(w-x)$ over $(0,w)$, using this identity and replacing $x$ by $w-x$, gives $M(w)\E X_{w}\le\int_0^{w}M(x)\,dx$. The partition $[0,w)=[0,x)\cup[x,w)$ gives $M(x)=M(w)-U[x,w)$, so the identity also gives $\int_0^{w}M(x)\,dx=wM(w)-M(w)\E X_{w}$. Hence, $2M(w)\E X_{w}\le wM(w)$; dividing by $2M(w)>0$ proves the mean bound $\E X_{w}\le w/2$.
\Halmos\endproof

\subsection{The Canonical Savings and the Mirror Cost}

Fix one retailer and suppress its index in~\eqref{eq:retailer-value}. Let $K$ be its actual setup cost, $H$ its holding rate, and $F$ its lead-time-demand cost. The single-location optimum at the actual setup cost is $J(K)=\J_D(F;K)$. For a planning setup cost $k>0$, define $G:=J(k)=\J_D(F;k)$ as the optimal value of the planning problem. We choose the replenishment cycle using this planning problem and evaluate it at the actual setup cost $K$ per order. Proposition~\ref{prop:scalar} will set $k=5K/3$ for $K>0$, so the planning problem charges each order two-thirds more than the actual setup cost $K$.

Apply Lemma~\ref{lem:scalar} at setup cost $k$ to choose an optimal pair $(\hat s,\hat S)$ for the planning problem. Its canonical conditions~\eqref{eq:canonical} use $G$ as the optimal value, so $F(\hat s)=G$. We call this pair the canonical optimum and its width $\hat w=\hat S-\hat s>0$ the canonical width. At depletion $x$ below the target, the holding and backlog cost $F(\hat S-x)$ falls short of the planning value $G$ by the cost savings
\begin{equation}\label{eq:reward}
 g(x):=G-F(\hat S-x)\quad(0\le x\le \hat w),
 \qquad \beta:=\frac{hp}{h+p}.
\end{equation}

\begin{lemma}[{\mdseries\scshape Canonical Savings and Width}]\label{lem:reward}
The function $g$ is continuous, concave, and satisfies
\begin{equation}\label{eq:reward-properties}
 0\le g\le G,\qquad g(\hat w)=0,\qquad
 \int_{[0,\hat w)}g(x)\,U(dx)=k,\qquad g(x)\le p(\hat w-x).
\end{equation}
The planning value $G$ is positive. If $H\ge h$, the canonical width satisfies $\beta \hat w\le G$.
\end{lemma}

\proof{Proof.} Canonical optimality and the renewal-reward identity~\eqref{eq:occupation} determine the shape and total amount of the savings $g$. Jensen's inequality then bounds the canonical cycle width $\hat w$.

\emph{Step 1. Establish the savings properties.} As $F$ is continuous and convex by~\eqref{eq:retailer-value}, $g$ is continuous and concave. By~\eqref{eq:canonical}, $F\le G$ on $[\hat s,\hat S]$ and $F(\hat s)=G$, so $g\ge0$ and $g(\hat w)=0$, while $F\ge0$ gives $g\le G$. By Lemma~\ref{lem:scalar}, the canonical policy, which is the width-$\hat w$ policy with target $\hat S$, has average cost $G$ when started at $\hat S$, so its order frequency $1/M(\hat w)$ and~\eqref{eq:occupation} give
\begin{equation}\label{eq:canonical-average-cost}
 G=\frac{k+\int_{[0,\hat w)}F(\hat S-x)\,U(dx)}{M(\hat w)}.
\end{equation}
The definition of $g$ and $U[0,\hat w)=M(\hat w)$ give $\int_{[0,\hat w)}g(x)\,U(dx)=GM(\hat w)-\int_{[0,\hat w)}F(\hat S-x)\,U(dx)$, which equals $k$ by~\eqref{eq:canonical-average-cost}. Equation~\eqref{eq:canonical-average-cost} also gives $G>0$, because $k>0$ and $F\ge0$. Since $F$ has slopes in $[-p,H]$ and $F(\hat s)=G$, $F(\hat S-x)\ge G-p(\hat w-x)$, so $g(x)\le p(\hat w-x)$ on $[0,\hat w]$.

\emph{Step 2. Bound the canonical width.} Since $H(y-z)^++p(z-y)^+$ is convex in $z$, Jensen's inequality applied to~\eqref{eq:retailer-value} gives $F(y)\ge H(y-\E Z)^++p(\E Z-y)^+$. At the two canonical endpoints $\hat s$ and $\hat S$, this implies $\hat s\ge\E Z-G/p$ and $\hat S\le\E Z+G/H$, so $\hat w=\hat S-\hat s\le(H+p)G/(Hp)$. Multiplying this width bound by $\beta=hp/(h+p)$ gives $\beta\hat w\le h(H+p)G/[H(h+p)]\le G$, where the final comparison follows from $H\ge h$, since this gives $h(H+p)\le H(h+p)$, that is, $hp\le Hp$.
\Halmos\endproof

We next account for the warehouse imbalance created by a retailer cycle through its mirror cost. The warehouse-rate loss is $\phi(x):=hx^++px^-$, with $x^+=\max\{x,0\}$ and $x^-=\max\{-x,0\}$: it charges the warehouse holding rate $h$ on a surplus and the backlog penalty $p$ on a shortfall. For width $w$, define its mirror cost by
\begin{equation}\label{eq:mirror}
 B_{w}:=\min_{z\in\R}\E\bigl[p(z-X_{w})^++h(X_{w}-z)^+\bigr].
\end{equation}
The mirror cost is the smallest long-run average warehouse-rate loss $\phi(a-\widehat P)$ over fixed centers $a=\Gamma-z$, where $\widehat P=\Gamma-X_{w}$ is the post-action position including outstanding requests. A minimizing depletion center exists in $[0,w]$.\footnote{Because $0\le X_{w}<w$, the expectation in~\eqref{eq:mirror} equals $h(\E X_{w}-z)$ for $z<0$ and $p(z-\E X_{w})$ for $z>w$, so it decreases up to $z=0$ and increases beyond $z=w$. Since it is continuous in $z$, it attains its minimum in $[0,w]$.} By the newsvendor condition, a minimizing $z$ is an $h/(h+p)$ quantile of $X_{w}$, so $a=\Gamma-z$ is a $p/(h+p)$ quantile of $\widehat P$. The per-period loss $\phi(a-\widehat P)$ is the mirror loss, and its optimized long-run average is the mirror cost. The mirror cost, therefore, charges the penalty rate $p$ per unit by which the retailer's position exceeds $a$, and the warehouse holding rate $h$ per unit by which it falls below $a$, as Figure~\ref{fig:mirror} in Appendix~\ref{app:scalar} shows. Adding this mirror cost to retailer setup, holding, and backlog cost gives the augmented cost in Lemma~\ref{lem:mirror} and Proposition~\ref{prop:scalar}.

\begin{lemma}[{\mdseries\scshape Mirror Bounds and Candidate Cost}]\label{lem:mirror}
For every $w>0$, $B_{w}\le\beta w$ and $B_{w}\le h\E X_{w}\le hw/2$. For $0<w\le \hat w$ and $0\le \upsilon\le \hat w-w$, the width-$w$ policy with target $\hat S-\upsilon$, started at its target and charged at actual setup cost $K$, has augmented average cost
\begin{equation}\label{eq:candidate}
 \cL(w,\upsilon)=G+\frac{K-\int_{[0,w)}g(\upsilon+x)\,U(dx)}{M(w)}+B_{w}.
\end{equation}
\end{lemma}

\proof{Proof.} At the center $z=hw/(h+p)$, the mirror loss is convex in $x\in[0,w]$ and equals $\beta w$ at both endpoints, so $B_w\le\beta w$. At $z=0$, it equals $hx$; hence $B_w\le h\E X_w\le hw/2$ by Lemma~\ref{lem:renewal}.

For the candidate with target $\hat S-\upsilon$, the renewal-reward identity~\eqref{eq:occupation} gives
\[
 \cL(w,\upsilon)=\frac{K}{M(w)}+\frac{1}{M(w)}\int_{[0,w)}F(\hat S-\upsilon-x)\,U(dx)+B_w=G+\frac{K-\int_{[0,w)}g(\upsilon+x)\,U(dx)}{M(w)}+B_w.
\]
The second equality uses $F(\hat S-\upsilon-x)=G-g(\upsilon+x)$ and $U[0,w)=M(w)$. The mirror term depends only on the width: for target $\Gamma=\hat S-\upsilon$ and center $a=\Gamma-z$, the loss is $\phi(a-(\Gamma-X_w))=p(z-X_w)^++h(X_w-z)^+$, whose expectation, minimized over $z$, is $B_w$.
\Halmos\endproof

\subsection{The Half-Window Shift Lemma}\label{sec:shift}

At a fixed width, the candidate-cost formula~\eqref{eq:candidate} shows that choosing the target amounts to maximizing retained savings. The next lemma shows that some shifted half window retains at least half of the full-cycle savings, which Proposition~\ref{prop:scalar} uses to bound the half-window cost.

\begin{lemma}[{\mdseries\scshape Half-Window Shift}]\label{lem:shift}
Let $b>0$ and let $g$ be continuous, nonnegative, and concave on $[0,2b]$. Define $\Phi(\upsilon):=\int_{[0,b)}g(\upsilon+x)\,U(dx)$ for $0\le \upsilon\le b$. Then
\begin{equation}\label{eq:shift}
 \max_{0\le \upsilon\le b}\int_{[0,b)}g(\upsilon+x)\,U(dx)
 \ge\frac12\int_{[0,2b)}g(x)\,U(dx).
\end{equation}
A deterministic shift $\upsilon^*\in[0,b]$ attains the maximum on the left-hand side of~\eqref{eq:shift}.
\end{lemma}

\begin{figure}[htbp]
\SingleSpaced
\centering
\includegraphics[width=.99\linewidth]{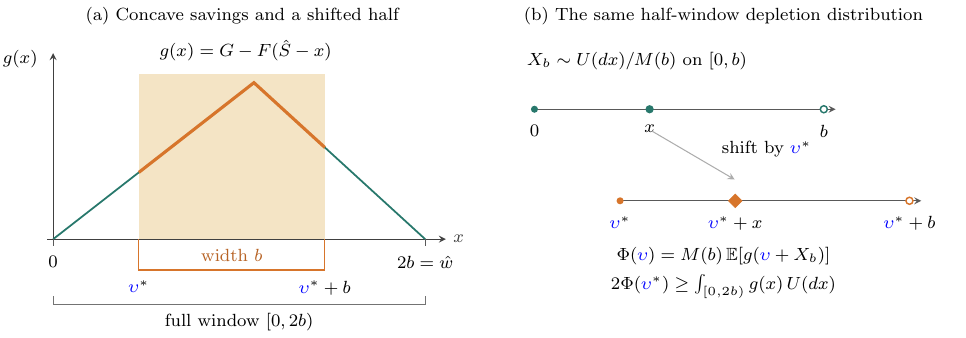}
\caption{The canonical savings window, its best shifted half, and the depletion distribution that the shift preserves.}
\label{fig:shift}
\begin{minipage}{0.98\linewidth}
\SingleSpaced\footnotesize\emph{Notes.} Shading marks a window of width $b$. Panel (b) shifts $x$ to $\upsilon^*+x$, preserving the distribution of $X_b$. Open right endpoints mark half-open windows. The curve is illustrative; the best shift depends on renewal measure $U$.
\end{minipage}
\end{figure}

Halving the width replaces $X_{2b}$ by $X_b$, weakly increasing the order frequency, whereas shifting the target by $\upsilon$ replaces $g(X_b)$ by $g(\upsilon+X_b)$ without changing the depletion distribution or the order frequency. The best shift $\upsilon^*$ depends only on $g$ and $U$, so it is fixed in advance. The proof splits a full cycle at the first period when depletion reaches $b$: the first part is the unshifted half window, and concavity yields a shifted half window with at least half the full-cycle savings (Figure~\ref{fig:shift}).

\proof{Proof.} We split the full cycle at the first passage of its half-width. Concavity then bounds the combined savings of the two pieces by twice the best shifted half-window savings.

\emph{Step 1. Decompose the full-cycle savings.} Set $\Phi^* :=\max_{0\le \upsilon\le b}\Phi(\upsilon)$. Uniform continuity of $g$ and finiteness of $M(b)$ make $\Phi$ continuous, so its maximum exists. Let $\tau_b=\inf\{n\ge1:D^{(n)}\ge b\}$ be the first period in which cumulative demand reaches $b$, and let $Y=D^{(\tau_b)}\ge b$ be the cumulative demand at that time, which is finite almost surely. Define $\Psi(y):=\int_{[0,2b-y)}g(y+x)\,U(dx)$ for $b\le y<2b$ and $\Psi(y)=0$ for $y\ge2b$. Split the partial sums counted by $U$ on $[0,2b)$ at $\tau_b$. Those before $\tau_b$ lie in $[0,b)$, so their expected contribution is $\Phi(0)$. Because $\tau_b$ is a stopping time, the later demands $D_{\tau_b+1},D_{\tau_b+2},\ldots$ are independent of $D_1,\ldots,D_{\tau_b}$ and i.i.d.\ with the same distribution as $D$, so the partial sums from $\tau_b$ onward equal $Y$ plus the cumulative sums of these fresh demands, starting with zero.\footnote{An exact hit at $b$ falls in the part from $\tau_b$ onward, and an exact hit at $2b$ is excluded; the fresh sums start at $Y$ itself and count repeated zero demands, in the same way that the renewal measure $U$ counts partial sums.} Such a sum $Y+x$ lies below $2b$ exactly when the fresh cumulative sum $x$ lies below $2b-Y$. Since the fresh sums are independent of $Y$ and $U$ gives their expected counts, the expected contribution given $Y<2b$ is $\int_{[0,2b-Y)}g(Y+x)\,U(dx)=\Psi(Y)$; if $Y\ge2b$, there is none. This gives the decomposition $\int_{[0,2b)}g(x)\,U(dx)=\Phi(0)+\E\Psi(Y)$ of the full-cycle savings.

\emph{Step 2. Compare the two pieces using concavity.} It remains to show that $\Phi(0)+\Psi(y)\le2\Phi^*$ for every value $y\ge b$ of $Y$; taking expectations in the decomposition proves~\eqref{eq:shift}. Let $m$ be a maximizer of $g$ on $[0,2b]$, which exists by continuity. Cases~1 and~2 treat $y\ge2b$ and $b\le y<2b$.

\emph{Case 1: $y\ge2b$.} Then $\Psi(y)=0$, so $\Phi(0)+\Psi(y)=\Phi(0)\le\Phi^*\le2\Phi^*$ because $\Phi^*\ge0$.

\emph{Case 2(a): $b\le y<2b$ and $m\le b$.} Then $g$ is nonincreasing on $[b,2b]$, since $g$ is concave and $m\le b$, so $g(y+x)\le g(b+x)$ for $0\le x<2b-y$. Also, $[0,2b-y)\subseteq[0,b)$ and $g\ge0$, so extending the integral gives $\Psi(y)\le\Phi(b)$. Hence, $\Phi(0)+\Psi(y)\le\Phi(0)+\Phi(b)\le2\Phi^*$, since $\Phi\le\Phi^*$ on $[0,b]$.

\emph{Case 2(b): $b\le y<2b$, $m>b$, and $y\le m$.} For $0\le x<2b-y$, where both $\Phi(0)$ and $\Psi(y)$ integrate, the concavity of $g$ on $[x,y+x]$ gives $g(y/2+x)\ge[g(x)+g(y+x)]/2$. For the remaining $2b-y\le x<b$, where only $\Phi(0)$ integrates, write $z=y/2+x$. If $z\le m$, then $g$ is nondecreasing up to $m$, so $g(z)\ge g(x)$. If $z>m$, then $m<z<b+y/2\le b+m/2$, because $x<b$ and $y\le m$. Hence, $m<2b$ and $2b-z\ge(2b-m)/2$, and the concavity of $g$ on $[m,2b]$ with $g(2b)\ge0$ gives $g(z)\ge g(m){(2b-z)}/{(2b-m)}\ge g(x)/2$, where the last step also uses $g(m)\ge g(x)$. Integrating these inequalities over the two ranges of $x$ yields $\Phi(y/2)\ge[\Phi(0)+\Psi(y)]/2$. Since $0\le y/2<b$, the shift $y/2$ is admissible in $\Phi^*$, so $\Phi(0)+\Psi(y)\le2\Phi(y/2)\le2\Phi^*$, which proves the claim of Step~2 in this case.

\emph{Case 2(c): $b\le y<2b$, $m>b$, and $y>m$.} Then $g$ is nonincreasing on $[m,2b]$, by the concavity of $g$, and nonnegative on $[0,2b]$, so $g(y+x)\le g(m+x)$ and $[0,2b-y)\subseteq[0,2b-m)$ give $\Psi(y)\le\Psi(m)$. The value $m$ itself satisfies the conditions of Case~2(b), because $b<m<y<2b$, so applying that case with $y$ replaced by $m$ gives $\Phi(0)+\Psi(m)\le2\Phi(m/2)\le2\Phi^*$. Combining this bound with $\Psi(y)\le\Psi(m)$ from the present case gives $\Phi(0)+\Psi(y)\le\Phi(0)+\Psi(m)\le2\Phi^*$.

The two cases cover every possible value of the cumulative demand $Y$. Taking expectations of their pointwise bound and using the decomposition of Step~1 proves the shift bound~\eqref{eq:shift}.
\Halmos\endproof

\subsection{The Tail-Mirror Inequality}\label{sec:tail}

The next estimate couples the savings in the upper half of a window with the mirror cost of its lower half. The quantities in the statement depend on the same renewal measure, and the proof uses this dependence to improve the coefficient obtained from separate bounds on the two terms.

\begin{lemma}[{\mdseries\scshape Tail-Mirror Inequality}]\label{lem:tail}
Let $h,p,b>0$, $r=p/h$, and $\beta=hp/(h+p)$. Set $\mathcal T_b:=\int_{[b,2b)}(2b-x)\,U(dx)$. Then
\begin{equation}\label{eq:tail-mirror}
 p\mathcal T_b+4rM(b)B_b\le(3r+1)\beta bM(b).
\end{equation}
\end{lemma}

The quantity $\mathcal T_b$ is the expected sum, over the periods of a cycle whose depletion lies in the upper half $[b,2b)$ of a width-$2b$ window, of the distance to the right endpoint $2b$. When $\hat w=2b$, Lemma~\ref{lem:reward} gives $g(x)\le p(\hat w-x)$, so $p\mathcal T_b$ bounds the savings $\int_{[b,2b)}g(x)\,U(dx)$ collected there.

To see why a joint bound is needed, compare separate estimates. The renewal interval bound alone gives $\mathcal T_b\le bM(b)$, while the endpoint mirror bound gives $B_b\le\beta b$. These separate bounds yield coefficient $r+1+4r=5r+1$ in place of $3r+1$ in~\eqref{eq:tail-mirror}, and at $r=1$ and maximal normalized width $\beta\hat w/G=1$, both candidate bounds of Proposition~\ref{prop:scalar} would exceed $3/2$ for every feasible normalized setup cost $\kappa=k/[GM(b)]\le2$. Thus, the retailer factor $3/2$ relies on this joint bound.

\proof{Proof.} We express the tail through integrated renewal counts, then compare it with the losses at two fixed mirror centers. Solving for these counts and substituting yields~\eqref{eq:tail-mirror}.

\emph{Step 1. Apply the doubling bound.} Tonelli's theorem gives
\[
 Q(y):=\int_0^yM(x)\,dx=\int_{[0,y)}(y-x)\,U(dx).
\]
Splitting the last integral at $b$ and using Lemma~\ref{lem:renewal}, respectively, yield
\[
 \mathcal T_b=Q(2b)-Q(b)-bM(b),\quad Q(2b)-Q(b)=2\int_{b/2}^bM(2v)\,dv\le4\int_{b/2}^bM(v)\,dv=4\{Q(b)-Q(b/2)\}.
\]
Consequently, $\mathcal T_b\le4\{Q(b)-Q(b/2)\}-bM(b)$, which involves $U$ only on $[0,b)$.

\emph{Step 2. Evaluate two mirror centers.} Define the cycle loss
\[
 V_z:=\int_{[0,b)}[p(z-x)^++h(x-z)^+]\,U(dx)\ge M(b)B_b.
\]
At $z=0$ and $z=b/2$, the identities $(x-z)^+=x-z+(z-x)^+$ and $\int_{[0,z)}(z-x)\,U(dx)=Q(z)$ give
\[
 V_0=h\{bM(b)-Q(b)\},\qquad V_{b/2}=\frac{hbM(b)}2-hQ(b)+(h+p)Q(b/2).
\]
These identities retain all point masses of the renewal measure, including any mass at zero.

\emph{Step 3. Combine the bounds.} Solving the preceding equations gives
\[
 Q(b)-Q(b/2)
 =\frac{(h+2p)bM(b)}{2(h+p)}
  -\frac{pV_0}{h(h+p)}-\frac{V_{b/2}}{h+p}.
\]
Substitution into the Step~1 bound, followed by $V_0,V_{b/2}\ge M(b)B_b$, yields
\[
 p\mathcal T_b\le\frac{p(h+3p)}{h+p}bM(b)-\frac{4p^2}{h(h+p)}V_0-\frac{4p}{h+p}V_{b/2}\le(3r+1)\beta bM(b)-4rM(b)B_b.
\]
The last inequality uses $r=p/h$, $\beta=hp/(h+p)$, and the fact that the two positive coefficients of the fixed-center losses sum to $4p/h=4r$. Rearranging this bound proves~\eqref{eq:tail-mirror}.
\Halmos\endproof

\subsection{The Two-Candidate Comparison}\label{sec:comparison}

\begin{proposition}[{\mdseries\scshape Single-Location Comparison}]\label{prop:scalar}
Suppose $H\ge h$ and $p\ge h$. For $K>0$, take a canonical optimum of the planning problem with setup cost $5K/3$, and let $\hat w$, $\hat S$, and $b:=\hat w/2$ be its width, target, and half-width. Operate either that full window or the width-$b$ window with target $\hat S-\upsilon^*$, where $\upsilon^*\in[0,b]$ is a maximizing shift from Lemma~\ref{lem:shift}. At actual setup cost $K$, their augmented costs $\cL_{\rm full}:=\cL(2b,0)$ and $\cL_{\rm half}:=\cL(b,\upsilon^*)$ from~\eqref{eq:candidate} satisfy
\begin{equation}\label{eq:scalar-bound}
 \min\{\cL_{\rm full},\cL_{\rm half}\}\le\frac32 J(5K/3).
\end{equation}
For $K=0$, a constant base-stock position at a minimizer of $F$ has augmented cost $J(0)$.
\end{proposition}

\proof{Proof.} The shift and tail-mirror lemmas bound the half-window cost, while the renewal bounds give the full-window bound. We combine these bounds to eliminate the normalized setup term and show that the cheaper candidate costs at most $\frac32J(5K/3)$. Base stock handles zero setup cost.

\emph{Step 1. Bound the half-window cost.} Let $k=5K/3$, $G:=J(k)$, and $\Phi^* =\int_{[0,b)}g(\upsilon^*+x)\,U(dx)$. Use the following ratios, which compare $p$ with $h$ and per-period costs with the planning value $G$: $r:=p/h\ge1$, $\omega:=\beta\hat w/G\le1$, $\kappa:=k/[GM(b)]$, and $R:=\Phi^*/[GM(b)]$, where $r\ge1$ because $p\ge h$, and $\omega\le1$ because Lemma~\ref{lem:reward} gives $\beta \hat w\le G$ when $H\ge h$. We obtain two lower bounds on retained savings. First, Lemma~\ref{lem:shift}, which applies to $g$ by Lemma~\ref{lem:reward}, gives $R\ge\kappa/2$, because $\Phi^*\ge\frac12\int_{[0,2b)}g(x)\,U(dx)=k/2$ by~\eqref{eq:shift} and~\eqref{eq:reward-properties}. Second, splitting $k=\int_{[0,2b)}g(x)\,U(dx)$ at $b$ gives $k=\Phi(0)+\int_{[b,2b)}g(x)\,U(dx)$. Here $\Phi(0)\le\Phi^*$ because $\upsilon^*$ maximizes $\Phi$ in Lemma~\ref{lem:shift}, and $g(x)\le p(2b-x)$ by~\eqref{eq:reward-properties}. Thus, $k\le\Phi^*+p\mathcal T_b$, and~\eqref{eq:tail-mirror} in Lemma~\ref{lem:tail} gives $p\mathcal T_b\le(3r+1)\beta bM(b)-4rM(b)B_b$. Dividing the resulting bound on $k$ by $GM(b)$ and using $\beta b/G=\omega/2$, since $b=\hat w/2$, gives the second lower bound, $R\ge\kappa-{(3r+1)}\omega/2+4rB_b/G$. By Lemma~\ref{lem:mirror}, the half-window cost is exactly $\cL_{\rm half}/G=1+3\kappa/5-R+B_b/G$. To cancel the mirror term in the half-window cost, take the convex combination of the two lower bounds on $R$ with respective weights $1-1/(4r)$ and $1/(4r)$, which lie in $[0,1]$ because $r\ge1/4$. This gives $R-B_b/G\ge(1/2+1/(8r))\kappa-(3r+1)\omega/(8r)$, and substitution into the half-window cost gives
\begin{equation}\label{eq:half-bound}
 \frac{\cL_{\rm half}}G\le
 1+\frac{3r+1}{8r}\omega+\frac{4r-5}{40r}\kappa.
\end{equation}

\emph{Step 2. Bound the full-window cost.} For the full window,~\eqref{eq:candidate} with $w=\hat w$ and $\upsilon=0$, together with $\int_{[0,\hat w)}g(x)\,U(dx)=k$ from~\eqref{eq:reward-properties}, gives $\cL_{\rm full}=G+(K-k)/M(\hat w)+B_{\hat w}$. Its actual setup cost $K=3k/5$ is $2k/5$ below the planning setup cost $k$. The bound $B_{\hat w}\le h\hat w/2$ from Lemma~\ref{lem:mirror} and the doubling bound $M(\hat w)\le2M(b)$ from Lemma~\ref{lem:renewal} then give
\begin{equation}\label{eq:full-bound}
 \frac{\cL_{\rm full}}G=1-\frac{2k}{5GM(\hat w)}+\frac{B_{\hat w}}G\le1-\frac{\kappa}{5}+\frac{h\hat w}{2G}=1+\frac{r+1}{2r}\omega-\frac{\kappa}{5},
\end{equation}
where the last equality uses $\omega=\beta \hat w/G$, so $h\hat w/G=(h/\beta)\omega$, and $h/\beta=(h+p)/p=(r+1)/r$.

\emph{Step 3. Compare the two candidates.} If $1\le r\le5/4$, the coefficient of $\kappa$ in~\eqref{eq:half-bound} is nonpositive, and $(3r+1)/(8r)\le1/2$ because $r\ge1$; this is where the comparison uses $p\ge h$. Using $\omega\le1$, we conclude that the half-window cost is at most $3G/2$. If $r\ge5/4$, the cheaper candidate is at most the convex combination of the two bounds with weights $(4r-5)/(12r-5)\in[0,1)$ on the full window and $8r/(12r-5)$ on the half window. These weights make the $\kappa$ terms in~\eqref{eq:half-bound} and~\eqref{eq:full-bound} cancel, because $(4r-5)/[5(12r-5)]=[8r/(12r-5)](4r-5)/(40r)$, giving $\min\{\cL_{\rm full},\cL_{\rm half}\}/G\le1+A(r)\omega$ with $A(r):=(10r^2+r-5)/[2r(12r-5)]$. On this range $A(r)>0$ and $1/2-A(r)={[2(r-3/2)^2+1/2]/[2r(12r-5)]}>0$, so $A(r)<1/2$, and $\omega\le1$ proves~\eqref{eq:scalar-bound}. The convex combination compares the two costs; the policy chooses the candidate with the lower augmented cost, breaking ties by a fixed rule. Figure~\ref{fig:comparison} plots the two bounds and shows where each candidate is cheaper.

\emph{Step 4. Handle zero setup cost.} If $K=0$, base stock at a minimizer $\Gamma$ of $F$ attains $J(0)=\min F$ by Lemma~\ref{lem:scalar}. It keeps the post-action position at $\Gamma$, so center $a=\Gamma$ gives mirror loss $\phi(a-\widehat P)=\phi(0)=0$ in every period. Hence, the mirror cost is zero, and the augmented cost is $F(\Gamma)=J(0)$.
\Halmos\endproof

\section{The Policy and the Lower Bound}\label{sec:fifo}

Theorem~\ref{thm:main} follows by comparing the physical upper bound of Proposition~\ref{prop:upper} with the shared-cost lower bound of Proposition~\ref{prop:lower}. Both use the same supplier and retailer single-location problems, and Section~\ref{sec:assembly} combines their coefficients uniformly over all finite lead times.

Each retailer runs the policy of Definition~\ref{def:fifo} on its observed demands, as if the warehouse filled every request at once, so that Proposition~\ref{prop:scalar} applies. Write $\widehat P_{i,t}$ for its post-action inventory position including outstanding requests,\footnote{After ordering, all inventory positions including outstanding requests stay in fixed bounded intervals. This keeps the transformed single-location charges bounded, which the proof of Proposition~\ref{prop:upper} in Appendix~\ref{app:initialization} uses for its renewal limits.} $\widehat q_{i,t}\ge0$ for its request in period $t$, and $a_i$ for its fixed inventory-position mirror center. For a positive-width policy of target $\Gamma_i$ and depletion center $z_i$, this means $a_i=\Gamma_i-z_i$. For a retailer with $K_i=0$ choose $a_i=\Gamma_i$, its base-stock level. Retailer $i$'s augmented average cost is
\begin{equation}\label{eq:local-upper}
 \cL_i:=\lim_{T\to\infty}\frac1T\E\sum_{t=0}^{T-1}
 \left[K_i\ind_{\{\widehat q_{i,t}>0\}}+F_i(\widehat P_{i,t})+\phi(a_i-\widehat P_{i,t})\right].
\end{equation}
By~\eqref{eq:occupation}, the limit exists and, for $K_i>0$, equals the augmented cost~\eqref{eq:candidate} of the selected candidate, because the fixed center $a_i$ attains the minimum in~\eqref{eq:mirror}. If $H_i\ge h$ and $p\ge h$, Definition~\ref{def:fifo} and~\eqref{eq:scalar-bound} give $\cL_i\le3J_i(5K_i/3)/2$. For $K_i=0$, the post-action position stays at the base-stock level $\Gamma_i=a_i$, a minimizer of $F_i$, so the mirror loss is zero and $\cL_i=F_i(\Gamma_i)=J_i(0)$ by Lemma~\ref{lem:scalar}.

Let $D^{(L_0)}$ be a sum of $L_0$ independent aggregate demands, with $D^{(0)}=0$. With the warehouse-rate loss $\phi$, define the supplier holding and backlog cost and its single-location value by
\begin{equation}\label{eq:source-value}
 f_{L_0}(y):=\E\phi(y-D^{(L_0)}),\qquad
 J_0(k):=\J_D(f_{L_0};k).
\end{equation}
Lemma~\ref{lem:scalar} applies to this supplier problem.\footnote{Jensen's inequality gives $f_{L_0}(y)\ge\phi(y-L_0\E D)$, so the nonnegative, convex, and Lipschitz function $f_{L_0}$ is coercive.} Supplier cost uses $L_0$ demands and retailer cost uses $L_i+1$ because they are measured at different points of a period. Supplier cost in period $t+L_0$, when an order placed in period $t$ arrives, is measured after dispatches and before demand is realized, so only the demands of periods $t$ through $t+L_0-1$ count. Retailer cost is measured after demand, so the demand of the arrival period $t+L_i$ also counts, giving $Z_i$ its $L_i+1$ terms in~\eqref{eq:retailer-value}.

\subsection{The Queue and Aggregate Inventory Position}\label{sec:queue}

\begin{figure}[htbp]
\SingleSpaced
\centering
\includegraphics[width=.99\linewidth]{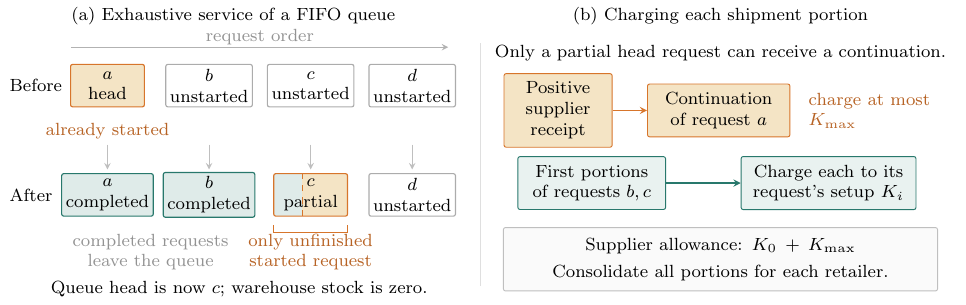}
\caption{Exhaustive FIFO service, a partly filled request, and the charge for each shipment portion.}
\label{fig:fifo}
\begin{minipage}{0.98\linewidth}
\SingleSpaced\footnotesize\emph{Notes.} After the prescribed pipelines arrive, only the started request at the head of the queue can receive a continuation at a receipt. Every other request reached receives its first portion, charged to its own setup. A partly filled request can span several receipts; each continuation is charged at most $K_{\max}$ to its supplier order.
\end{minipage}
\end{figure}

Under Definition~\ref{def:fifo}, exhaustive FIFO service completes successive requests and allows a partial final request, with all portions sent to one retailer in a period sharing one setup cost. Figure~\ref{fig:fifo} shows the queue, a partly filled request, and the continuation charge that Lemma~\ref{lem:continuation} bounds.

Let $\Delta_{i,t}$ be the quantity of retailer $i$'s requests awaiting dispatch at the warehouse, $\Delta_t=\sum_i \Delta_{i,t}$, and $O_t$ the outstanding supplier pipeline after the current ordering and receipts. For $L_0=0$, $O_t=0$. The invariants and the aggregate positions are
\begin{align}
 P_{i,t}&=\widehat P_{i,t}-\Delta_{i,t},\qquad W_t\Delta_t=0,\label{eq:queue-invariants}\\
 E_t&:=W_t+\sum_i(P_{i,t}-a_i)
     =W_t-\Delta_t+\sum_i(\widehat P_{i,t}-a_i),\qquad Y_t:=E_t+O_t.
 \label{eq:aggregate}
\end{align}
Every quantity indexed by $t$ in~\eqref{eq:queue-invariants} and~\eqref{eq:aggregate} is measured after dispatches and before the demand of period $t$ is realized. The retailer and supplier rules apply to their positions before the current request and supplier order, respectively. Requests, dispatches, receipts, and demand each preserve $P_i=\widehat P_i-\Delta_i$,\footnote{A request raises $\widehat P_i$ and $\Delta_i$ equally, and a dispatch raises $P_i$ while reducing $W$ and $\Delta_i$ equally. A retailer receipt leaves both $P_i$ and $\widehat P_i$ unchanged. Demand lowers both retailer positions $P_i$ and $\widehat P_i$ equally.} and exhaustive service gives $W\Delta=0$. The definition of $E_t$ rearranges to the mirror identity $W_t-\Delta_t=E_t+\sum_i(a_i-\widehat P_{i,t})$: warehouse stock net of waiting requests equals $E_t$, which the supplier rule controls, plus each retailer's position reflected about its mirror center. Lemma~\ref{lem:running} charges each reflected position its mirror loss, using the same fixed center as the augmented cost.

On the aggregate position $Y_t$ from~\eqref{eq:aggregate}, run the optimal $(s,S)$ policy of Lemma~\ref{lem:scalar} for the single-location problem with cost $f_{L_0}$ from~\eqref{eq:source-value} and setup cost $K_0+K_{\max}$. The pre-order value of $Y$ is known before any dispatch, since dispatches leave it unchanged. Supplier orders raise $Y$, supplier receipts move units from $O$ to $W$, and demand in period $t$ lowers $Y$ by $D_t$. Consequently, $Y_{t+1}=Y_t-D_t+q_{0,t+1}$, and $E_{t+L_0}=Y_t-\sum_{\ell=t}^{t+L_0-1}D_\ell$ because $E=Y-O$, every order already in $O_t$ arrives by $t+L_0$, and every order placed after $t$ arrives later. At zero supplier lead time, this identity reads $E_t=Y_t$. Since the next $L_0$ aggregate demands are independent of the history $\mathcal F_t$ through the dispatches of period $t$, which determines $Y_t$, the identity for $E_{t+L_0}$ and~\eqref{eq:source-value} give
\begin{equation}\label{eq:source-expectation}
 \E[\phi(E_{t+L_0})\mid\mathcal F_t]=f_{L_0}(Y_t).
\end{equation}
A dispatch lowers $W$ and raises one $P_i$ by the same amount, so it changes neither $E$ nor $Y$. Hence, the identities for $Y$ and $E$ hold for every dispatch decision and every finite supplier lead time. In particular, the supplier's orders do not depend on how the warehouse serves the FIFO queue.

\subsection{Bounds on Setup, Holding, and Backlog Costs}

\begin{lemma}[{\mdseries\scshape Continuation Charge}]\label{lem:continuation}
Suppose the initial queue holds at most $n$ requests, none partly filled, and let $N_{\rm init}$ count positive receipts in the initial supplier pipeline. For every $T\ge1$,
\begin{equation}\label{eq:continuation-charge}
 \sum_{t<T}\sum_iK_i\ind_{\{q_{i,t}>0\}}
 \le\sum_{t<T}\sum_iK_i\ind_{\{\widehat q_{i,t}>0\}}
 +K_{\max}\left(n+N_{\rm init}+\sum_{t<T}\ind_{\{q_{0,t}>0\}}\right).
\end{equation}\par
\end{lemma}

\proof{Proof.} Exhaustive FIFO service permits at most one continuation per positive supplier receipt. Charging first portions to requests and continuations to receipts bounds dispatch setup cost.

\emph{Step 1. Bound the charge in each period.} Charge the first positive portion of each request to that request's setup cost $K_i$. At the end of exhaustive service, at most one started request can remain unfinished. Such a request is at the head of the queue, and the warehouse has no stock left. Filling any further portion of that request requires a positive supplier receipt. In any receipt period, only the previously started request at the head of the queue can receive a continuation; every additional request reached in that period receives its first portion. Hence, the dispatch setup cost $\sum_iK_i\ind_{\{q_{i,t}>0\}}$ of period $t$ is at most the total setup cost of the requests receiving their first positive portion in period $t$, plus $K_{\max}$ if a positive supplier receipt arrives. A request may span many receipts, each continuation charged separately to that receipt's supplier order. Shipping all portions for a retailer in a period as one dispatch can only reduce the actual setup cost.

\emph{Step 2. Sum the charges over the horizon.} First portions of the at most $n$ initial requests cost at most $nK_{\max}$, and continuations at receipts from the initial pipeline cost at most $K_{\max}N_{\rm init}$. Every other positive receipt in the horizon comes from a positive supplier order placed in the horizon. Each such order produces at most one receipt, which can trigger at most one continuation. Thus, summing the per-period bound gives~\eqref{eq:continuation-charge}, with each receipt charged at most once.
\Halmos\endproof

\begin{lemma}[{\mdseries\scshape FIFO Holding and Backlog Cost Bound}]\label{lem:running}
Every period after the initialization in Appendix~\ref{app:initialization} satisfies $hW_t+\sum_iF_i(P_{i,t})\le\phi(E_t)+\sum_i{[F_i(\widehat P_{i,t})+\phi(a_i-\widehat P_{i,t})]}$.
\end{lemma}

\proof{Proof.} The retailer slope bound charges waiting requests at the backlog rate. Exhaustive service and subadditivity express this charge through the pooled position and the individual mirror losses.

\emph{Step 1. Compare the retailer positions.} By~\eqref{eq:queue-invariants}, $P_i=\widehat P_i-\Delta_i$, and $F_i$ has slopes in $[-p,H_i]$, so $F_i(P_i)\le F_i(\widehat P_i)+p\Delta_i$.

\emph{Step 2. Combine the warehouse and mirror losses.} Since $W,\Delta\ge0$ and $W\Delta=0$, $hW+p\Delta=\phi(W-\Delta)$. Since $\phi$ is subadditive and the mirror identity gives $W-\Delta=E+\sum_i(a_i-\widehat P_i)$, $\phi(W-\Delta)\le\phi(E)+\sum_i\phi(a_i-\widehat P_i)$. Summing the retailer bounds and using $\sum_i\Delta_i=\Delta$ gives $hW+\sum_i F_i(P_i)\le\sum_i F_i(\widehat P_i)+hW+p\Delta$. Substitute $\phi(W-\Delta)$ for $hW+p\Delta$ and apply the bound on $\phi(W-\Delta)$ above to obtain the stated bound in every subsequent period. This holding and backlog bound uses only exhaustive service; FIFO yields the continuation bound in Lemma~\ref{lem:continuation}.
\Halmos\endproof

\begin{proposition}[{\mdseries\scshape Physical Upper Bound}]\label{prop:upper}
From every valid initial state, the deterministic policy $\pi_{\rm FIFO}$ satisfies
\begin{equation}\label{eq:upper}
 C(\pi_{\rm FIFO})\le J_0(K_0+K_{\max})+\sum_i\cL_i,
\end{equation}
where $J_0(K_0+K_{\max})$ is the optimal average cost of the single-location problem with cost $f_{L_0}$ in~\eqref{eq:source-value} and setup cost $K_0+K_{\max}$, and $\cL_i$ is retailer $i$'s augmented average cost in~\eqref{eq:local-upper}.
\end{proposition}

Appendix~\ref{app:initialization} proves the proposition by combining three ingredients. First, Lemma~\ref{lem:continuation} bounds physical setup cost by the setup cost of requests plus $K_{\max}$ per supplier order. Second, Lemma~\ref{lem:running} bounds each period's holding and backlog cost. Third,~\eqref{eq:source-expectation} turns the resulting warehouse term $\phi(E_t)$ into the expected one-period cost $f_{L_0}(Y_{t-L_0})$ of the supplier problem. Taking long-run averages then gives $J_0(K_0+K_{\max})$ for the supplier terms and $\cL_i$ for each retailer's terms. A preparation phase of finite expected duration and cost extends this argument to every valid initial state.

\subsection{A Shared-Cost Lower Bound}\label{sec:lower}

\begin{proposition}[{\mdseries\scshape Shared-Cost Lower Bound}]\label{prop:lower}
Suppose $H_i\ge h$ for every retailer and the backlog penalty $p>0$ is common. Then for every real $\gamma>1$,
\begin{equation}\label{eq:lower}
 J_0(\gamma K_0)
 +(\gamma-1)\sum_iJ_i\left(\frac{\gamma K_i}{\gamma-1}\right)\le\gamma\OPT.
\end{equation}
\end{proposition}

\begin{figure}[htbp]
\SingleSpaced
\centering
\includegraphics[width=.99\linewidth]{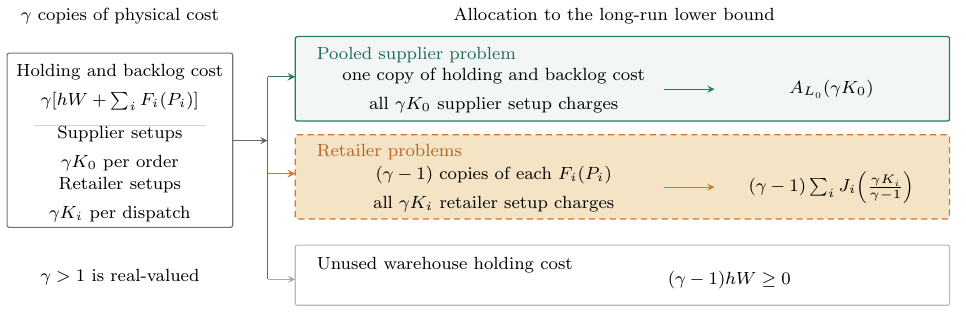}
\caption{Allocating $\gamma$ copies of physical cost to pooled supplier and individual retailer lower bounds.}
\label{fig:allocation}
\begin{minipage}{0.98\linewidth}
{\SingleSpaced\footnotesize\emph{Notes.} Each retailer cost $F_i(P_i)$ is the expected physical cost $L_i$ periods later. One copy of the holding and backlog cost, with all supplier setups, is bounded through the supplier problem; the other $\gamma-1$ copies of retailer costs, with all retailer setups, are bounded through the retailer problems. The arrows give long-run lower bounds after dividing bounded endpoint terms by the horizon. Unused warehouse holding cost is nonnegative, and the same allocation is valid for any real multiplier $\gamma>1$, which Theorem~\ref{thm:main} uses for every $\gamma\ge\gamma_m$.}
\end{minipage}
\end{figure}

Figure~\ref{fig:allocation} shows the split of $\gamma$ copies of physical cost. In each period, $\gamma{[hW_t+\sum_iF_i(P_{i,t})]}={[hW_t+\sum_iF_i(P_{i,t})]}+{(\gamma-1)}\sum_iF_i(P_{i,t})+{(\gamma-1)}hW_t$. On the right, the first term, with $\gamma K_0$ per supplier order, is bounded through the pooled supplier problem by Step~1 below and~\eqref{eq:source-lower}. The second term, with $\gamma K_i$ per dispatch, is bounded through $\gamma-1$ copies of each retailer problem at setup cost $\gamma K_i/(\gamma-1)$ by~\eqref{eq:local-lower}. Since $\gamma>1$, the last term $(\gamma-1)hW_t$ is nonnegative and is dropped.

\proof{Proof.} We allocate physical costs between the retailer problems and a pooled supplier problem. The pooled position identities of Step~1 justify the supplier lead-time cost, and the finite-horizon bound~\eqref{eq:finite-horizon} keeps the endpoint correction uniform over all feasible policies and horizons.

\emph{Step 1. Construct the pooled supplier comparison.} Fix an arbitrary nonanticipative feasible policy, possibly history-dependent or independently randomized. Let $\mathcal F_t$ contain the observations through the dispatches of period $t$, before its demand is realized, and the auxiliary randomness used so far. Because $H_i\ge h$, the definitions~\eqref{eq:retailer-value} of $F_i$ and of $\phi$ imply $\E\phi(P_i-Z_i)\le F_i(P_i)$, and Jensen's inequality yields $\phi(P_i-\E Z_i)\le\E\phi(P_i-Z_i)$. Since $W\ge0$, $hW=\phi(W)$. Adding these bounds and using the subadditivity of $\phi$ gives the pointwise inequality $\phi({W_t+\sum_i(P_{i,t}-\E Z_i)})\le hW_t+\sum_iF_i(P_{i,t})$.

Define the pooled positions $\widetilde E_t=W_t+\sum_i(P_{i,t}-\E Z_i)$ and $\widetilde Y_t=\widetilde E_t+O_t$. Every dispatch leaves both unchanged, so the argument of Section~\ref{sec:queue} applies to every feasible nonanticipative policy and gives $\widetilde Y_{t+1}=\widetilde Y_t-D_t+q_{0,t+1}$ and $\widetilde E_{t+L_0}=\widetilde Y_t-\sum_{\ell=t}^{t+L_0-1}D_\ell$. Since the next $L_0$ aggregate demands are independent of $\mathcal F_t$ and of the policy's auxiliary randomness, the identity for $\widetilde E_{t+L_0}$ and~\eqref{eq:source-value} give $\E[\phi(\widetilde E_{t+L_0})\mid\mathcal F_t]=f_{L_0}(\widetilde Y_t)$. This counterpart of~\eqref{eq:source-expectation} for an arbitrary policy is used in~\eqref{eq:source-lower} below.

\emph{Step 2. Apply the retailer lower bounds.} For an explicit finite-horizon allocation over $T$ periods, write $N=T+L_0+\max_iL_i$. The finite-horizon bound~\eqref{eq:finite-horizon} of Lemma~\ref{lem:scalar} applies to the positions $P_{i,t}$ and $\widetilde Y_t$.\footnote{Each $P_{i,t}$ rises only through dispatches and $\widetilde Y_t$ only through supplier orders, and both fall by demand, so each follows a nonanticipative single-location policy that also observes the network history, as Lemma~\ref{lem:scalar} allows. Its bound holds from every initial position, so centering the supplier positions at $\E Z_i$ instead of $a_i$ leaves the value $J_0$ unchanged.} For retailer $i$, apply~\eqref{eq:finite-horizon} with cost $F_i$, setup cost $\gamma K_i/(\gamma-1)$, and horizon $T+L_0$. Multiply the resulting bound by $\gamma-1$, thereby using $\gamma-1$ copies of each retailer's holding and backlog cost, and sum over retailers to obtain
\begin{equation}\label{eq:local-lower}
 (\gamma-1)(T+L_0)\sum_iJ_i\left(\frac{\gamma K_i}{\gamma-1}\right)-\gamma\sum_iK_i
 \le\E\sum_{t=0}^{T+L_0-1}\sum_i
 \left[\gamma K_i\ind_{\{q_{i,t}>0\}}+(\gamma-1)F_i(P_{i,t})\right].
\end{equation}\par

\emph{Step 3. Apply the supplier lower bound.} It remains to bound one copy of $hW_t+\sum_iF_i(P_{i,t})$, together with $\gamma K_0$ per supplier order. To apply~\eqref{eq:finite-horizon} to that problem, whose one-period cost is $f_{L_0}(\widetilde Y_t)$, we first convert this physical cost into that form. By Step~1, $\phi(\widetilde E_t)\le hW_t+\sum_iF_i(P_{i,t})$ in every period, and by the conditional expectation in Step~1, $\E f_{L_0}(\widetilde Y_t)=\E\phi(\widetilde E_{t+L_0})$. Summing over $t<T+L_0$ and dropping the nonnegative terms $\phi(\widetilde E_0),\ldots,\phi(\widetilde E_{L_0-1})$ gives $\E\sum_{t=0}^{T-1}f_{L_0}(\widetilde Y_t)\le\E\sum_{t=0}^{T+L_0-1}(hW_t+\sum_iF_i(P_{i,t}))$. Applying~\eqref{eq:finite-horizon} with cost $f_{L_0}$, setup cost $\gamma K_0$, and horizon $T$, and then adding the supplier setup charges $\gamma K_0$ over $T$ periods to this inequality, gives
\begin{equation}
 T J_0(\gamma K_0)-\gamma K_0\le\E\Bigl[\sum_{t=0}^{T-1}\gamma K_0\ind_{\{q_{0,t}>0\}}+\sum_{t=0}^{T+L_0-1}\Bigl(hW_t+\sum_iF_i(P_{i,t})\Bigr)\Bigr].\label{eq:source-lower}
\end{equation}

\emph{Step 4. Combine the bounds and take long-run averages.} By~\eqref{eq:retailer-lead}, each charged term $F_i(P_{i,t})$ with $t<T+L_0$ is the conditional expected holding and backlog cost of retailer $i$ in period $t+L_i<N$, and all charged setups and warehouse costs also occur before period $N$. Hence, the sum of the right sides of~\eqref{eq:local-lower} and~\eqref{eq:source-lower} is at most $\gamma$ times the expected physical cost over periods $0$ to $N-1$: it charges each retailer setup $\gamma K_i$, each supplier setup $\gamma K_0$, and each retailer holding and backlog cost $(\gamma-1)+1=\gamma$ times, while the warehouse holding cost $hW_t$ appears only once, in~\eqref{eq:source-lower}, and $hW_t\le\gamma hW_t$. Physical costs in the periods that the right-side sums do not cover are nonnegative, so adding the left sides of~\eqref{eq:local-lower} and~\eqref{eq:source-lower} gives\footnote{If the expected physical cost $\E_\pi\sum_{t=0}^{N-1}c_t$ on the right side of the following inequality, with $c_t$ the period cost in~\eqref{eq:cost}, is infinite, the inequality holds trivially. Otherwise, the expectations on the right sides of~\eqref{eq:local-lower} and~\eqref{eq:source-lower} are finite, because each is at most $\gamma$ times that cost, so adding the two lower bounds is legitimate.} $TJ_0(\gamma K_0)+{(\gamma-1)(T+L_0)}\sum_iJ_i({\gamma K_i/(\gamma-1)})-\gamma{(K_0+\sum_iK_i)}\le\gamma\E_\pi\sum_{t=0}^{N-1}c_t$.
The subtracted constant $\gamma(K_0+\sum_iK_i)$ comes from the bounded functions in Lemma~\ref{lem:scalar} and does not depend on terminal stocks or pipelines. Since the lead times are fixed, dividing by $N$, taking the limit superior as $T\to\infty$, and then the infimum over nonanticipative policies proves the shared-cost lower bound~\eqref{eq:lower} for every real multiplier $\gamma>1$.
\Halmos\endproof

\subsection{Assembly of the Guarantee}\label{sec:assembly}

\proof{Proof of Theorem~\ref{thm:main}.} We retain the holding-rate dependence in the retailer comparison, then align the resulting upper bound with the shared-cost lower bound of Proposition~\ref{prop:lower}.

\emph{Step 1. Bound each retailer's augmented cost.} Fix a retailer with $K>0$. Set $k=5K/3$, let $G=J(k)>0$ and $\hat w$ be its canonical value and width, and write $r=p/h\ge1$ and $\omega=\beta\hat w/G$, where $\beta=hp/(h+p)$. The endpoint calculation in Step~2 of the proof of Lemma~\ref{lem:reward} gives
\begin{equation}\label{eq:holding-width}
 \omega\le\frac{h(H+p)}{H(h+p)}
 \le\frac{m+r}{m(1+r)}
 \le\frac{m+1}{2m},
\end{equation}
where the second inequality uses $H\ge mh$, and the third follows from $(m+1)/(2m)-(m+r)/[m(1+r)]=(m-1)(r-1)/[2m(1+r)]\ge0$. The proof of Proposition~\ref{prop:scalar} gives $\min\{\cL_{\rm full},\cL_{\rm half}\}/G\le1+\omega/2$: for $1\le r\le5/4$, the setup-cost coefficient in~\eqref{eq:half-bound} is nonpositive and its width coefficient $(3r+1)/(8r)$ is at most $1/2$; for $r\ge5/4$, Step~3 of that proof gives $1+A(r)\omega$ with $A(r)<1/2$. Thus,
\begin{equation}\label{eq:holding-local}
 \cL=\min\{\cL_{\rm full},\cL_{\rm half}\}
 \le\left(1+\frac\omega2\right)G
 \le\left(1+\frac{m+1}{4m}\right)G
 =(\gamma_m-1)J(5K/3).
\end{equation}
For $K=0$, base stock gives $\cL=J(0)\le(\gamma_m-1)J(0)$, so this bound also covers zero setup costs.

\emph{Step 2. Align the retailer setup costs.} Single-location values are nondecreasing in setup cost, and multiplying a setup cost by $a\ge1$ multiplies its value by at most $a$: the policy optimal at the lower setup cost remains feasible and has nonnegative running costs. If $\gamma_m\le\gamma\le5/2$, then $\gamma/(\gamma-1)\ge5/3$ and $\gamma_m-1\le\gamma-1$, so monotonicity and~\eqref{eq:holding-local} give
\[
 \cL_i\le(\gamma_m-1)J_i(5K_i/3)
 \le(\gamma-1)J_i\!\left(\frac{\gamma K_i}{\gamma-1}\right).
\]
If $\gamma\ge5/2$, scaling by $a=5(\gamma-1)/(3\gamma)\ge1$ gives the same comparison:
\[
 \cL_i\le(\gamma_m-1)\frac{5(\gamma-1)}{3\gamma}
 J_i\!\left(\frac{\gamma K_i}{\gamma-1}\right)
 \le(\gamma-1)J_i\!\left(\frac{\gamma K_i}{\gamma-1}\right),
\]
since $\gamma_m-1\le3/2$ and $5/(2\gamma)\le1$. Hence the retailer comparison holds for every $\gamma\ge\gamma_m$.

\emph{Step 3. Combine the network bounds.} The setup assumption and monotonicity give $J_0(K_0+K_{\max})\le J_0(\gamma K_0)$. Substituting this and Step~2 into Proposition~\ref{prop:upper}, then applying Proposition~\ref{prop:lower}, yields
\[
 C(\pi_{\rm FIFO})
 \le J_0(K_0+K_{\max})+\sum_i\cL_i
 \le J_0(\gamma K_0)+(\gamma-1)\sum_iJ_i\!\left(\frac{\gamma K_i}{\gamma-1}\right)
 \le\gamma\OPT.
\]
The unchanged policy and initialization of Definition~\ref{def:fifo} satisfy Proposition~\ref{prop:upper} from every valid initial state and at all finite lead times. If $K_0=0$, the setup assumption forces every $K_i=0$; both network bounds apply to the resulting base-stock policies, so the same comparison covers this case.
\Halmos\endproof

\section{Arbitrary Setup Costs}\label{sec:arbitrary}

For positive supplier setup cost, the FIFO factor $\max\{\gamma_m,1+K_{\max}/K_0\}$ with $m=\min_iH_i/h$ from Theorem~\ref{thm:main} grows with $K_{\max}/K_0$. To obtain a factor that remains bounded as retailer setup costs increase, we consolidate dispatches. The consolidation policy keeps the FIFO queue of Definition~\ref{def:fifo} but holds the stock allocated to a request until its cumulative allocation reaches the next of at most two dispatch thresholds, so each request ships in at most three portions. This is a quantity-based counterpart of the time-based shipment consolidation studied by \citet{Marklund2011}, in which the warehouse dispatches at regular intervals; the analysis below assumes associated demands.

We assume that each period's demand vector $\mathbf D_t=(D_{1,t},\ldots,D_{n,t})$ is \emph{associated}, that is, $\operatorname{Cov}(u(\mathbf D_t),v(\mathbf D_t))\ge0$ for all bounded, measurable, coordinatewise nondecreasing functions $u,v$ \citep{EsaryProschanWalkup1967}. This includes independent demands and nondecreasing functions of independent common and retailer-specific shocks. Demand vectors remain i.i.d.\ over periods.

Fix $m\ge1$ with $H_i\ge mh$ for every retailer. Define
\begin{equation}\label{eq:con-parameters}
 \bar\gamma_m:=\frac{13}{3}+\frac{32}{57m},\quad
 g_m:=\frac{\bar\gamma_m}{2},\quad
 \zeta_m:=\frac{g_m}{g_m-1},\quad
 B_m:=2(g_m-1),\quad
 u_m:=\frac{3\zeta_m(1+1/m)}{4(4+\zeta_m/m)}.
\end{equation}
Plan retailer $i$ at setup cost $k_i=\zeta_mK_i$ and write $G_i=J_i(k_i)$ and $\Lambda_i=\min_xF_i(x)$. If $K_i>0$, let $\hat w_i=\hat S_i-\hat s_i$ be the canonical width of its optimum at setup cost $k_i$. Its setup-cost share is $u_i=k_i/[G_iM_i(\hat w_i)]$, with $M_i$ defined in~\eqref{eq:renewal}. These planning costs depend on $m$ and differ from those of the FIFO policy in Definition~\ref{def:fifo}, which plans every retailer at $5K_i/3$.

\begin{definition}[{\mdseries\scshape Consolidation Policy}]\label{def:con}
For the parameters in~\eqref{eq:con-parameters}, the deterministic policy $\pi_{\rm con}^{(m)}$ follows four rules after the initialization in Appendix~\ref{sec:arb-initial} and breaks ties by fixed rules.

(i)~\emph{Retailer policies.} Retailer $i$ uses base stock at a minimizer of $F_i$ if $\Lambda_i+K_i\le B_mG_i$. Otherwise, it uses the half window, the $(s,S)$ pair with the smallest average holding and backlog cost subject to $\hat s_i\le s<S\le\hat S_i$ and $S-s=\hat w_i/2$, if $u_i\le u_m$, and the full window $[\hat s_i,\hat S_i]$ otherwise.

(ii)~\emph{Requests and allocation.} Retailer $i$ applies its policy to $\widehat P_{i,t}=P_{i,t}+\Delta_{i,t}+R_{i,t}$, where $\Delta_{i,t}$ is the quantity of its requests not yet allocated and $R_{i,t}$ is their reserve, the stock allocated to them but still held at the warehouse. Requests join one FIFO queue, and in each period, after the supplier order, the warehouse allocates its unallocated stock to requests in queue order, the last one possibly in part. A request's cumulative allocation is the total allocated to it so far.

(iii)~\emph{Dispatch.} Allocations to base-stock retailers are dispatched at once. Every other request, of size $Q$, receives, when it is created, an offset $\vartheta\in(0,1)$ set by the rule of Appendix~\ref{sec:arb-derandomization}, which fixes its dispatch thresholds for the life of the request: $3Q\vartheta/4$ and $3Q(1+\vartheta)/4$ if $\vartheta<1/3$, and only $3Q\vartheta/4$ if $1/3\le\vartheta<1$. In each period in which the request's cumulative allocation rises to or past a threshold or reaches $Q$, the warehouse ships the request's entire reserve in one portion; otherwise, it holds the reserve. All portions sent to the same retailer in a period form one dispatch.

(iv)~\emph{Supplier orders.} Let $(\hat s_0,\hat S_0)$ be the $(s,S)$ policy of Lemma~\ref{lem:scalar} that attains the supplier value $J_0(g_mK_0)$ in~\eqref{eq:source-value}, and let $a_0=(-\hat s_0)^+$. The warehouse orders from the supplier with the $(s,S)$ policy $(\hat s_0+a_0,\hat S_0+a_0)$, or with base stock at a nonnegative minimizer of $f_{L_0}$ if $K_0=0$. It applies this policy to $W_t+O_t+\sum_i(P_{i,t}-\Gamma_i)$, where $\Gamma_i$ is the target selected for retailer $i$ in (i).\Halmos
\end{definition}

Retailers with $K_i=0$ use base stock, because $\Lambda_i+K_i=J_i(0)=G_i\le B_mG_i$. We call base-stock retailers \emph{immediate} and the others, whose requests ship at dispatch thresholds, \emph{batch} retailers. Since stock is allocated in queue order, only the head of the queue can be partly allocated, so the warehouse holds a reserve for at most one request at a time. The threshold spacing keeps this reserve at most $3Q/4$ for a request of size $Q$; Appendix~\ref{sec:arb-policy} gives a numerical example.

To control how long allocated stock can remain at the warehouse, rule (iv) raises both levels of $(\hat s_0,\hat S_0)$ by $a_0$, which is positive only when $\hat s_0<0$, so the supplier position is nonnegative after each ordering decision and every request is fully allocated within $L_0$ periods of its creation. Appendices~\ref{sec:arb-policy} and~\ref{sec:arb-source-budget} use this to bound the reserve cost by the supplier's lead-time newsvendor cost, and the raise of $a_0$ adds at most $ha_0$ to the supplier's long-run average cost.

\begin{theorem}[{\mdseries\scshape Arbitrary Setup Costs}]\label{thm:arbitrary}
Let $m\ge1$. Suppose the retailer demands in each period are associated, $H_i\ge mh$ for every retailer, and $p\ge h$. For arbitrary nonnegative supplier and retailer setup costs and all finite supplier and retailer lead times, the deterministic policy $\pi_{\rm con}^{(m)}$ of Definition~\ref{def:con} satisfies, from every valid initial state,
\[
 C(\pi_{\rm con}^{(m)})\le \bar\gamma_m\OPT,
 \qquad \bar\gamma_m=\frac{13}{3}+\frac{32}{57m}.
\]
\end{theorem}

At $m=1$, the parameters are $g_1=93/38$, $\zeta_1=93/55$, $B_1=55/19$, and $u_1=279/626$, giving $\bar\gamma_1=93/19<4.9$. As $m$ increases, the guarantee decreases to $13/3$; for example, $m=5/3$ gives $1331/285\approx4.6702$, and $m=2$ gives $263/57\approx4.6140$. The policy retunes its planning costs and retailer-selection thresholds with $m$; its dispatch spacing remains $3Q/4$ for every $m$.

Appendix~\ref{sec:arb-source-budget} proves the same guarantee under a weaker condition on the total demands of batch and immediate retailers. Positive quadrant dependence of these two totals suffices. This condition restricts only the dependence between the two totals, not within either group.

Retaining allocated stock raises warehouse holding cost and can delay retailer service. The proof pays for this extra cost with slack in the supplier comparison. When a request waits for allocation, the warehouse has no unallocated stock, so its backlog loss is $(h+p)X_t$ below the supplier-problem loss plus the charge $hX_t$ on the retailers' total depletion $X_t$. After spending this slack on the reserve $R_t$, only the credited reserve $(R_t-X_t)^+$ in~\eqref{eq:arb-credited} remains to be bounded. The supplier cost and part of this reserve cost are at most $2J_0(g_mK_0)$. Retaining $h/H_i\le1/m$ in the width comparison bounds each retailer's remaining contribution by $B_mG_i$. Hence,
\[
 C(\pi_{\rm con}^{(m)})
 \le2J_0(g_mK_0)+B_m\sum_iJ_i(\zeta_mK_i)
 \le2g_m\OPT=\bar\gamma_m\OPT,
\]
where the last inequality doubles~\eqref{eq:lower} at $\gamma=g_m$, using $B_m=2(g_m-1)$ and $\zeta_m=g_m/(g_m-1)$, and the equality uses $\bar\gamma_m=2g_m$. Appendix~\ref{app:arbitrary} proves the retailer comparison for every $m\ge1$, together with initialization and deterministic offset selection under finite mean demand. The choice of $\bar\gamma_m$ gives a simple sufficient bound; we do not claim that this guarantee is tight.

To prove the first inequality above, we assign an independent uniform offset $\vartheta$ to each request. This gives a request with thresholds at most $7/3$ portions in expectation and, for every allocation path, an expected reserve of at most $3Q/8$, half the spacing. With the preceding cost comparison, these estimates yield that inequality for the randomized policy. Appendix~\ref{sec:arb-derandomization} replaces the uniform offset by a deterministic rule based on the state at request creation; this nonanticipative rule preserves the expected-cost bound up to errors vanishing in the long-run average. This deterministic rule establishes existence: for each new request, it selects an offset whose expected setup and reserve cost, averaged over all retailers' demands in the next $L_0$ periods, is at most its average over uniform offsets plus a margin that halves from one request to the next. The randomized version requires only one uniform draw per request and satisfies the same expected-average-cost guarantee.

\begin{corollary}[{\mdseries\scshape Selected Policy}]\label{cor:select}
Suppose the retailer demands in each period are associated, $H_i\ge h$ for every retailer, $p\ge h$, and $K_0>0$. Set $m=\min_iH_i/h$, $\gamma_m=9/4+1/(4m)$, and $\bar\gamma_m=13/3+32/(57m)$. Let $\pi_{\rm select}$ use $\pi_{\rm FIFO}$ if $1+K_{\max}/K_0\le \bar\gamma_m$ and $\pi_{\rm con}^{(m)}$ otherwise. Then
\[
 C(\pi_{\rm select})\le
 \min\!\left\{\bar\gamma_m,\;\max\!\left\{\gamma_m,\;1+\frac{K_{\max}}{K_0}\right\}\right\}\OPT.
\]
\end{corollary}
The corollary follows from Theorems~\ref{thm:main} and~\ref{thm:arbitrary}, since $\gamma_m\le5/2<13/3<\bar\gamma_m$. The selected policy retains the factor $\gamma_m$ when $K_{\max}\le(\gamma_m-1)K_0$, including $5/2$ at $m=1$. Above this setup threshold, the FIFO factor $1+K_{\max}/K_0$ applies until it reaches $\bar\gamma_m$; beyond that point, consolidation gives the factor $\bar\gamma_m$, independent of setup-cost ratios. If $K_0=0$, use $\pi_{\rm FIFO}$ when $K_{\max}=0$ and $\pi_{\rm con}^{(m)}$ otherwise, giving factors $\gamma_m$ and $\bar\gamma_m$, respectively, since Theorems~\ref{thm:main} and~\ref{thm:arbitrary} cover this case.

\section{Numerical Illustration on Retail Data}\label{sec:numerics}

Using an integer-rounded FIFO implementation on daily retail sales data with a common penalty, we assess three aspects of performance. First, we compare the FIFO policy's simulated cost with its analytical upper bound and the shared-cost lower bound of Proposition~\ref{prop:lower}. Second, we compare it with an uncoordinated $(s,S)$ policy and tuned echelon $(R,nQ)$ and integer-ratio policies. Third, we compare it with dynamic-programming benchmarks for two stores. All instances have $K_{\max}<K_0$, where $\pi_{\rm FIFO}$ has the smaller factor in Corollary~\ref{cor:select}, so we simulate its integer-rounded implementation. Figure~\ref{fig:m5} shows the three store networks and their daily sales. Appendix~\ref{app:numerics} gives the item selection, dynamic-programming computation, and further checks.

\subsection{Instances and Bounds}\label{sec:numerics-instances}

\begin{figure}[htbp]
\SingleSpaced
\centering
\includegraphics[width=.99\linewidth]{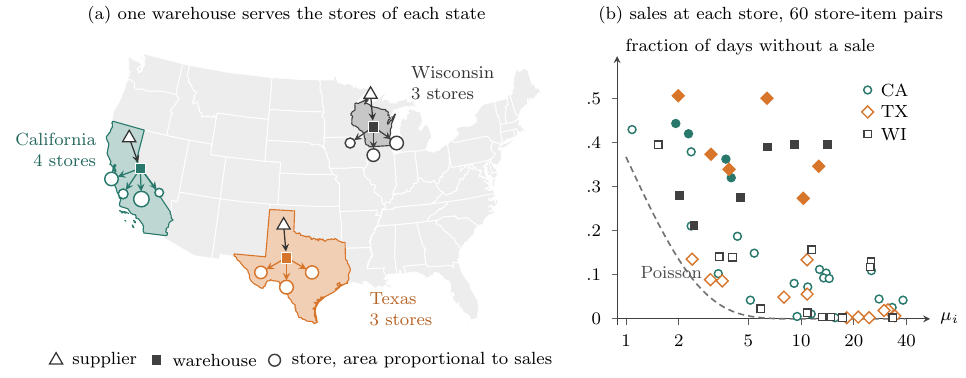}
\caption{The three store networks of the M5 data and their daily sales.}
\label{fig:m5}
\begin{minipage}{.98\linewidth}
{\SingleSpaced\footnotesize\emph{Notes.} Panel (a) shows the supplier, the warehouse, and the stores of each state's network. Positions within each state are illustrative because the data identify stores only by state; each store marker's area is proportional to the store's total mean daily unit sales of the six items selected for its state. Panel (b) plots, for each store and selected item over the last 730 days, the fraction of days without a sale against mean daily sales $\mu_i$ on a logarithmic scale. Filled markers show the five items whose every store records no sale on at least a fifth of days. The dashed curve is $e^{-\mu_i}$, the probability of zero demand under a Poisson demand with the same mean $\mu_i$.}
\end{minipage}
\end{figure}

The M5 forecasting competition data contain 1,941 days of Walmart unit sales and weekly prices at ten stores: four in California, three in Texas, and three in Wisconsin \citep{MakridakisSpiliotisAssimakopoulos2022}. One warehouse serves each state's stores, as in Figure~\ref{fig:m5}(a). In each state we select two items in each of three sales tiers, giving 18 state-item pairs (Appendix~\ref{app:numerics}). We construct stationary test instances by drawing each store's daily demand from the empirical distribution of its daily sales over the last 730 days, independently across stores and days, although Theorem~\ref{thm:main} also covers demand that is correlated across stores within a day. We use daily sales as demand without correcting for censoring by stockouts: the guarantee applies to every demand distribution, and these empirical distributions provide test instances drawn from retail data. These instances also include intermittent demand: in Figure~\ref{fig:m5}(b), every store-item pair has more days without any sale than a Poisson demand with the same mean would have, and some pairs have no sale on half of all days.

We calibrate three cost components. First, annual holding costs at the stores and the warehouse are 25 percent of unit value. Store values are average M5 prices; we assume warehouse value is $3/4$ of their mean, giving $\min_i H_i/h\simeq1.26$. Second, we calibrate the common penalty by the nominal service target $p/(p+\bar H)\in\{0.90,0.97\}$, where $\bar H$ is the mean store holding rate. Third, we set $K_i=H_i\mu_iT_s^2/2$ and $K_0=h\mu T_w^2/2$, with mean daily store demand $\mu_i$ and aggregate demand $\mu$, so that the economic order quantity (EOQ) cycles are $T_s\in\{3,7\}$ days for stores and $T_w\in\{7,14,28\}$ days for the warehouse. These cycles serve only to calibrate setup costs; our $(s,S)$ policies order when inventory positions reach reorder levels, not at fixed intervals. The 18 state-item pairs, two service targets, six EOQ-cycle pairs $(T_s,T_w)$, supplier lead times $L_0\in\{1,3,7,14\}$, and retailer lead times $L_i=1$ give 864 instances. Each satisfies the conditions for the $5/2$ specialization of Theorem~\ref{thm:main}: $\min_iH_i/h>1$, so $\bar H>h$; $p>\bar H$ because both service targets exceed $1/2$; and $K_{\max}/K_0<0.69<3/2$.

Single-location calculations evaluate, without simulation, convolutions of the integer-valued demand distributions and their renewal measures. Appendix~\ref{app:numerics-methods} explains why integer enumeration recovers the single-location optima of the divisible-quantity model and how the search range is chosen to contain them. The implementation uses half width $\lceil \hat w/2\rceil$ and integer shifts, with each candidate's augmented cost evaluated for that implemented policy. Proposition~\ref{prop:scalar} is not proved for these candidates, but the physical upper bound~\eqref{eq:upper} of Proposition~\ref{prop:upper} holds whenever each retailer uses an $(s,S)$ or base-stock policy with a fixed mirror center. Dividing~\eqref{eq:lower} in Proposition~\ref{prop:lower} by $\gamma$ and writing $\theta=(\gamma-1)/\gamma$ gives $\mathrm{LB}(\theta):=\theta\sum_iJ_i(K_i/\theta)+{(1-\theta)}J_0({K_0/(1-\theta)})\le\OPT$ for $0<\theta<1$, and~\eqref{eq:upper} bounds the implemented policy's cost by $\mathrm{UB}:=J_0({K_0+K_{\max}})+\sum_i\cL_i$, where $J_i$ and $J_0$ are single-location values in~\eqref{eq:retailer-value} and~\eqref{eq:source-value}, and $\cL_i$ is the augmented cost~\eqref{eq:local-upper} of retailer $i$'s chosen candidate. Maximizing $\mathrm{LB}(\theta)$ over 23 values of $\theta$ in $[0.02,0.95]$ gives $\mathrm{LB}^*$, whereas the $5/2$ guarantee uses $\mathrm{LB}(3/5)$; thus, $\mathrm{UB}/\mathrm{LB}^*$ is an instance-specific guarantee for the implemented policy, independently of the exact-half-window comparison in Proposition~\ref{prop:scalar}.

For each policy, simulated cost $C_{\rm sim}$ is its average daily cost over 100,000 days after 10,000 warm-up days, with common demand paths across policies within each experiment.\footnote{The replication files retain deterministic seeds, their instance and stream labels, and the generator version.} We report nominal 95 percent Student's $t$ intervals, $C_{\rm sim}$ plus or minus $t_{9,0.975}\simeq2.262$ times the standard error of ten batch means, each over 10,000 consecutive days; coverage is approximate because these means need not be independent and normal. For our policy, this margin has mean $0.20$ and maximum $1.68$ percent of $C_{\rm sim}$ over the 864 instances. Rankings compare point estimates; for small cost differences between our policy and a benchmark, the paired batch-mean interval can include zero.

\subsection{Results with a Common Penalty}

\begin{figure}[htbp]
\SingleSpaced
\centering
\includegraphics[width=.99\linewidth]{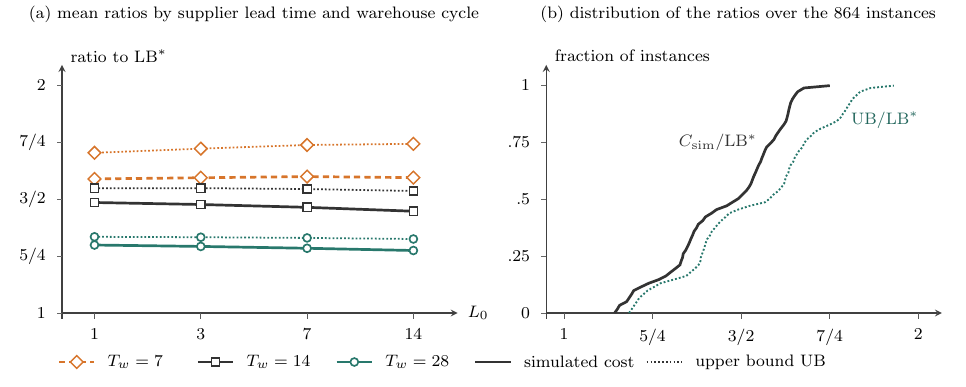}
\caption{Simulated policy costs and analytical upper bounds relative to the best evaluated lower bound.}
\label{fig:numerics}
\begin{minipage}{.98\linewidth}
{\SingleSpaced\footnotesize\emph{Notes.} Panel (a) averages across items, service targets, and calibrating store EOQ cycles $T_s$ for each supplier lead time $L_0$ and warehouse EOQ cycle $T_w$. Heavy solid or dashed curves show $C_{\rm sim}/\mathrm{LB}^*$; dotted curves show $\mathrm{UB}/\mathrm{LB}^*$. Panel (b) shows both empirical cumulative distributions across all 864 instances. Section~\ref{sec:numerics-instances} defines $\mathrm{UB}$ and $\mathrm{LB}^*$.}
\end{minipage}
\end{figure}

Figure~\ref{fig:numerics} reports simulated costs and analytical upper bounds relative to $\mathrm{LB}^*$, and Table~\ref{tab:numerics-main} gives them by group. Since $\mathrm{LB}^*\le\OPT$, $\mathrm{UB}/\mathrm{LB}^*$ bounds the policy's ratio to $\OPT$; it is below $2$ in every instance, against the factor $5/2$ of Theorem~\ref{thm:main}. Across the 864 instances, $C_{\rm sim}/\mathrm{LB}^*$ has mean $1.45$, median $1.49$, and range $1.14$ to $1.75$. The ratio $\mathrm{UB}/\mathrm{LB}^*$ averages $1.53$ and is below $1.94$; $\mathrm{UB}/C_{\rm sim}$ ranges from $1.01$ to $1.14$, averaging $1.06$. The proof of the $5/2$ guarantee compares $\mathrm{UB}$ with $\mathrm{LB}(3/5)$; $\mathrm{UB}/\mathrm{LB}(3/5)$ is below $1.94$, while $\mathrm{LB}(3/5)/\mathrm{LB}^*$ averages $0.91$ and is at least $0.74$.

As $L_0$ rises from $1$ to $14$, mean ratios change little: $C_{\rm sim}/\mathrm{LB}^*$ goes from $1.46$ to $1.44$ and $\mathrm{UB}/\mathrm{LB}^*$ from $1.53$ to $1.54$, consistent with the guarantee, whose two bounds share the supplier cost $f_{L_0}$ at each $L_0$. Their maxima rise from $1.65$ to $1.75$ and from $1.82$ to $1.93$. The setup calibration matters more: the mean simulated-cost ratio ranges from $1.21$ at $(T_s,T_w)=(3,28)$, where store setups are smallest relative to the supplier's, to $1.64$ at $(7,7)$. Raising the service target from $0.90$ to $0.97$ raises this mean from $1.43$ to $1.47$; state and tier means are $1.44$ to $1.46$ and $1.42$ to $1.47$.

Every retailer selects a half window. For cycles planned at setup cost $5K_i/3$, the full window's augmented cost, evaluated at actual setup cost $K_i$, exceeds the half window's by $5$ to $17$ percent, and using only full windows raises simulated cost by $5.1$ percent on average and up to $10.1$ percent. The planning setup cost $5K_i/3$ is also nearly the best in simulation: on the 72 instances with service target $0.90$ and $(T_s,T_w)=(3,7)$, six planning setup costs from $K_i$ to $3K_i$ give simulated costs within $1.5$ percent of one another, and $5K_i/3$ exceeds the cheapest by less than $0.09$ percent.

The uncoordinated $(s,S)$ policy is also an echelon stock $(s,S)$ policy with FIFO service, but it plans single-location $(s,S)$ policies at actual setup costs, $K_i$ for each store and $K_0$ with supplier cost $f_{L_0}$ for the warehouse on its centered echelon position, so the policies differ only in parameters. Our policy has the lower estimated cost in all 864 instances, by $0.07$ to $2.9$ percent.

We next compare our policy, which requires only single-location calculations, with two benchmarks tuned by simulation search for each instance on a separate 50,000-day path and evaluated on the same 100,000-day path as our policy. The first is a periodic-review echelon $(R,nQ)$ adaptation of the batch-transfer structure of \citet{ChenZheng1997}, which dispatches complete queued batches in FIFO order. The second is an integer-ratio policy \citep{ChuShen2010,ShangTaoZhou2015}: the warehouse orders up to an echelon level every $T_0$ periods and store $i$ up to its level every $T_i$ periods, with each $T_i$ an integer multiple or divisor of $T_0$, store orders timed to warehouse receipts, and FIFO service with partial fills. We report the best policies found by the searches specified in Appendix~\ref{app:numerics-methods}; these searches do not establish optimality within either policy class.

On the 216 instances with service target $0.90$ and $(T_s,T_w)$ equal to $(3,7)$, $(3,28)$, or $(7,7)$, our policy's estimated cost exceeds the $(R,nQ)$ benchmark's by $8.2$ percent on average ($2.1$ to $17.3$ percent) and the integer-ratio benchmark's by $12.2$ percent on average (the range runs from our policy being $8.3$ percent cheaper to $45.6$ percent costlier). The cost differences depend on the setup calibration. At $(3,28)$, where supplier setups are largest relative to store setups, these means fall to $3.4$ and $0.8$ percent; at $(7,7)$, where store setups are expensive relative to the supplier's, they rise to $12.9$ and $26.8$ percent. Our policy's estimated cost is lower than the integer-ratio benchmark's in 31 instances, all for five items whose every store records no sale on at least a fifth of days (filled markers in Figure~\ref{fig:m5}(b)): a fixed interval orders at every review with any demand since the last order, whereas our cycles wait until demand has accumulated. The same benchmarks also help assess the lower bound. Both are feasible, so their simulated costs, averaging $1.33\mathrm{LB}^*$ and $1.29\mathrm{LB}^*$, indicate that the lower bound is on average within these factors of $\OPT$; Section~\ref{sec:numerics-exact} checks this bound against bounded-state dynamic-programming values on two-store instances.

Appendix~\ref{app:numerics-checks} reports three further checks. First, varying retailer lead times from $0$ to $2$ gives mean ratios $C_{\rm sim}/\mathrm{LB}^*$ of $1.48$ to $1.55$. Second, a national network of all ten stores gives $1.49$ to $1.56$. Third, on the next year's actual sales, the policy costs on average $1.03$ times its design-period upper bound, which assumes i.i.d.\ demand; item-level ratios from $0.79$ to $1.38$ mainly track changes in demand between the two years, and its realized cost averages $1.07$ times that of a hindsight-refitted policy in the same class (range $1.02$ to $1.18$), with the original cost parameters held fixed.

\subsection{Dynamic Programming for Two Stores}\label{sec:numerics-exact}

For 144 two-store instances built from the six low-tier items, we compute bounded-state dynamic-programming values $C_{\rm DP}$ as numerical benchmarks for $\OPT$; Appendix~\ref{app:numerics-dp} describes the instances, the computation, and its checks. On the same demand paths, the ratio of our policy's simulated cost to $C_{\rm DP}$ averages $1.22$ and ranges from $1.08$ to $1.45$, while $\mathrm{LB}^*$ averages $0.82$ times $C_{\rm DP}$ and is at least $0.71$ times it. The simulated-cost ratio to $C_{\rm DP}$ rises with store setup costs relative to the supplier's, from a mean of $1.10$ at $(T_s,T_w)=(3,28)$ to $1.38$ at $(7,7)$. The dispatch patterns help explain this increase: at $(7,7)$, the computed DP policies send 74 to 97 percent of store shipments in periods with supplier receipts, whereas our policy keeps warehouse stock to serve half-window store cycles, and its store setup costs are $1.75$ to $1.93$ times those of the DP policies.

\section{Concluding Remarks}\label{sec:conclusion}

We construct a deterministic FIFO policy with a $5/2$ guarantee for periodic-review OWMR systems with fixed costs at both echelons and arbitrary finite lead times. Demand vectors are i.i.d.\ over periods with finite positive means and may have dependent components across stores. The retailer holding rates and common backlog penalty are at least the warehouse holding rate, and each retailer setup cost is at most $3/2$ times the supplier's. For arbitrary setup costs and associated retailer demands, the consolidation policy of Theorem~\ref{thm:arbitrary} achieves a factor of $93/19<4.9$, which falls toward $13/3$ as store holding rates rise above the warehouse rate. Both guarantees compare physical cost with the optimal cost over all nonanticipative policies and hold from every valid initial state.

Under the FIFO policy, retailers select replenishment windows using their own demand and lead time, while supplier orders follow one aggregate inventory position, i.e., the echelon inventory position less the constant sum of the retailers' mirror centers. This separation permits local planning while charging every portion of a request: the first portion to that request's setup cost, and each continuation to the extra $K_{\max}$ in the setup cost of the supplier order that enables it. The mirror cost prices the warehouse imbalance created by each retailer's batching. Implementation requires one optimal single-location solve per retailer and a comparison of its full and half windows, which also yields the selected window's mirror center; one supplier $(s,S)$ rule planned at setup cost $K_0+K_{\max}$; and FIFO dispatches with partial fills and a fixed order for same-period requests.

Higher store holding rates sharpen the FIFO guarantee: when they are at least $5/3$ times the warehouse rate and each retailer setup cost is at most $7/5$ times the supplier's, Theorem~\ref{thm:main} gives $12/5$, and it approaches $9/4$ as this ratio grows. If each store's positive demand has a continuous distribution with a decreasing failure rate, retuning the retailer windows lowers the factor to about $2.335$ when each retailer setup cost is at most $1.335$ times the supplier's (Theorem~\ref{thm:dfr}), but this class excludes integer-valued demand. These refinements lead us to conjecture a deterministic $12/5$ guarantee for all demand distributions when each retailer setup cost is at most the supplier's. A separate open question is whether these bounds, including $5/2$ and $93/19$, are tight.

Another extension concerns retailer-specific penalties. There, the warehouse must also decide whom to serve first in a shortage, and our bounds no longer price every shortage at one rate. A guarantee uniform in penalties and setup costs remains open, even for two penalty classes.

On retail data, the integer-rounded FIFO implementation's simulated cost averages $1.45$ times the best evaluated lower bound across 864 instances, and without simulation search its estimated cost comes within $8.2$ and $12.2$ percent on average of tuned echelon and integer-ratio benchmarks. The two-store dynamic-programming benchmarks show that jointly chosen replenishment decisions can reduce store setups and warehouse holding, particularly when store setup costs are high relative to the supplier's. These comparisons quantify the cost of using local planning and FIFO fulfillment. The experiments use stationary, known demand distributions and a common penalty; in all of them, every retailer selected the half window. Together, these comparisons suggest several directions for further work: combining such a guarantee with simulation search, testing larger retailer setups, where the full window can be chosen and consolidation applies, and extending the guarantees to integer dispatch quantities, demand correlated over time, and demand learning.

\label{end:manuscript}\clearpage
\renewcommand\normalsizeXIyesCE{\OriginalManuscriptNormalSize\fontsize{11}{16.3}\selectfont\abovedisplayskip=5pt plus 2pt minus 1pt\belowdisplayskip=5pt plus 2pt minus 1pt\abovedisplayshortskip=2pt plus 1pt minus 1pt\belowdisplayshortskip=4pt plus 1pt minus 1pt}
\renewcommand\footnotesizeIXyesCE{\OriginalManuscriptFootnoteSize\fontsize{9}{11.7}\selectfont}\renewcommand\smallXyesCE{\OriginalManuscriptSmallSize\fontsize{10}{13}\selectfont}\def\bibfont{\fontsize{10}{12.5}\selectfont}
\normalsize\raggedbottom
\predisplaypenalty=10000 %
\setcounter{page}{1}
\renewcommand\thepage{ec\arabic{page}}
\setcounter{footnote}{0}
\renewcommand\thefootnote{E.\arabic{footnote}}
\renewcommand\thesection{\Alph{section}}
\renewcommand\thefigure{\Alph{section}.\arabic{figure}}\setcounter{figure}{0}
\renewcommand\thetable{\Alph{section}.\arabic{table}}\setcounter{table}{0}
\renewcommand\theequation{\Alph{section}.\arabic{equation}}
\renewcommand\thedefinition{\Alph{section}.\arabic{definition}}
\renewcommand\thelemma{\Alph{section}.\arabic{lemma}}
\renewcommand\theproposition{\Alph{section}.\arabic{proposition}}
\renewcommand\thetheorem{\Alph{section}.\arabic{theorem}}
\renewcommand\thecorollary{\Alph{section}.\arabic{corollary}}
\renewcommand\theremark{\Alph{section}.\arabic{remark}}
\makeatletter
\@addtoreset{figure}{section}\@addtoreset{table}{section}\@addtoreset{corollary}{section}\@addtoreset{equation}{section}\@addtoreset{definition}{section}\@addtoreset{lemma}{section}\@addtoreset{proposition}{section}\@addtoreset{theorem}{section}\@addtoreset{remark}{section}
\makeatother
\setcounter{section}{0}\setcounter{equation}{0}\setcounter{definition}{0}\setcounter{lemma}{0}\setcounter{proposition}{0}\setcounter{theorem}{0}\setcounter{corollary}{0}\setcounter{remark}{0}
\setlength{\jot}{2pt}
\begin{center}
\textbf{\large Online Appendices to\\``A $5/2$ Guarantee for Echelon Stock $(s,S)$ Policies\\in One-Warehouse Multi-Retailer Systems''}
\end{center}

\noindent Appendix~\ref{app:scalar} proves the single-location Lemma~\ref{lem:scalar}; Section~\ref{sec:assembly} includes the holding-rate refinement in the proof of Theorem~\ref{thm:main}. Appendix~\ref{app:initialization} gives the initialization and random-horizon arguments for Proposition~\ref{prop:upper}. Appendix~\ref{app:arbitrary} proves the holding-rate-dependent guarantee of Theorem~\ref{thm:arbitrary} for arbitrary setup costs under $\pi_{\rm con}^{(m)}$, including $93/19$ at $m=1$. Appendix~\ref{app:smooth} develops the DFR extension, and Appendix~\ref{app:numerics} gives the items, results by group, and further checks of the numerical study. All appendices use the main text's model and notation, and their equations are numbered by appendix.

\section{Single-Location Inventory Problems}\label{app:scalar}

\proof{Proof of Lemma~\ref{lem:scalar}.} Recall that $f$ is nonnegative, convex, coercive, and globally Lipschitz and that $0<\E D<\infty$. We first establish that the average-cost optimality equation~\eqref{eq:scalar-average} has a solution with a bounded relative value function. Discounting gives a contraction, so each discount factor has a well-defined value function. Subtracting its minimum yields relative values whose bounds and continuity estimates are uniform in the discount factor, allowing a subsequential limit as that factor approaches one. We then use the bounded relative value function to control the endpoint terms in a telescoping sum, which verifies optimality and gives the finite-horizon lower bound. Finally, we identify the canonical thresholds $\hat s$ and $\hat S$ and refine the relative value function to obtain~\eqref{eq:account}, treating zero setup cost at the end. Steps 1 to 5 below carry out these parts in order.

\emph{Step 1. Solve the discounted reset problem.} Assume first that $K>0$. The original problem is the single-location problem of Section~\ref{sec:scalar}, in which an order can only raise the position. The reset problem relaxes it by letting a setup move the position to any real number, also below its current value; we call such a move a downward reset. This relaxation makes the reset branch constant across positions and allows a second process to copy the first process's reset target in the coupling below. For $0<\delta<1$ and a bounded continuous function $v$, define
\[
 (\mathcal B_\delta v)(x):=\min\Bigl\{f(x)+\delta\E v(x-D),\ K+\min_{y\in\R}\bigl[f(y)+\delta\E v(y-D)\bigr]\Bigr\},
\]
where the first term is the cost of placing no order and facing one demand, and the second that of a reset before demand is realized. Since $v$ is bounded and continuous, bounded convergence makes $y\mapsto\E v(y-D)$ bounded and continuous; as $f$ is also continuous and coercive, the inner minimum over $y$ is attained. The operator maps bounded continuous functions into themselves and is a contraction with modulus $\delta$. Let $v_\delta$ be its unique bounded continuous fixed point, so $v_\delta(x)=(\mathcal B_\delta v_\delta)(x)$ for every $x$. Since $\mathcal B_\delta$ preserves nonnegativity, iteration from zero gives $v_\delta\ge0$. Iterating the fixed-point inequality along any policy and using $\delta^T\sup|v_\delta|\to0$ shows that $v_\delta(x)$ is at most that policy's expected discounted cost from $x$. A stationary policy that takes a minimizing branch, resetting to a fixed minimizer of $f(y)+\delta\E v_\delta(y-D)$, attains equality. Thus, $v_\delta(x)$ is the minimum expected discounted cost in the reset problem when the first period starts at position $x$. The policy constructed after~\eqref{eq:scalar-average} resets only upward, which the original problem allows.

\emph{Step 2. Normalize and pass to average cost.} Write $\bar v_\delta(y):=f(y)+\delta\E v_\delta(y-D)$ for the cost of skipping the first order at $y$ and then acting optimally from the next position $y-D$. Since $v_\delta=\mathcal B_\delta v_\delta$, $v_\delta(y)$ is the smaller of $\bar v_\delta(y)$ and the reset cost $K+\min\bar v_\delta$, which is the same at every position. Since $v_\delta\ge0$ and $f$ is coercive, $\bar v_\delta\ge f$ exceeds the reset cost outside some bounded closed interval, where $v_\delta$ then equals the reset cost, its largest value. Hence, the continuous function $v_\delta$ attains its minimum over $\R$ on that interval: $a_\delta:=\min_{y\in\R}v_\delta(y)$ is attained at some $S_\delta$. For each position $y$, set $u_\delta(y):=v_\delta(y)-a_\delta$, the excess cost of starting at $y$ over the minimum $a_\delta$. Set $J_\delta:=(1-\delta)a_\delta$.\footnote{Equivalently, $J_\delta=a_\delta/\bigl(\sum_{t=0}^{\infty}\delta^t\bigr)$ is the constant per-period cost whose discounted total is $a_\delta$.}

Since $\bar v_\delta\ge v_\delta\ge a_\delta$ and $K>0$, the reset cost satisfies $K+\min_y \bar v_\delta(y)\ge K+a_\delta>a_\delta=v_\delta(S_\delta)$. Thus, not ordering at $S_\delta$ attains $a_\delta$, so $\min_y \bar v_\delta(y)=a_\delta$ and the reset cost equals $K+a_\delta$. We obtain the following bounds and discounted optimality equation.\footnote{\label{fn:discounted-cost-bounds}The discounted optimality equation follows by subtracting $\delta a_\delta$ from $v_\delta=\mathcal B_\delta v_\delta$. Since $a_\delta\le v_\delta\le a_\delta+K$ and not ordering at $S_\delta$ attains $a_\delta$, we have $a_\delta=f(S_\delta)+\delta\E v_\delta(S_\delta-D)\ge\min f+\delta a_\delta$ and, at a minimizer $x^*$ of $f$, $a_\delta\le v_\delta(x^*)\le f(x^*)+\delta\E v_\delta(x^*-D)\le\min f+\delta(a_\delta+K)$. Rearranging with $J_\delta=(1-\delta)a_\delta$ gives $\min f\le J_\delta\le\min f+\delta K\le\min f+K$.}
\begin{gather*}
 0\le u_\delta\le K,\qquad \min f\le J_\delta\le\min f+K,\\
 J_\delta+u_\delta(x)=\min\bigl\{f(x)+\delta\E u_\delta(x-D),\,J_\delta+K\bigr\},\qquad u_\delta(S_\delta)=0.
\end{gather*}
If not ordering is optimal at $x$, the no-order branch $f(x)+\delta\E u_\delta(x-D)$ is at most the reset branch $J_\delta+K$, and its continuation term $\delta\E u_\delta(x-D)$ is nonnegative, so $f(x)\le J_\delta+K\le\min f+2K$, where the second inequality uses $J_\delta\le\min f+K$ from footnote~\ref{fn:discounted-cost-bounds}. Because $f$ is convex and coercive, the sublevel set $\{x:f(x)\le\min f+2K\}$ is a bounded interval $[\underline x,\bar x]$, which contains the no-order states of every discounted problem and does not depend on the discount factor $\delta$.

To obtain a continuity bound uniform in the discount factor, choose $d>0$ with $\Prob(D\ge d)>0$. Since the no-order states lie in $[\underline x,\bar x]$, an optimal discounted policy resets at once from any state outside it. From inside, a process that does not order leaves the interval after $\lfloor(\bar x-\underline x)/d\rfloor+1$ demands of at least $d$, and the policy then resets. Its first reset time is, therefore, dominated by a negative-binomial waiting time, uniformly in $\delta$ and in the initial state; let $\bar\tau$ bound its mean. For any initial states $x$ and $y$, couple the processes started there through common demands, and let both follow the optimal reset decisions and targets of the process from $x$. Their positions differ by $x-y$ until the first reset and agree afterward, so, with $\operatorname{Lip}(f)$ the Lipschitz constant of $f$, their costs differ by at most $\operatorname{Lip}(f)|x-y|$ per period before that reset. This comparison and its reversal, with $x$ and $y$ exchanged, give $|u_\delta(x)-u_\delta(y)|\le\operatorname{Lip}(f)\bar\tau|x-y|$, uniformly in $\delta$.

Thus, the functions $u_\delta$ are bounded and equicontinuous, and $S_\delta\in[\underline x,\bar x]$ because not ordering at $S_\delta$ is optimal. The Arzel\`a--Ascoli theorem gives, on each interval $[-m,m]$, a subsequence $\delta\to1$ along which $u_\delta$ converges uniformly. Choosing these subsequences for $m=1,2,\ldots$, each inside the previous one, and taking the $m$th term of the $m$th subsequence (a diagonal extraction) gives one sequence $\delta\to1$ along which $u_\delta\to u$ locally uniformly; since $\min f\le J_\delta\le\min f+K$ by footnote~\ref{fn:discounted-cost-bounds} and $S_\delta\in[\underline x,\bar x]$, a further subsequence also gives $J_\delta\to\bar J$ and $S_\delta\to \hat S$. Bounded convergence in the preceding discounted optimality equation gives the optimality equation below, and equicontinuity gives $u(\hat S)=0$:
\begin{equation}\label{eq:scalar-average}
 \bar J+u(x)=\min\bigl\{f(x)+\E u(x-D),\,\bar J+K\bigr\},\qquad 0\le u\le K,\qquad u(\hat S)=0.
\end{equation}
At $\hat S$, $u(\hat S)=0<K$, so the no-order branch attains the minimum in~\eqref{eq:scalar-average}: $f(\hat S)+\E u(\hat S-D)=\bar J$. Since the continuation term $\E u(\hat S-D)$ is nonnegative because $u\ge0$, this equality gives $f(\hat S)\le\bar J$.

\emph{Step 3. Verify optimality and the finite-horizon bound.} Start at $\hat S$, reset to $\hat S$ whenever $u=K$ at the pre-action position, and otherwise do not order. Pre-action positions never exceed $\hat S$, and $u(\hat S)=0<K$, so every reset raises the position and the policy is feasible for the original problem. At each post-action position, $u<K$, so the no-order branch of~\eqref{eq:scalar-average} holds with equality, and each reset lowers $u$ by exactly $K$. Telescoping $u\in[0,K]$ shows that the policy's expected $T$-period cost differs from $T\bar J$ by at most $K$, so dividing by $T$ shows that its long-run average cost is $\bar J$.

For a lower bound under any nonanticipative policy, let $\mathcal G_t$ contain everything observed before the demand of period $t$ is realized, including network history and auxiliary randomness. Suppose its pre-action position $y_t$, order $q_t\ge0$, and post-action position $x_t=y_t+q_t$ are $\mathcal G_t$-measurable, $D_t$ is independent of $\mathcal G_t$ with the distribution of $D$, and $y_{t+1}=x_t-D_t$. The minimum in~\eqref{eq:scalar-average} is at most its first branch, so $\bar J+u(x)\le f(x)+\E u(x-D)$ for every $x$. By these measurability and independence properties, taking this bound at $x=x_t$ and then expectations gives the first inequality below; since $0\le u\le K$, an order lowers $u$ by at most $K$, which gives the second: $\bar J+\E u(x_t)-\E u(y_{t+1})\le\E f(x_t)$ and $u(y_t)-u(x_t)\le K\ind_{\{q_t>0\}}$. Summing over $t<T$ yields $T\bar J+\E u(y_0)-\E u(y_T)\le\E\sum_{t<T}[K\ind_{\{q_t>0\}}+f(x_t)]$, where the endpoint difference is at least $-K$, and every expectation is well defined because $u$ is bounded and $f\ge0$. Thus, the finite-horizon lower bound holds from every initial position and with the enlarged history $\mathcal G_t$. Comparing it with the constructed policy gives $\bar J=\J_D(f;K)=J$. Hence, that policy is optimal, and the bound $f(\hat S)\le\bar J$ from Step~2 becomes the canonical condition $f(\hat S)\le J$ in~\eqref{eq:canonical} of Lemma~\ref{lem:scalar}.

It remains to attain this value from every initial position. From any prescribed position $x\le \hat S$, an initial order reaches $\hat S$ at cost at most $K$. From $x>\hat S$, wait without ordering until demand brings the position to $\hat S$ or below, then raise it to $\hat S$. Before that crossing, positions lie in $[\hat S,x]$, where $f$ is bounded, and the waiting time is bounded by a negative-binomial waiting time. Thus, preparation has finite expected cost and duration and yields average cost $J$ from any position, so $\J_D(f;K)=J$ does not depend on the initial position, which completes Step~3.

\emph{Step 4. Identify the canonical thresholds.} The constructed policy's no-order states, where $u<K$, lie in $\{x:f(x)\le\min f+2K\}$ by~\eqref{eq:scalar-average} and $J\le\min f+K$, so its cycles have finite mean length, and the renewal-reward theorem \citep[Section~3.6]{Ross1996EC} makes its setup cost per period positive. Since every period also costs at least $\min f$, $J>\min f$. The convexity and coercivity of $f$ then make $\{x:f(x)<J\}=(\hat s,\bar s)$ a nonempty bounded interval with $f(\hat s)=f(\bar s)=J$. For $x\le \hat s$, $f(x)\ge J$, so whichever branch attains the minimum in~\eqref{eq:scalar-average} gives $u(x)\ge\E u(x-D)$. Since $x-D\le \hat s$, iteration gives $u(x)\ge\E u(x-D^{(n)})$, where the cumulative demand $D^{(n)}=D_1+\cdots+D_n$ diverges almost surely. At all sufficiently negative positions, $f>J+K$, so~\eqref{eq:scalar-average} gives $u=K$; bounded convergence then yields $u(x)\ge K$, and $u\le K$ gives $u(x)=K$ for $x\le \hat s$, so ordering is optimal there.

Since $u(\hat S)=0<K$ while $u=K$ on $(-\infty,\hat s]$, we have $\hat S>\hat s$, and $f(\hat S)\le J$ gives $\hat S\le\bar s$ because $f>J$ to the right of $\bar s$. For $\hat s<x<\hat S$, $f(x)<J$ and $u\le K$ give $f(x)+\E u(x-D)<J+K$, so the no-order branch attains the minimum in~\eqref{eq:scalar-average} and $u(x)=f(x)+\E u(x-D)-J<K$. Thus, on states reachable from $\hat S$, the policy orders to $\hat S$ exactly when its pre-action position is at most $\hat s$, does not order on $(\hat s,\hat S]$, and satisfies the canonical conditions~\eqref{eq:canonical}, since $f(\hat s)=J$, $f(\hat S)\le J$, and $f<J$ on $(\hat s,\hat S)$. Post-action positions lie in $(\hat s,\hat S]$. Its cycle length is at most the waiting time for a fixed number of demands of at least $d$, whose geometric tail gives finite moments of every order.

\emph{Step 5. Prove the refined cost inequality and handle zero setup cost.} To add the excess $(f(x)-J)^+$ in~\eqref{eq:account}, which~\eqref{eq:scalar-average} does not give, set $\tilde u:=u-K$, with values in $[-K,0]$, and define $\check u:=0$ on $(-\infty,\hat s]$, $\check u:=\tilde u$ on $(\hat s,\bar s)$, and $\check u:=-K$ on $[\bar s,\infty)$. Then $\osc(\check u)\le K$, and~\eqref{eq:scalar-average} gives $J+\tilde u(x)\le f(x)+\E\tilde u(x-D)$ for every $x$. We verify the refined inequality in three regions. First, at and below $\hat s$, $\check u$ is zero at $x$ and at every position $x-D$ after one demand, and $f(x)\ge J$. Second, inside $(\hat s,\bar s)$, $\check u(x)=\tilde u(x)$ and $\check u(x-D)\ge\tilde u(x-D)$ because $\tilde u\le0$, while $f(x)<J$. Third, at and above $\bar s$, the current relative value is the minimum $-K$ and $f(x)\ge J$. In each case $J+(f(x)-J)^+\le f(x)+\E[\check u(x-D)-\check u(x)]$, so~\eqref{eq:account} holds with $\check u$ in place of $u$, proving the refined inequality.

Finally, if $K=0$, every period costs at least $\min f$, and base stock at a minimizer $x^*$ of $f$ costs exactly $\min f$ in every period once the position has fallen to $x^*$ or below. From $x>x^*$, that takes a wait dominated by a negative-binomial waiting time, as in Step~3, during which the position stays in $[x^*,x]$, where $f\le f(x)$, so the wait adds finite expected cost. Thus, $J=\min f$ is attained from every initial position, and~\eqref{eq:account} and~\eqref{eq:finite-horizon} hold with $u=0$ because $f\ge J$ everywhere.
\Halmos\endproof

\begin{figure}[htbp]
\SingleSpaced
\centering
\includegraphics[width=.76\linewidth]{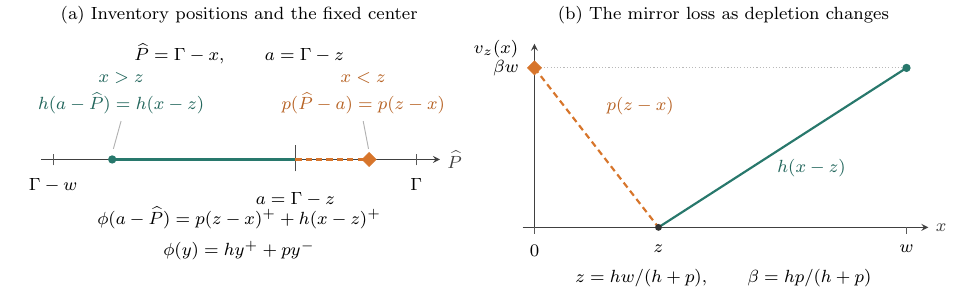}
\caption{Retailer inventory positions and the reversed rates in the warehouse mirror loss.}
\label{fig:mirror}
\begin{minipage}{0.98\linewidth}
\SingleSpaced\footnotesize\emph{Notes.} For $\widehat P=\Gamma-x$ and $a=\Gamma-z$, $v_z(x):=\phi(a-\widehat P)=p(z-x)^++h(x-z)^+$. We draw $p=2h$ and $z=w/3=hw/(h+p)$. The dashed and solid branches have slopes $-p$ and $h$. Both endpoint losses equal $\beta w$, so $B_{w}\le\beta w$.
\end{minipage}
\end{figure}

\begin{figure}[htbp]
\SingleSpaced
\centering
\includegraphics[width=.76\linewidth]{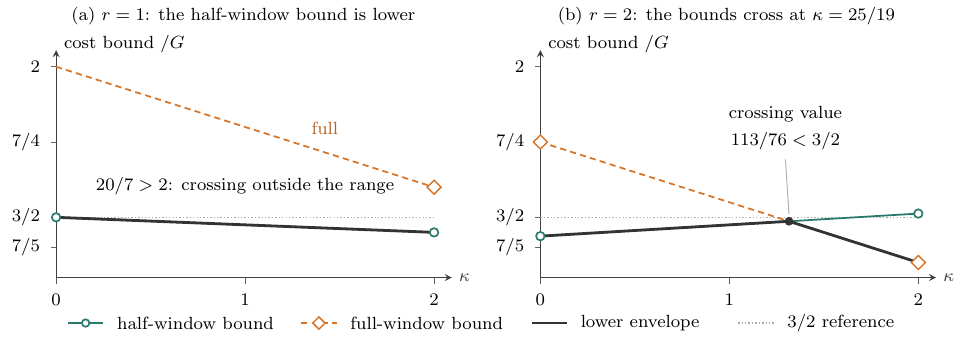}
\caption{The two candidate cost bounds and their lower envelope at $\omega=1$, for $r=1$ and $r=2$.}
\label{fig:comparison}
\begin{minipage}{0.98\linewidth}
\SingleSpaced\footnotesize\emph{Notes.} Common axes show $0<\kappa\le2$. The heavy curve is the lower envelope. Open circles mark the half-window bound and open diamonds the full-window bound. The two bounds cross at $\kappa=20/7>2$ for $r=1$, outside the plotted range, and at $\kappa=25/19$ for $r=2$, where both equal $113/76<3/2$. Values at $\kappa=0$ are limits.
\end{minipage}
\end{figure}

Figure~\ref{fig:comparison} plots the half-window bound~\eqref{eq:half-bound} and the full-window bound~\eqref{eq:full-bound} of Proposition~\ref{prop:scalar} at the maximal normalized width. Since $F\ge0$,~\eqref{eq:canonical-average-cost} gives $k\le GM(2b)$, and $M(2b)\le2M(b)$, so the feasible range is $0<\kappa\le2$. The lower envelope bounds the cheaper augmented cost. The half-window bound is lower when setup cost is a small share of the planning value, while the full-window bound can be lower only when that share is large within that range, which requires $r>25/19$.

\section{Initialization and the Random Horizon}\label{app:initialization}

\proof{Proof of Proposition~\ref{prop:upper}.} We prepare the system, bound its expected $T$-period cost by~\eqref{eq:finite-upper}, and let $T\to\infty$. Bounds uniform in the time left after preparation extend to the original horizon.

\emph{Step 1. Prepare the system and start the policies.} Place no new orders or dispatches until all prescribed pipelines arrive, retaining warehouse stock. This takes finite time and has finite expected cost. Continue withholding supply until every retailer position is at most its selected target $\Gamma_i$. Positions after these arrivals are bounded above by deterministic constants. Hence, after the arrivals, retailer $i$ reaches $\Gamma_i$ or below once its cumulative demand reaches a constant $d_i\ge0$. Count periods from the arrivals. Since $\ind\{\sum_{s<t}D_{i,s}<d_i\}\le e^{d_i}\prod_{s<t}e^{-D_{i,s}}$\footnote{The right side equals $e^{d_i-\sum_{s<t}D_{i,s}}$, which exceeds $1$ when $\sum_{s<t}D_{i,s}<d_i$ and is positive when the indicator is zero. This is the pointwise step in Chernoff's bound for sums of demands that are independent over periods.} and demand vectors are independent over periods, preparation lasts more than $t$ periods with probability at most $\sum_ie^{d_i}(\E e^{-D_i})^t$, where $\E e^{-D_i}<1$ because $\E D_i>0$. Thus, the preparation duration $\tau_{\rm p}$ has finite moments of every order.

The same product bound controls preparation cost using only first demand moments, even with dependence across retailers. For every earlier period $s$ and retailer $j$, independence over periods and $e^{-D_{i,s}}\le1$ give $\E[D_{j,s}\prod_{s'<t}e^{-D_{i,s'}}]\le\E D_j(\E e^{-D_i})^{t-1}$, so the expected demand accumulated on the event that preparation lasts more than $t$ periods is at most a constant times $(t+1)\sum_ie^{d_i}(\E e^{-D_i})^{t-1}$,\footnote{Preparation lasts more than $t$ periods only if some retailer $i$ has $\sum_{s<t}D_{i,s}<d_i$, so the indicator of this event is at most $\sum_ie^{d_i}\prod_{s'<t}e^{-D_{i,s'}}$. Multiplying it by the demand $\sum_jD_{j,s}$ of each of the $t+1$ periods $s\le t$ and applying the preceding bound to each term gives the constant $\sum_j\E D_j$, which is finite and does not depend on $t$.} which is summable over $t$ because $\E e^{-D_i}<1$. This bounds the backlog accumulated during preparation and its cost. Retained warehouse stock and positive retailer inventories are bounded by deterministic constants, so their holding cost is at most a constant times $\E\tau_{\rm p}$. Therefore, preparation has finite expected duration and cost, and its backlogs are integrable.

With preparation complete, if $\Gamma_i-P_i>0$, add one initial request of that size to the queue, so that each inventory position including outstanding requests starts deterministically at $\Gamma_i$. At most $n$ requests are added, with integrable total size. Start the supplier rule on the aggregate position $Y_t=W_t+O_t+\sum_i(P_{i,t}-a_i)$ of~\eqref{eq:aggregate}, with mirror centers $a_i$, and serve the queue. All later dispatches use purchased units or initial stock. The initial requests make the first identity in~\eqref{eq:queue-invariants} hold at the start. Requests, dispatches, receipts, and demand change both of its sides equally, and exhaustive service gives the second identity after each allocation, so both identities hold throughout.

The initial pre-order value of $Y$ is bounded above by the retained warehouse stock plus $\sum_i(\Gamma_i-a_i)$, because $O=0$ and $P_i\le\Gamma_i$. Its negative part is integrable because the initial request sizes $\Gamma_i-P_i$ are. If it is at or below the supplier reorder point, the supplier rule immediately raises $Y$ to its target. Otherwise, the rule leaves $Y$ unchanged. Thus, all post-action values lie in a fixed bounded interval. Any initial excess above the supplier target is removed in finite expected time by the argument that proves $M(w)<\infty$ in Section~\ref{sec:renewal}. Since $\tau_{\rm p}$ is a stopping time, later demands are independent of the preparation history and have their original distribution.

\emph{Step 2. Bound physical costs after preparation.} Let period zero start at the end of preparation. During the first $L_0$ periods, expected costs are finite because the initial backlog and order size are integrable and the delay is fixed. Thereafter, by Section~\ref{sec:queue}, $E$ is bounded $Y$ minus a sum of $L_0$ fresh aggregate demands. Set $\Xi_t:=\sum_i(\widehat P_{i,t}-a_i)$, which is bounded. Since $W_t-\Delta_t=E_t-\Xi_t$ by~\eqref{eq:aggregate} and $W_t\Delta_t=0$ by~\eqref{eq:queue-invariants}, $W_t=(E_t-\Xi_t)^+$ and $\Delta_t=(\Xi_t-E_t)^+$. The bounded positions and nonnegative demand bound warehouse stock above. Finite first moments of fresh lead-time demand give the waiting quantity $\Delta_t$ a uniformly bounded first moment. The lead-time identity~\eqref{eq:retailer-lead} then gives uniformly bounded expected physical inventory costs after the first $L_i$ periods, because $P_{i,t}=\widehat P_{i,t}-\Delta_{i,t}$ by~\eqref{eq:queue-invariants} and $F_i$ is Lipschitz. These bounds use only first moments of fresh demands.

\emph{Step 3. Derive the finite-horizon comparison.} Write $\E_{\rm p}$ for expectation conditional on the preparation history. By~\eqref{eq:retailer-lead}, physical retailer holding and backlog costs through period $T-1$, after their first $L_i$ periods, equal the expected position costs $\E_{\rm p}F_i(P_{i,t})$ over $0\le t<T-L_i$. Since each term is nonnegative, extend the sum to $0\le t<T$, add $hW_t$, and apply Lemma~\ref{lem:running}. After the first $L_0$ periods, the expected supplier cost $\phi(E_t)$ equals $f_{L_0}(Y_{t-L_0})$ by~\eqref{eq:source-expectation}; extend this nonnegative sum to length $T$ as well. Since the prescribed pipelines are now empty and no initial request is partly filled, Lemma~\ref{lem:continuation} applies with $N_{\rm init}=0$ and charges continuations at most $K_{\max}$ per positive supplier order, which with the supplier's own setup cost $K_0$ gives the setup cost $K_0+K_{\max}$ of the supplier problem; the initial requests add at most $nK_{\max}$, one charge $K_{\max}$ per retailer.

To account for the initial periods omitted by these lead-time shifts, let $I_i^0$ be retailer $i$'s net inventory after preparation, $W^0$ the retained warehouse stock, $r_i=\Gamma_i-I_i^0\ge0$ its initial request, and $E^0=W^0+\sum_i(I_i^0-a_i)$ the initial aggregate position before any supplier receipt. The uncovered retailer periods have expected cost at most $H_i(I_i^0)^++p(I_i^0)^-+p(t+1)\E D_i$, and the uncovered supplier terms have expected cost at most $\max\{h,p\}(|E^0|+t\E D)$. Thus, an explicit initialization term is
\[
 c_{\rm init}:=nK_{\max}+\underbrace{\sum_i\sum_{t=0}^{L_i-1}\left[H_i(I_i^0)^++p(I_i^0)^-+p(t+1)\E D_i\right]}_{c_{\rm R}}+\underbrace{\max\{h,p\}\sum_{t=0}^{L_0-1}\left[|E^0|+t\E D\right]}_{c_{\rm S}}.
\]
Empty sums are zero, and the braces mark the retailer sum $c_{\rm R}$ and the supplier sum $c_{\rm S}$. Step 1 bounds $W^0$ and $(I_i^0)^+$ by constants and shows that $(I_i^0)^-$ and $r_i$ have finite means, and $|E^0|\le W^0+\sum_i|\Gamma_i-a_i|+\sum_i r_i$. Hence, $c_{\rm init}$ has a finite mean; it is determined by the preparation history and does not depend on $T$. By~\eqref{eq:retailer-lead}, Lemma~\ref{lem:running}, and~\eqref{eq:source-expectation}, in that order,
\[
 \begin{split}
 \E_{\rm p}\sum_{t<T}\Bigl[hW_t+\sum_i\bigl(H_iI_{i,t}^++pI_{i,t}^-\bigr)\Bigr]
 &\le\E_{\rm p}\sum_{t<T}\Bigl[hW_t+\sum_iF_i(P_{i,t})\Bigr]+c_{\rm R}\\
 &\le\E_{\rm p}\sum_{t<T}\Bigl[\phi(E_t)+\sum_i\bigl(F_i(\widehat P_{i,t})+\phi(a_i-\widehat P_{i,t})\bigr)\Bigr]+c_{\rm R}\\
 &\le\E_{\rm p}\sum_{t<T}\Bigl[f_{L_0}(Y_t)+\sum_i\bigl(F_i(\widehat P_{i,t})+\phi(a_i-\widehat P_{i,t})\bigr)\Bigr]+c_{\rm R}+c_{\rm S}.
 \end{split}
\]
Adding the setup bound of Lemma~\ref{lem:continuation} gives
\begin{equation}\label{eq:finite-upper}
 \E_{\rm p}\sum_{t<T}c_t\le\E_{\rm p}\sum_{t<T}\Bigl[(K_0+K_{\max})\ind_{\{q_{0,t}>0\}}+f_{L_0}(Y_t)+\sum_i\bigl(K_i\ind_{\{\widehat q_{i,t}>0\}}+F_i(\widehat P_{i,t})+\phi(a_i-\widehat P_{i,t})\bigr)\Bigr]+c_{\rm init}.
\end{equation}

\emph{Step 4. Control the remaining random horizon.} We need bounds uniform in the time left after preparation. For $\ell=0,1,\ldots,n$, let $\mathcal C_\ell(m)$ be the expected $m$-period cost of the supplier policy or retailer $\ell$'s policy started at its target, and set $\bar c_0=J_0(K_0+K_{\max})$, the supplier policy's average cost by Lemma~\ref{lem:scalar}, and $\bar c_i=\cL_i$ for $i\ge1$, the limit in~\eqref{eq:local-upper}. Bounded per-period charges and the renewal-reward identity~\eqref{eq:occupation} give $\mathcal C_\ell(m)/m\to\bar c_\ell$. Therefore, $e_\ell(T):=\max_{0\le m\le T}(\mathcal C_\ell(m)-m\bar c_\ell)^+=o(T)$ for $\ell=0,1,\ldots,n$. The retailer policies start at their targets, but the supplier policy may not. From any supplier position in its bounded interval, the first period in which the policy orders up to its target, the reset time $\sigma$, has uniformly bounded expectation, and its per-period cost before $\sigma$ is at most $c_{\rm run}<\infty$. Conditioning at that reset bounds its expected cost over $m\le T$ periods by $m\bar c_0+e_0(T)+c_{\rm tr}$, where $c_{\rm tr}:=c_{\rm run}\sup_y\E_y\sigma+K_0+K_{\max}<\infty$ includes the reset setup. For zero setup cost, use its first passage to the target or below; later post-action positions equal the target.

In the original time index, the first $T$ periods consist of $\min\{\tau_{\rm p},T\}$ preparation periods followed by $m=(T-\tau_{\rm p})^+$ periods of the policy. Apply~\eqref{eq:finite-upper} with $m$ in place of $T$, condition on the preparation history, and use $m\le T$ and the nonnegativity of the limiting costs. If $C_{\rm p}$ is the total preparation cost, then $T^{-1}\E\sum_{t<T}c_t\le\bar c_0+\sum_{i=1}^n\bar c_i+[\E C_{\rm p}+\E c_{\rm init}+c_{\rm tr}+e_0(T)+\sum_{i=1}^n e_i(T)]/T$. Since $\E C_{\rm p}$ and $\E c_{\rm init}$ are finite and each $e_\ell(T)=o(T)$, taking the limit superior proves~\eqref{eq:upper}. These renewal limits use each retailer's marginal demand distribution and the aggregate demand distribution, and require only independence over periods, not across retailers. The preparation, policy choices, and tie rules are deterministic functions of the initial state and demand history.
\Halmos\endproof

\section{The Guarantee for Arbitrary Setup Costs}\label{app:arbitrary}

This appendix proves Theorem~\ref{thm:arbitrary}, so $H_i\ge h$ for every retailer and $p\ge h$ throughout. Demand vectors are i.i.d.\ over periods and associated within each period. Section~\ref{sec:arb-source-budget} proves the required newsvendor comparison under a weaker condition on aggregate demands and applies Lemma~\ref{lem:arb-max} to independent period totals. The proof proceeds in eight steps. Steps 1 and 2 establish the single-location bounds and the allocation identities that express physical cost in terms of supplier, retailer, and reserve costs. Steps 3 to 6 bound these terms for independent uniform offsets and match their coefficients to the shared-cost lower bound~\eqref{eq:lower}. Step 7 selects deterministic offsets using only the observed history, with a summable additional cost. Step 8 extends the guarantee from the clean start to every valid initial state using finite expected preparation and coupling costs.

\subsection{Step 1. Single-Location Foundations}\label{sec:arb-foundations}

Recall the warehouse-rate loss $\phi(x)=hx^++px^-$ and the single-location values $J_i$ and $J_0$ in~\eqref{eq:retailer-value} and~\eqref{eq:source-value}. Write $\Lambda_i:=\min_x F_i(x)$ for the holding and backlog cost of retailer $i$'s best base stock. For an integrable random variable $Y$ and positive rates $a,b$, write the newsvendor cost with overage rate $a$ and underage rate $b$ as $\NV_{a,b}(Y):=\min_z\E[a(z-Y)^++b(Y-z)^+]$.\footnote{The minimum is attained because the expected loss is finite, convex, and coercive in $z$.}

We use three properties of the single-location problem. First, at every positive planning setup cost $k$, Lemma~\ref{lem:scalar} provides a canonical optimum $(\hat s,\hat S)$ with planning value $G=J(k)$ and width $\hat w=\hat S-\hat s$; it satisfies $F(\hat s)=G$ and $F\le G$ on $[\hat s,\hat S]$. Second, for a width-$w$ policy, the renewal-reward identity~\eqref{eq:occupation} gives order frequency $1/M(w)$ and depletion distribution $U(dx)/M(w)$ on $[0,w)$, and Lemma~\ref{lem:renewal} gives mean depletion at most $w/2$ and $M(\hat w)\le2M(\hat w/2)$. Third, by Lemma~\ref{lem:reward}, the canonical savings $g$ in~\eqref{eq:reward} are continuous, concave, and nonnegative on $[0,\hat w]$, with $\int_{[0,\hat w)}g\,dU=k$, and $\beta \hat w\le G$, so $h\hat w\le2G$ because $p\ge h$ gives $\beta\ge h/2$. Retailers with $K_i=0$ use an optimal base-stock policy, at a minimizer of $F_i$, and have optimal value $J_i(0)=\Lambda_i$.

Started at its target, a width-$w$ policy has i.i.d.\ cycles whose length has finite moments of every order, by the argument of Section~\ref{sec:renewal}. The renewal-reward theorem then gives the limits of expected averages of its order counts, holding and backlog costs, and bounded functions of its depletion as expected cycle totals divided by $M(w)$; we call these limits expected renewal averages.

The next two lemmas bound the demand that triggers a request: its mean and the largest demand in a supplier lead-time window. Both enter the reserve bounds of Sections~\ref{sec:arb-source-budget} and~\ref{sec:arb-certificate}.

\begin{lemma}[{\mdseries\scshape Setup Cost and Mean Demand}]\label{lem:arb-mean}
For an $(s,S)$ policy of positive width $w=S-s$ in a retailer's single-location problem, let $\mu:=\E D$ be mean demand, $\nu$ the order frequency, and $\overline F$ the average holding and backlog cost. Then $\overline F\ge\min(H,p)\mu(1-\nu)$. Consequently, if the policy is a canonical optimum with planning value $G=J(k)<k$, then
\begin{equation}\label{eq:arb-mean}
 h\mu\le\frac{G-k\nu}{1-\nu}<G.
\end{equation}
\end{lemma}
\proof{Proof.} We lower-bound cycle cost by the demand between its first and last periods. Conditional Jensen's inequality removes future lead-time demand, and Wald's identity converts the resulting demand span into a bound involving mean demand $\mu$ and order frequency $\nu$.

\emph{Step 1. Reduce to cumulative demand within a cycle.} Start at $S$, and let $\tau=\inf\{j\ge1:D^{(j)}\ge w\}$, with the partial sums $D^{(j)}=D_1+\cdots+D_j$ of Section~\ref{sec:renewal}; a cycle's expected holding and backlog cost $M(w)\overline F$ is the expected sum of $F(S-D^{(j-1)})$ over $j\le\tau$. Let $L$ be the retailer's lead time. The event $\{\tau\ge j\}$ is determined by $D_1,\ldots,D_{j-1}$, so by~\eqref{eq:retailer-value} the next $L+1$ demands can serve as the $L+1$ fresh demands in $F(S-D^{(j-1)})$. Conditional Jensen's inequality over $D_{j+1},\ldots,D_{j+L}$ gives
\[
 \begin{split}
 \E\!\left[\ind_{\{\tau\ge j\}}F(S-D^{(j-1)})\right]
 &=\E\!\left[\ind_{\{\tau\ge j\}}\bigl\{H(S-D^{(j+L)})^++p(S-D^{(j+L)})^-\bigr\}\right]\\
 &\ge\E\!\left[\ind_{\{\tau\ge j\}}\bigl\{H(S-L\mu-D^{(j)})^++p(S-L\mu-D^{(j)})^-\bigr\}\right].
 \end{split}
\]
Summing these nonnegative terms over $j\ge1$ yields $M(w)\overline F\ge\E\sum_{j=1}^{\tau}[H(S-L\mu-D^{(j)})^++p(S-L\mu-D^{(j)})^-]$.

\emph{Step 2. Use the first and last cycle terms.} If $\tau\ge2$, the sum in this bound contains the distinct terms $j=1$ and $j=\tau$. Since $Hx^++px^-\ge\min(H,p)|x|$, the triangle inequality gives
\[
 \begin{split}
 &H(S-L\mu-D_1)^++p(S-L\mu-D_1)^-+H(S-L\mu-D^{(\tau)})^++p(S-L\mu-D^{(\tau)})^-\\
 &\qquad\ge\min(H,p)\bigl(|S-L\mu-D_1|+|S-L\mu-D^{(\tau)}|\bigr)\ge\min(H,p)(D^{(\tau)}-D_1).
 \end{split}
\]
If $\tau=1$, then $D^{(\tau)}=D_1$ and this lower bound is zero. Wald's identity, justified by Tonelli's theorem for nonnegative demands, gives $\E(D^{(\tau)}-D_1)=\mu(M(w)-1)$; dividing by $M(w)=1/\nu$ proves the first claim of the lemma. For the canonical optimum, splitting its average cost gives $\overline F=G-k\nu\ge0$, so $\nu\le G/k<1$, and the first claim and $\min(H,p)\ge h$ give $h\mu(1-\nu)\le G-k\nu$, which is the first inequality in~\eqref{eq:arb-mean}; the second holds because $k>G$ and $\nu>0$ give $G-k\nu<G(1-\nu)$.
\Halmos\endproof

\begin{lemma}[{\mdseries\scshape Independent-Demand Maximum}]\label{lem:arb-max}
Let $V_j$ form a finite nonempty family of independent nonnegative integrable variables, $V_{\max}=\max_jV_j$, $V_{\rm sum}=\sum_jV_j$, $V_*=\max_j\E V_j$, and $\alpha=h/(h+p)$. For any event $B$ with $\Prob(B)\le\alpha$,
\begin{equation}\label{eq:arb-max}
 (h+p)\E[V_{\max}\ind_B]\le(1+h/p)\NV_{h,p}(V_{\rm sum})+hV_*,
\end{equation}
where $B$ may depend on the $V_j$ and on other randomness. For long-run expected averages, if each period has a family with this distribution and $\limsup_{T\to\infty}T^{-1}\sum_{t<T}\Prob(B_t)\le\alpha$, then $\limsup_{T\to\infty}T^{-1}(h+p)\sum_{t<T}\E[V_{\max,t}\ind_{B_t}]\le(1+h/p)\NV_{h,p}(V_{\rm sum})+hV_*$.
\end{lemma}
\proof{Proof.} We express the newsvendor cost as an upper-tail average of the sum $V_{\rm sum}$ and use independence to compare this with the upper tail of the maximum $V_{\max}$. Quantile bounds then control arbitrary events of probability at most $\alpha$ and their long-run averages.

\emph{Step 1. Express the newsvendor cost through an upper-tail average.} Since Steps 2 and 3 use $B$ for threshold events, write $\mathcal B$ for the event $B$ of~\eqref{eq:arb-max}. Define the upper-tail average $\mathcal C_\alpha(Y):=\min_{y\in\R}\{y+\E(Y-y)^+/\alpha\}$. For nonnegative integrable $Y$, the one-sided derivatives in $y$ show that a minimizing $y$ can be chosen nonnegative and satisfies $\Prob(Y>y)\le\alpha\le\Prob(Y\ge y)$; we call such a $y$ an upper-$\alpha$ quantile of $Y$. Since $\alpha=h/(h+p)$, the identity $(z-V_{\rm sum})^+=z-V_{\rm sum}+(V_{\rm sum}-z)^+$ gives $\NV_{h,p}(V_{\rm sum})=h\{\mathcal C_\alpha(V_{\rm sum})-\E V_{\rm sum}\}$. Put $e:=\NV_{h,p}(V_{\rm sum})/h$. We first bound weighted deviations of the sum; Step 3 will choose weights from threshold events of the maximum. For any random weight $\chi$ with $0\le\chi\le1/\alpha$ and $\E\chi=1$, and any $y\in\R$, we have $\E[V_{\rm sum}\chi]=y+\E[(V_{\rm sum}-y)\chi]\le y+\E(V_{\rm sum}-y)^+/\alpha$. Taking the infimum over $y$ yields $\E[V_{\rm sum}\chi]\le\mathcal C_\alpha(V_{\rm sum})$. Since $\E\chi=1$ and $\mathcal C_\alpha(V_{\rm sum})=\E V_{\rm sum}+e$ by the definition of $e$, this gives
\begin{equation}\label{eq:arb-weighted}
 \E[(V_{\rm sum}-\E V_{\rm sum})\chi]\le e\qquad(0\le\chi\le1/\alpha,\ \E\chi=1).
\end{equation}

\emph{Step 2. Relate threshold events to covariance.} Fix $y\ge V_*$. Let $B_j$ contain $\{V_j>y\}$ and a possibly randomized part of $\{V_j=y\}$, decided by a tie coin for $j$; the coins are independent of one another and of the variables $V_j$. Set $r_j:=\Prob(B_j)$, $B:=\bigcup_jB_j$, and $q:=\Prob(B)$. If $q=1$, then $\ind_B=1$ almost surely, so both sides of~\eqref{eq:arb-tail-covariance} are zero. For $q<1$, we have $r_j<1$ for every $j$. Each $B_j$ is determined by $V_j$ and its coin, and these pairs are independent across $j$, so $V_j$ and $B_j$ are independent of $B_{-j}:=\bigcup_{\ell\ne j}B_\ell$. Since $\ind_B=1-\ind_{B_j^c}\ind_{B_{-j}^c}$ and this independence gives $1-q=(1-r_j)\Prob(B_{-j}^c)$,
\[
 \operatorname{Cov}(V_j,\ind_B)
 =-\operatorname{Cov}\bigl(V_j,\ind_{B_j^c}\ind_{B_{-j}^c}\bigr)
 =-\Prob(B_{-j}^c)\operatorname{Cov}(V_j,\ind_{B_j^c})
 =\frac{1-q}{1-r_j}\left\{\E[V_j\ind_{B_j}]-(\E V_j)r_j\right\}.
\]
Since $B_j$ contains $\{V_j>y\}$ and $V_j=y$ on the rest of $B_j$, $V_j\ind_{B_j}=y\ind_{B_j}+(V_j-y)^+$, so $\E[V_j\ind_{B_j}]=yr_j+\E(V_j-y)^+$. Hence, summing the last display over $j$ and using linearity of covariance gives
\[
 \frac{\operatorname{Cov}(V_{\rm sum},\ind_B)}{1-q}
 =\sum_j\frac{(y-\E V_j)r_j+\E(V_j-y)^+}{1-r_j}
 \ge(y-V_*)q+\E(V_{\max}-y)^+
 =\E[(V_{\max}-V_*)\ind_B],
\]
where the inequality uses $y\ge V_*\ge\E V_j$, $1/(1-r_j)\ge1$, the union bound $\sum_jr_j\ge q$, and $\sum_j(V_j-y)^+\ge(V_{\max}-y)^+$; for the last equality, $B$ contains $\{V_{\max}>y\}$, because $V_{\max}>y$ forces some $V_j>y$, and $V_{\max}\ge V_j\ge y$ on each $B_j$, so $(V_{\max}-y)\ind_B=(V_{\max}-y)^+$, and adding $(y-V_*)q=\E[(y-V_*)\ind_B]$ to its expectation gives $\E[(V_{\max}-V_*)\ind_B]$. Thus, multiplying by $1-q$ gives
\begin{equation}\label{eq:arb-tail-covariance}
 \operatorname{Cov}(V_{\rm sum},\ind_B)\ge(1-q)\E[(V_{\max}-V_*)\ind_B].
\end{equation}

\emph{Step 3. Compare the two upper-tail averages.} Fix an upper-$\alpha$ quantile $y_\alpha\ge0$ of $V_{\max}$, so that $\Prob(V_{\max}>y_\alpha)\le\alpha\le\Prob(V_{\max}\ge y_\alpha)$. Two cases, according to whether $y_\alpha\ge V_*$, prove~\eqref{eq:arb-tail-comparison} below.

\emph{Case 1: $y_\alpha\ge V_*$.} Take $y=y_\alpha$, and let each tie coin include $\{V_j=y\}$ in $B_j$ with a common probability $\theta\in[0,1]$. Because the coins are independent, $q=\Prob(B)$ is a polynomial in $\theta$, so it is continuous. At $\theta=0$, $B=\{V_{\max}>y\}$, and at $\theta=1$, $B=\{V_{\max}\ge y\}$, so the quantile property $\Prob(V_{\max}>y_\alpha)\le\alpha\le\Prob(V_{\max}\ge y_\alpha)$ and the intermediate value theorem give some $\theta$ with $q=\alpha$. For this $\theta$, $B$ contains every outcome with $V_{\max}>y$, and $V_{\max}=y$ elsewhere on $B$, so $\E[V_{\max}\ind_B]=y\alpha+\E(V_{\max}-y)^+=\alpha\mathcal C_\alpha(V_{\max})$. The weight $\chi=\ind_B/\alpha$ satisfies $0\le\chi\le1/\alpha$ and $\E\chi=q/\alpha=1$, and~\eqref{eq:arb-tail-covariance} applies because $y=y_\alpha\ge V_*$, so
\[
 e\ge\frac{\E[(V_{\rm sum}-\E V_{\rm sum})\ind_B]}{\alpha}
 =\frac{\operatorname{Cov}(V_{\rm sum},\ind_B)}{\alpha}
 \ge\frac{1-\alpha}{\alpha}\E[(V_{\max}-V_*)\ind_B]
 =(1-\alpha)\{\mathcal C_\alpha(V_{\max})-V_*\},
\]
where the first inequality is~\eqref{eq:arb-weighted}, the first equality is the definition of covariance, the second inequality is~\eqref{eq:arb-tail-covariance} with $q=\alpha$, and the last equality uses $\E[V_{\max}\ind_B]=\alpha\mathcal C_\alpha(V_{\max})$.

\emph{Case 2: $y_\alpha<V_*$.} Take $y=V_*$ and $B_j=\{V_j>V_*\}$, so $q=\Prob(V_{\max}>V_*)\le\alpha$. When $q=0$, $V_{\max}\le V_*$ almost surely, so $\mathcal C_\alpha(V_{\max})\le V_*$, and~\eqref{eq:arb-tail-comparison} follows from $e\ge0$. Otherwise the weight $\chi:=\ind_B/\alpha+(\alpha-q)\ind_{B^c}/[\alpha(1-q)]$ lies in $[0,1/\alpha]$ and has mean one. Applying~\eqref{eq:arb-weighted} to this weight and~\eqref{eq:arb-tail-covariance} at threshold $V_*$ gives $e\ge(1-\alpha)\operatorname{Cov}(V_{\rm sum},\ind_B)/[\alpha(1-q)]\ge(1-\alpha)\E(V_{\max}-V_*)^+/\alpha$. The choice $y=V_*$ in the definition of $\mathcal C_\alpha(V_{\max})$ gives $\mathcal C_\alpha(V_{\max})\le V_*+\alpha^{-1}\E(V_{\max}-V_*)^+$, and the preceding bound on $e$ caps the second term at $e/(1-\alpha)$. Thus, both cases establish
\begin{equation}\label{eq:arb-tail-comparison}
 \mathcal C_\alpha(V_{\max})\le V_*+\frac{e}{1-\alpha}.
\end{equation}

\emph{Step 4. Apply the comparison to events and time averages.} We now apply the upper-tail bound to the event in the lemma: $\E[V_{\max}\ind_{\mathcal B}]\le y_\alpha\Prob(\mathcal B)+\E(V_{\max}-y_\alpha)^+\le\alpha\mathcal C_\alpha(V_{\max})$, where the first inequality takes expectations in the pointwise bound $V_{\max}\ind_{\mathcal B}\le y_\alpha\ind_{\mathcal B}+(V_{\max}-y_\alpha)^+$. To obtain this bound, multiply $V_{\max}\le\max\{V_{\max},y_\alpha\}=y_\alpha+(V_{\max}-y_\alpha)^+$ by $\ind_{\mathcal B}\ge0$, and then bound $(V_{\max}-y_\alpha)^+\ind_{\mathcal B}$ by $(V_{\max}-y_\alpha)^+$. For the second inequality, $y_\alpha\ge0$ and $\Prob(\mathcal B)\le\alpha$ give $y_\alpha\Prob(\mathcal B)\le\alpha y_\alpha$, and $\alpha y_\alpha+\E(V_{\max}-y_\alpha)^+=\alpha\mathcal C_\alpha(V_{\max})$ because the upper-$\alpha$ quantile $y_\alpha$ attains the minimum defining $\mathcal C_\alpha(V_{\max})$. Multiply this bound by $h+p$, bound its right side by~\eqref{eq:arb-tail-comparison} multiplied by $(h+p)\alpha=h$, and use $he=\NV_{h,p}(V_{\rm sum})$ from Step 1 and $1/(1-\alpha)=1+h/p$ to obtain~\eqref{eq:arb-max}.

For the long-run assertion, which needs no joint stationarity, suppose each period $t$ has a family as in the lemma whose distribution does not depend on $t$, so that its maximum $V_{\max,t}$ and its sum both have fixed integrable distributions, and suppose the events $B_t$ have long-run average probability at most $\alpha$. With an upper-$\alpha$ quantile $y\ge0$ of the distribution of $V_{\max,t}$, every period satisfies
\begin{equation}\label{eq:arb-quantile}
 \E[V_{\max,t}\ind_{B_t}]\le y\Prob(B_t)+\E(V_{\max,t}-y)^+.
\end{equation}
Averaging gives the upper-tail bound $\alpha\mathcal C_\alpha(V_{\max,t})$ of the proof's last step, since the average of $\Prob(B_t)$ is at most $\alpha$ and $y$ is a minimizer in $\mathcal C_\alpha(V_{\max,t})$. Since the sum also has a fixed distribution, the right side of~\eqref{eq:arb-tail-comparison} is the same in every period, so~\eqref{eq:arb-max} holds for these averages.\Halmos\endproof

\subsection{Step 2. Policy and Allocation}\label{sec:arb-policy}
Sections~\ref{sec:arb-policy} to~\ref{sec:arb-derandomization} start with every retailer at net inventory and position $\Gamma_i$, nothing in transit, and an empty warehouse and queue. The supplier's first order raises its position to $\hat S_0+a_0$; at zero supplier setup cost, it raises the position to the selected base-stock target. Section~\ref{sec:arb-initial} reduces every valid initial state to this clean start. Fix $m\ge1$ and suppress its subscript in the policy parameters. Set
\begin{equation}\label{eq:arb-parameters}
 e=\frac1m,\quad g=\frac{13}{6}+\frac{16e}{57},\quad
 \zeta=\frac g{g-1},\quad
 B=2(g-1),\quad
 u_*=u_m=\frac{3\zeta(1+e)}{4(4+e\zeta)}.
\end{equation}
The threshold $u_*$ is where the half- and full-cycle bounds of Lemma~\ref{lem:arb-certificate} meet. Its proof verifies that these bounds are at most $BG$ for all $e\in(0,1]$; $e=0$ is used only as a limiting endpoint in that verification. Plan the supplier at setup cost $gK_0$ and retailer $i$ at $k_i=\zeta K_i$, with values $G_0=J_0(gK_0)$ and $G_i=J_i(k_i)$. We write $\pi_{\rm con}$ and $\pi_{\rm random}$ for the deterministic policy $\pi_{\rm con}^{(m)}$ and its uniform-offset counterpart at this fixed $m$, which draws each offset uniformly from $(0,1)$.

For a retailer with $K_i>0$, let $\hat w_i=\hat S_i-\hat s_i$ be the canonical width of its canonical optimum $(\hat s_i,\hat S_i)$ at setup cost $k_i$, and $u_i=k_i/[G_iM_i(\hat w_i)]$ the share of $G_i$ spent on setups, with $M_i$ the function $M$ of~\eqref{eq:renewal} for retailer $i$. Choose a policy at each retailer by the following three rules.
\begin{itemize}
\item If $\Lambda_i+K_i\le BG_i$, use base stock at a minimizer of $F_i$ and dispatch every allocation at once.
\item Otherwise, if $u_i\le u_*$, use the width-$\hat w_i/2$ policy with a maximizing shift in~\eqref{eq:shift}.
\item Otherwise, use the canonical optimum $(\hat s_i,\hat S_i)$, of width $\hat w_i$ and target $\hat S_i$.
\end{itemize}
We call retailers of the first type \emph{immediate} and the others \emph{batch} retailers. The warehouse ships every allocation to an immediate retailer at once, but it holds the stock allocated to each request of a batch retailer until the dispatch thresholds defined below, so the request ships in at most three portions, regardless of how many periods its allocation spans. Retailers with $K_i=0$ are immediate because this gives $\Lambda_i+K_i=J_i(0)=G_i\le BG_i$. Let $\Gamma_i$ and $w_i$ be the selected target and width, so $w_i$ is $0$, $\hat w_i/2$, or $\hat w_i$, the last two corresponding to the \emph{half-cycle} and \emph{full-cycle} policies.

Append positive requests to one FIFO queue with fixed tie-breaking. In each period, allocate unassigned warehouse stock exhaustively to requests in queue order, the last one possibly in part; allocated stock not yet dispatched is the request's \emph{reserve}. When a batch retailer places a request of size $Q$, the request receives $\vartheta$ uniform on $(0,1)$ and independent of all demands and other offsets, and uses the \emph{offset} $\xi=b\vartheta$, where $b=3Q/4$, and \emph{dispatch thresholds}, levels of its cumulative allocation, at $\xi$ and at $\xi+b$ if the latter is below $Q$. Since $b=3Q/4$, these are the thresholds $3Q\vartheta/4$ and $3Q(1+\vartheta)/4$ of Definition~\ref{def:con}, the second present only when $\vartheta<1/3$; Section~\ref{sec:arb-derandomization} replaces this uniform $\vartheta$ by a deterministic choice. Whenever its allocation reaches or passes a new threshold, and when it is fully allocated, dispatch its entire reserve.\footnote{For example, a request with $Q=100$ and $\vartheta=0.2$ has thresholds $3Q\vartheta/4=15$ and $3Q(1+\vartheta)/4=90$. Suppose that, at the head of the queue, it is allocated $10$, $30$, $55$, and $5$ units in four periods, so its cumulative allocation is $10$, $40$, $95$, and $100$. The warehouse holds the first $10$ units, ships $40$ when the cumulative allocation passes $15$, ships $55$ when it passes $90$, and ships the last $5$ when it reaches $100$. If a single receipt instead raised the cumulative allocation from $40$ to $100$, the $60$ units then held would ship together, in a single shipment with a single setup.} Intuitively, the thresholds cap a request's shipments at three, and the random offset keeps its reserve small on average: for every allocation path, the expected reserve is at most $3Q/8$, half the spacing $b$. All shipments to a retailer in a period form one dispatch, counted once; a shipment uses only stock already allocated to its request.

For the supplier, raise both levels of a canonical optimum $(\hat s_0,\hat S_0)$ attaining $G_0$ by $a_0=(-\hat s_0)^+$, so the reorder level is nonnegative. At zero setup cost, use base stock at a nonnegative minimizer of $f_{L_0}$ and $a_0=0$. Apply this policy to the supplier position $W+O+\sum_i(P_i-\Gamma_i)$.

Let $\widehat P_i$ be the post-action inventory position including outstanding requests, $X_i=\Gamma_i-\widehat P_i\ge0$ depletion, $\Delta_i$ unallocated request volume, $R_i$ reserve, $\overline W$ unassigned warehouse stock, and $O$ supplier stock in transit. If $P_i$ is the physical position (net inventory plus stock in transit), then $\widehat P_i=P_i+\Delta_i+R_i$. Write $X=\sum_iX_i$, $\Delta=\sum_i\Delta_i$, and $R=\sum_iR_i$. The \emph{ideal} position $P_i^{\rm id}$, which retailer $i$ would have if every allocation shipped at once, and the other quantities satisfy
\begin{align}
 P_i^{\rm id}&=\widehat P_i-\Delta_i,&
 \overline W\Delta&=0,&W&=\overline W+R,&
 P_i&=P_i^{\rm id}-R_i,\label{eq:arb-physical}\\
 E&=\overline W-\Delta-X,&Y&=E+O,&
 E_{t+L_0}&=Y_t-\sum_{\ell=t}^{t+L_0-1}D_\ell.\label{eq:arb-source}
\end{align}
Thus, $E=W+\sum_i(P_i-\Gamma_i)$, as in~\eqref{eq:aggregate} with centers $\Gamma_i$; retention, which moves $R_i$ from $P_i$ to warehouse stock, leaves $E$ and $Y$ unchanged. Hence, supplier orders, requests, allocations, and the head of the queue are the same as if every allocation shipped at once, independent of the offsets; we also call them ideal. Section~\ref{sec:arb-derandomization} uses this to choose offsets from the ideal state.

After initialization, $Y_t\ge0$ because the shifted reorder level and base-stock target are nonnegative. Since $\overline W_t+O_t=Y_t+\Delta_t+X_t\ge\Delta_t$ and stock in transit arrives by $t+L_0$, all requests pending at $t$ are fully allocated by $t+L_0$. Hence, if the request at the head of the queue at $t$ is unfinished, its retailer $i$ created it in a period $\tau$ from $t-L_0+1$ to $t$. Its size is $X_{i,\tau-1}+D_{i,\tau-1}\le w_i+D_{i,\tau-1}$, because it raises the position of retailer $i$ to $\Gamma_i$ and post-action depletion in period $\tau-1$ is at most $w_i$. We call $D_{i,\tau-1}$ its triggering demand, which Section~\ref{sec:arb-source-budget} later bounds.

The same allocation rule also determines where reserve can be held. Only the unfinished request at the head of the queue can hold a reserve, since completion dispatches it and later requests have no allocation. If the queue is nonempty, exhaustive allocation gives $\overline W=0$ and~\eqref{eq:arb-source} gives $E=-\Delta-X<0$, implying a supplier shortage; hence, a reserve is held only when $E<0$.

We bound the cost $hW+\sum_iF_i(P_i)$ in each period. Since $F_i$ has slopes at least $-p$, $F_i(P_i)\le F_i(\widehat P_i)+p(\Delta_i+R_i)$; also $hW=h\overline W+hR$, and $h\overline W+p\Delta=\phi(\overline W-\Delta)=\phi(E+X)$ because $\overline W\Delta=0$. Hence, this cost is at most
\begin{equation}\label{eq:arb-running}
 \sum_iF_i(\widehat P_i)+\phi(E+X)+(h+p)R.
\end{equation}
To separate depletion from the supplier loss in this bound, use the identity $\phi(E+X)=\phi(E)+hX-(h+p)\min\{X,E^-\}$. It follows by checking the signs of $E$ and $E+X$. On a nonempty queue, $E^-=\Delta+X\ge X$, so the minimum is $X$; on an empty queue, $R=0$ and the subtracted term is nonnegative. Thus,~\eqref{eq:arb-running} is at most
\begin{equation}\label{eq:arb-credited}
 \phi(E)+\sum_i[F_i(\widehat P_i)+hX_i]+(h+p)(R-X)^+.
\end{equation}
This bound holds pathwise in every period; we call $(R-X)^+$, the reserve in excess of depletion, the \emph{credited reserve}. Section~\ref{sec:arb-initial} relates $F_i(P_{i,t})$ to physical retailer costs in period $t+L_i$.

\subsection{Step 3. Reserve After a Dispatch}\label{sec:arb-phase}
Every batch retailer's request has a threshold at $\xi<b<Q$ and a second at $\xi+b$ exactly when $\xi<Q/4$, with probability $(Q/4)/b=1/3$. Charge one shipment at each threshold and one at completion, even if allocation crosses several thresholds in one period. There are three charged shipments if the second threshold exists and two otherwise; the expected number is, therefore, $7/3$ per request. Charging $K_i$ for each covers every physical setup, including when several charged shipments form one dispatch, and has expectation $7K_i/3$ per request for every allocation path.

\begin{lemma}[{\mdseries\scshape Credited Reserve}]\label{lem:arb-phase}
Fix the entire ideal path, which determines an unfinished request's size $Q$ and cumulative allocation $x<Q$, and let $z\ge0$ also be determined by that path. With $b=3Q/4$ and $R$ the request's reserve,
\begin{equation}\label{eq:arb-phase}
 \E_{\xi}(R-z)^+\le\frac{[(b-z)^+]^2}{2b}.
\end{equation}
At completion, when the cumulative allocation $x$ reaches $Q$, the reserve is zero.
\end{lemma}
\proof{Proof.} Conditional on the ideal path, $Q$, $x$, and $z$ are fixed and $\xi$ is uniform on $(0,b)$. Let $A=(x-\xi)\bmod b\in[0,b)$, which is also uniform on $[0,b)$. A comparison rule that ships only the stock up to each crossed threshold leaves reserve $\min\{x,A\}$: it leaves $x$ before the first threshold and the distance to the preceding threshold thereafter. Actual dispatches clear the entire reserve, so $R\le\min\{x,A\}\le A$. Since $x<Q<2b$, thresholds beyond $Q$ do not affect this comparison. Hence,
\[
 \E_\xi(R-z)^+
 \le\frac1b\int_0^b(a-z)^+\,da
 =\frac{[(b-z)^+]^2}{2b}.\tag*{\Halmos}
\]
\endproof

To separate the cycle-width contribution from the demand that triggers a request, write $\delta=3/8$ and $\psi(z)=[(1-4z/3)^+]^2/(8/3)$, which is nonincreasing with $\psi(0)=\delta$. The elementary inequality
\begin{equation}\label{eq:arb-perspective}
 \frac{[(a+e-z)^+]^2}{a+e}\le e+\frac{[(a-z)^+]^2}{a}
 \qquad(a,e,z\ge0)
\end{equation}
uses the value zero whenever a ratio's denominator is zero. This is the continuous extension, since $0\le[(a-z)^+]^2/a\le a$ for $a>0$ and $z\ge0$. For $a>0$, if $z\le a$, the right minus left is $z^2[1/a-1/(a+e)]$; if $z>a$, the left is at most $e$. Taking $a\downarrow0$ gives the boundary case, including $a=e=0$.

We first apply this inequality when the request at the head of the queue is from retailer $i$, which uses the full-cycle policy of width $w_i=\hat w_i$. Let $d$ be the request's triggering demand, as defined in Section~\ref{sec:arb-policy}. The request's size satisfies $Q\le \hat w_i+d$ by that section, so $b=3Q/4\le3(\hat w_i+d)/4$, and
\begin{equation}\label{eq:arb-split}
 \E_{\xi}(R-X)^+
 \le \delta d+\hat w_i\psi(X/\hat w_i)\le \delta d+\hat w_i\psi(X_i/\hat w_i),
\end{equation}
where the first inequality takes three steps. First,~\eqref{eq:arb-phase} with $z=X$, which the ideal path determines, bounds the left side by $[(b-X)^+]^2/(2b)$. Second, this bound is nondecreasing in $b$, because it is zero for $b\le X$ and its derivative $(b^2-X^2)/(2b^2)$ is nonnegative for $b>X$, so $b\le3(\hat w_i+d)/4$ allows replacing $b$ by $3(\hat w_i+d)/4$. Third, dividing~\eqref{eq:arb-perspective} with $(a,e,z)=(3\hat w_i/4,3d/4,X)$ by two bounds the result by $3d/8+[(3\hat w_i/4-X)^+]^2/(3\hat w_i/2)=\delta d+\hat w_i\psi(X/\hat w_i)$. The second inequality holds because $\psi$ is nonincreasing and $X\ge X_i$, as depletions are nonnegative.

For the other batch policy, suppose the request at the head of the queue is from retailer $i$, which uses the half-cycle policy of width $w_i=\hat w_i/2$. With $d$ again denoting the request's triggering demand,~\eqref{eq:arb-phase} with $z=0$ and $Q\le w_i+d$ give the weaker bound $\E_{\xi}(R-X)^+\le\E_{\xi}R\le b/2=\delta Q\le\delta(w_i+d)$.

To combine the two cases, let $B_{i,t}$ be the event that retailer $i$'s request remains unfinished at the head of the queue after period $t$'s allocation. If the request at the head of the queue is from a retailer using a full-cycle or half-cycle policy, let $d_t$ be its triggering demand. Set $d_t=0$ if the queue is empty or that retailer uses base stock. Write $\mathcal P_f$ and $\mathcal P_h$ for the sets of retailers using the full-cycle and half-cycle policies, respectively. Only the unfinished request at the head of the queue can hold a reserve. If its retailer uses base stock, that reserve is zero. Hence,~\eqref{eq:arb-split} for full-cycle requests and the weaker bound $\delta(w_i+d)$ for half-cycle requests give the decomposition
\begin{equation}\label{eq:arb-decomposition}
 \E_\xi(R_t-X_t)^+
 \le\delta d_t
 +\sum_{i\in\mathcal P_f}\hat w_i\psi(X_{i,t}/\hat w_i)\ind_{B_{i,t}}
 +\sum_{i\in\mathcal P_h}\delta w_i\ind_{B_{i,t}}.
\end{equation}
We bound the trigger term $\delta d_t$ by a supplier lead-time demand cost and retailer mean demands, and bound the remaining width terms at each batch retailer. The indicators keep every width term tied to the period in which that retailer's request remains unfinished at the head of the queue.

\subsection{Step 4. Supplier Cost and Retailer Widths}\label{sec:arb-source-budget}
Let $C_0=\NV_{h,p}(D^{(L_0)})=\min f_{L_0}$ be the newsvendor cost of supplier lead-time demand, where $\NV_{h,p}$ is defined in Section~\ref{sec:arb-foundations} and $f_{L_0}$ in~\eqref{eq:source-value}. Recall from Section~\ref{sec:arb-policy} that $G_0=J_0(gK_0)$ is the supplier's planning value and $a_0=(-\hat s_0)^+$ the upward shift of both its levels.
\begin{lemma}[{\mdseries\scshape Supplier Cost and Shortages}]\label{lem:arb-source}
Supplier setup costs at the actual setup cost $K_0$ plus the loss $\phi(E_t)$ have long-run expected average at most $G_0+ha_0$; moreover, $ha_0+C_0\le G_0$, and the shortage event $E_t<0$ has long-run expected frequency at most $\alpha=h/(h+p)$.
\end{lemma}
\proof{Proof.} The upward shift adds bounded holding cost and can only reduce shortages. The canonical endpoint conditions and target optimality bound supplier cost and shortage frequency.

\emph{Step 1. Bound supplier costs.} By~\eqref{eq:arb-source}, $\E[\phi(E_{t+L_0})\mid\mathcal F_t]=f_{L_0}(Y_t)$, so their long-run expected averages agree. The upward shift $a_0$ keeps the order frequency and, as $f_{L_0}$ has slopes at most $h$, raises this average by at most $ha_0$; since actual setup is below planning setup, the supplier costs average at most $G_0+ha_0$. If $\hat s_0<0$, nonnegative lead-time demand and~\eqref{eq:canonical} give $G_0=f_{L_0}(\hat s_0)=p(\E D^{(L_0)}+a_0)$, and $C_0\le f_{L_0}(0)=p\E D^{(L_0)}$, so $h\le p$ gives $ha_0+C_0\le G_0$. If $a_0=0$, the claim $C_0\le G_0$ holds because every policy, including the supplier's, pays at least $\min f_{L_0}=C_0$ per period.

\emph{Step 2. Bound the shortage frequency.} Let $X_0$ be the long-run depletion of the unshifted canonical policy, independent of $D^{(L_0)}$. By Lemma~\ref{lem:scalar}, the canonical policy is optimal, including among policies of its width, which share its order frequency, so $\hat S_0$ minimizes $\E\phi(S-X_0-D^{(L_0)})$ over $S$, and the right derivative $h-(h+p)\Prob(X_0+D^{(L_0)}>\hat S_0)$ is nonnegative. Thus, shortage frequency is at most $\alpha$; the upward shift only decreases it. Expected renewal averages give this bound for the average of $\Prob(E_t<0)$. The same argument applies at a zero-setup base-stock minimizer.
\Halmos\endproof

We next bound the demand that triggers the unfinished request at the head of the queue. We first compare the newsvendor costs of the batch retailers' demand and of total demand, then apply Lemma~\ref{lem:arb-max} to the batch retailers' demand totals across periods. Let $\mathcal P$ be the set of batch retailers, $\mu_i=\E D_i$, $D_t^{\mathcal P}:=\sum_{i\in\mathcal P}D_{i,t}$, and $D_t^{\mathcal P^c}:=\sum_{i\notin\mathcal P}D_{i,t}$, where the first total is the demand of the batch retailers, which use cycle policies, and the second that of the immediate retailers, which use base stock. The following aggregate condition suffices for the entire proof:
\begin{equation}\label{eq:arb-aggregate-dependence}
 \operatorname{Cov}\!\left(D_t^{\mathcal P^c},\ind_{\{D_t^{\mathcal P}>x\}}\right)\ge0
 \qquad\text{for every }x.
\end{equation}
Association implies~\eqref{eq:arb-aggregate-dependence}: take $u(\mathbf D_t)=\min\{D_t^{\mathcal P^c},b\}$ and $v(\mathbf D_t)=\ind_{\{D_t^{\mathcal P}>x\}}$, which are bounded and coordinatewise nondecreasing, and let $b\to\infty$.\footnote{Recall from Section~\ref{sec:arbitrary} that $\mathbf D_t$ is associated if $\operatorname{Cov}(u(\mathbf D_t),v(\mathbf D_t))\ge0$ for all bounded coordinatewise nondecreasing functions $u$ and $v$. The finite first moment of $D_t^{\mathcal P^c}$ justifies the limit as $b\to\infty$ by dominated convergence.} Association can be weakened to positive quadrant dependence of the two totals, meaning $\Prob(D_t^{\mathcal P}>x,D_t^{\mathcal P^c}>y)\ge\Prob(D_t^{\mathcal P}>x)\Prob(D_t^{\mathcal P^c}>y)$ for every $x,y$. Integrating the difference between the two sides over $y\ge0$ gives the covariance in~\eqref{eq:arb-aggregate-dependence}. Thus, this condition restricts only the dependence between the two totals, not within either group. If every retailer is a batch retailer, $D_t^{\mathcal P^c}=0$ and the condition holds for arbitrary within-period dependence.

To prove the newsvendor comparison, write $X=D_t^{\mathcal P}$, $Y=D_t^{\mathcal P^c}$, and $f_z(x)=h(z-x)^++p(x-z)^+$. The convex function $f_z$ has bounded increasing subgradient $s_z(x)=-h+(h+p)\ind_{\{x>z\}}$. Its supporting-line inequality, the identity $\E[s_z(X)(Y-\E Y)]=(h+p)\operatorname{Cov}(\ind_{\{X>z\}},Y)$, and~\eqref{eq:arb-aggregate-dependence} give, for every target $z$,
\begin{equation}\label{eq:arb-dependent-loss}
 \E f_z(X+Y-\E Y)
 \ge\E f_z(X)+(h+p)\operatorname{Cov}(\ind_{\{X>z\}},Y)
 \ge\E f_z(X).
\end{equation}
All terms are integrable because $f_z$ is Lipschitz, $s_z$ is bounded, and demands have finite first moments. Since $f_z(x+c)=f_{z-c}(x)$, applying~\eqref{eq:arb-dependent-loss} at target $z-c-\E Y$ gives $\E f_z(X+Y+c)\ge\E f_z(X+\E Y+c)$ for every constant $c$, so replacing the complementary demand $Y$ by its mean cannot raise the expected loss. Replacing the complementary demand of each period $\ell=1,\ldots,L_0$ by its mean, one period at a time,\footnote{Because periods are independent, conditional on the demands of all other periods, the two totals of a period keep their joint distribution, and the other periods add a constant $c$, so the replacement can be made one period at a time.} therefore, gives $\E f_z(\sum_{\ell=1}^{L_0}D_\ell)\ge\E f_z(\sum_{\ell=1}^{L_0}D_\ell^{\mathcal P}+L_0\E D_t^{\mathcal P^c})$. Minimizing over $z$ and absorbing the constant $L_0\E D_t^{\mathcal P^c}$ into the target proves
\begin{equation}\label{eq:arb-protected-newsvendor}
 C_{\mathcal P}:=\NV_{h,p}\!\left(\sum_{\ell=1}^{L_0}D_\ell^{\mathcal P}\right)
 \le\NV_{h,p}\!\left(\sum_{\ell=1}^{L_0}D_\ell\right)=C_0.
\end{equation}

Suppose $L_0>0$ and $\mathcal P$ is nonempty. When the unfinished request at the head of the queue comes from a retailer in $\mathcal P$, its triggering demand $d_t$ in~\eqref{eq:arb-split} is that retailer's demand at a date in $\{t-L_0,\ldots,t-1\}$, so it is at most that date's total over $\mathcal P$; otherwise, $d_t=0$. Thus, $d_t\le V_{\max,t}:=\max_{t-L_0\le\ell<t}D_\ell^{\mathcal P}$, where the $L_0$ period totals in this maximum are independent and have common mean $\sum_{i\in\mathcal P}\mu_i$. The queue is nonempty only when $E_t<0$, so Lemma~\ref{lem:arb-max}, through~\eqref{eq:arb-quantile} and Lemma~\ref{lem:arb-source}, bounds the long-run expected average trigger cost $(h+p)\delta d_t$ by $\delta(1+h/p)C_{\mathcal P}+\delta h\sum_{i\in\mathcal P}\mu_i\le(3/4)C_0+\delta\sum_{i\in\mathcal P}h\mu_i$, where the inequality uses~\eqref{eq:arb-protected-newsvendor}, $\delta=3/8$, and $h\le p$. We call $\delta h\mu_i$ the mean terms and charge them to their retailers. The supplier cost and the pooled term $\frac34C_0$ are at most
\begin{equation}\label{eq:arb-source-budget}
 G_0+ha_0+\frac34C_0\le2G_0,
\end{equation}
where the last step uses $ha_0+C_0\le G_0$ from Lemma~\ref{lem:arb-source}. If $L_0=0$ or $\mathcal P$ is empty, set $V_{\max,t}=0$; the reserve is zero and no within-period dependence condition is needed, as the trigger term vanishes.

For a retailer $i$ using the full-cycle policy, suppress its index in $G=J(k)$, $k$, and $\hat w$, and put $v=h\hat w/G$, $\widetilde X_t=X_{i,t}/\hat w$, and $\widetilde X=X_i/\hat w$, its normalized long-run depletion, with $X_i$ distributed as $X_{\hat w}$ in~\eqref{eq:occupation}. Two cost terms belong to this retailer: the mirror loss $hX_{i,t}$ in~\eqref{eq:arb-credited}, and the width term $(h+p)\hat w\psi(\widetilde X_t)\ind_{B_{i,t}}$, which bounds its share of the credited reserve cost by~\eqref{eq:arb-decomposition}, multiplied by $h+p$. The expected width term is at most $(h+p)\delta\hat w\Prob(E_t<0)$ and at most $(h+p)\hat w\E\psi(\widetilde X_t)$, so its long-run average is at most $\min\{h\delta\hat w,h\hat w\E\psi(\widetilde X)/\alpha\}$,\footnote{The first bound on the expected width term uses $\psi\le\psi(0)=\delta$, for the nonincreasing $\psi$ of Section~\ref{sec:arb-phase}, and $B_{i,t}\subseteq\{E_t<0\}$; the second drops the indicator. Since the expected width term is at most each bound in every period, its long-run average is at most the smaller of their long-run averages. Lemma~\ref{lem:arb-source} bounds the long-run average of $\Prob(E_t<0)$ by $\alpha$, and $(h+p)\alpha=h$, so the first bound averages at most $h\delta\hat w$. Since $\psi$ is bounded,~\eqref{eq:occupation} makes $\E\psi(\widetilde X)$ the long-run average of $\E\psi(\widetilde X_t)$, so the second averages $(h+p)\hat w\E\psi(\widetilde X)=h\hat w\E\psi(\widetilde X)/\alpha$.} while the mirror loss averages $h\E X_i=h\hat w\E\widetilde X$. Since $v=h\hat w/G$, the sum of the two averages, divided by $G$, is at most
\begin{equation}\label{eq:arb-W}
 \mathcal W(v,\widetilde X)=v\left[\E\widetilde X+\min\left\{\delta,\frac{\E\psi(\widetilde X)}{\alpha}\right\}\right].
\end{equation}
These bounds use event containment and separate expectations, so the shortage event, the identity of the request at the head, and depletion need not be independent. Since only the unfinished request at the head holds a reserve,~\eqref{eq:arb-split} and the bound $\delta(w_i+d)$ from Section~\ref{sec:arb-phase} for a retailer using the half-cycle policy bound the expectation over the offset of the credited reserve cost $(h+p)(R-X)^+$ by the trigger term $(h+p)\delta d$ plus a width term for retailer $i$, incurred only on the event $B_{i,t}$ that its request is unfinished at the head of the queue. For a retailer using the half-cycle policy, this width term is $(h+p)\delta w_i\ind_{B_{i,t}}$, and the argument for the first bound above shows it averages at most $h\delta w_i$. The trigger bound uses the sum of batch retailers' mean demands, whereas each width term is charged to the retailer whose request is at the head. We therefore bound the trigger and width terms separately and charge every batch retailer both its mean term $\delta h\mu_i$ from the trigger bound and its own width term, regardless of which retailer's request is at the head.

\subsection{Step 5. A Holding-Rate Width Bound}\label{sec:arb-uniform}
The next bound retains the retailer holding-rate condition $H\ge mh$, with $e=1/m\in(0,1]$. It links the width and mirror bound $\mathcal W$ to the share $k/[GM(\hat w)]$ of the planning value spent on setups. The ratios $v=h\hat w/G$ and $\widetilde X=X_i/\hat w$ are the normalized quantities of~\eqref{eq:arb-W} in Step~4.

\begin{lemma}[{\mdseries\scshape Canonical Width Bound}]\label{lem:arb-uniform}
Suppose $H\ge mh$ and $p\ge h$. For a retailer using the full canonical policy with positive planning setup cost $k$, the bound~\eqref{eq:arb-W} satisfies
\begin{equation}\label{eq:arb-uniform}
 \mathcal W(v,\widetilde X)
 \le\frac78(1+e)-\frac{7e}{12}\frac{k}{GM(\hat w)}.
\end{equation}
\end{lemma}

\proof{Proof.} The canonical endpoints bound the cycle savings. We keep the bound $e$ on the holding-rate ratio and use the depletion's support and mean to reduce the claim to three endpoint checks.

\emph{Step 1. Retain the holding-rate ratio in the savings cap.} Set $\eta=h/H\le e$ and $\rho=h/p\le1$, so $\alpha=h/(h+p)=\rho/(1+\rho)$. For every real $y$,
\[
 \frac{H(\hat s-y)^++p(\hat s-y)^-}{p}
 +\frac{H(\hat S-y)^++p(\hat S-y)^-}{H}
 \ge\hat S-\hat s=\hat w,
\]
where the two terms are $y-\hat s$ and $\hat S-y$ if $\hat s\le y\le\hat S$, and otherwise one term alone is at least $\hat w$. Taking expectations at lead-time demand gives $F(\hat s)/p+F(\hat S)/H\ge\hat w$. Since $F(\hat s)=G$ and $F(\hat S)\le G$, this implies $v\le\eta+\rho$ and $F(\hat S)/G\ge(v-\rho)/\eta$. Nonnegativity of $F$ and its slopes in $[-p,H]$ therefore give
\begin{equation}\label{eq:arb-cap}
 \frac{G-F(\hat S-\hat wz)}{G}
 \le\bar g_\eta(z):=
 \min\left\{1,\frac{v(1-z)}{\rho},
                  \frac{\eta+\rho-v+vz}{\eta}\right\}
 \qquad(0\le z\le1).
\end{equation}
When the third term is below $1$, its numerator after subtracting $1$ is $\rho-v+vz<0$, so increasing $\eta$ to $e$ increases that term. When this numerator is nonnegative, the third term is at least $1$ and does not affect the minimum. Hence $\bar g_\eta\le\bar g_e$, and $v\le e+\rho$. Lemma~\ref{lem:reward} and~\eqref{eq:occupation} give
\[
 u:=\frac{k}{GM(\hat w)}
 =\E\frac{G-F(\hat S-\hat w\widetilde X)}G
 \le\E\bar g_e(\widetilde X),
 \qquad \E\widetilde X\le\frac12,
\]
where the mean bound is Lemma~\ref{lem:renewal}. It suffices to prove
$\mathcal W+(7e/12)\E\bar g_e(\widetilde X)\le C_e$, where $C_e:=7(1+e)/8$, for every random variable $\widetilde X$ on $[0,1)$ that satisfies this mean bound.

\emph{Step 2. Dispose of the small-width cases.} Since $\delta=3/8$, equation~\eqref{eq:arb-W} gives $\mathcal W\le7v/8$. If $\rho\le1-2e/3$, use $v\le e+\rho$ and $\bar g_e\le1$ to obtain
\[
 \mathcal W+\frac{7e}{12}\E\bar g_e(\widetilde X)
 \le\frac78(e+\rho)+\frac{7e}{12}\le C_e.
\]
If $v\le\rho$, the same argument gives a bound of $7\rho/8+7e/12\le7/8+7e/12\le C_e$. It remains to consider
\[
 \rho>1-\frac{2e}{3}\ge\frac13,
 \qquad \rho<v\le e+\rho.
\]

\emph{Step 3. Reduce the remaining case to three endpoints.} Put $x=v\widetilde X$, and define
\[
 \psi_v(x):=v\psi(x/v)=\frac{[(v-4x/3)^+]^2}{8v/3},
 \qquad
 \bar g_v(x):=\bar g_e(x/v)
 =\min\left\{1,\frac{v-x}{\rho},
                   \frac{e+\rho-v+x}{e}\right\}.
\]
Bound the minimum in~\eqref{eq:arb-W} by its convex combination with weights $2\rho/(1+\rho)$ on $\E\psi(\widetilde X)/\alpha$ and $(1-\rho)/(1+\rho)$ on $\delta$. Both weights are nonnegative and sum to one, and this gives
\begin{equation}\label{eq:arb-blend}
 \mathcal W\le\E x+2\E\psi_v(x)
                  +\frac{3v(1-\rho)}{8(1+\rho)}.
\end{equation}
The function
$3x/4+2\psi_v(x)+(7e/12)\bar g_v(x)$ is convex on each of $[0,v-\rho]$ and $[v-\rho,v]$: $\psi_v$ is convex and the cap is affine on each interval. Its maximum is therefore attained at $0$, $v-\rho$, or $v$. To each of these three endpoint values, add $\E x/4\le v/8$ and the last term of~\eqref{eq:arb-blend}.

At $0$, the resulting total is
$v(8-\rho)/[12(1+\rho)]+7(e+\rho)/12$; at $v$, it is $v(5+2\rho)/[4(1+\rho)]$. Using $v\le e+\rho$, both are at most
\[
 A_e(\rho):=\frac{(e+\rho)(5+2\rho)}{4(1+\rho)}
 \le C_e,
 \qquad
 C_e-A_e(\rho)=
 \frac{(1-\rho)(4\rho+7-3e)}{8(1+\rho)}\ge0.
\]

\emph{Step 4. Check the middle endpoint.} At $x=v-\rho$, where $\bar g_v(x)=1$, the total is
\[
 T_e(v,\rho)=\frac{7e}{12}-\frac{3\rho}{4}
  +\frac{v(5+2\rho)}{4(1+\rho)}
  +\frac{[(4\rho-v)^+]^2}{12v},
\]
which is convex in $v$: its first three terms are affine, and the last term is convex, since $(4\rho-v)^2/v=16\rho^2/v-8\rho+v$ on $v<4\rho$ and joins zero with zero slope at $v=4\rho$. It suffices to check $v=\rho$ and $v=e+\rho$. At the first endpoint,
\[
 T_e(\rho,\rho)
 =\frac{7e}{12}+\frac{\rho(5+2\rho)}{4(1+\rho)}
 \le\frac{7e}{12}+\frac78\le C_e,
\]
because the fraction increases on $0<\rho\le1$. At the second endpoint, fix $\rho\in[1/3,1]$ and define $\Delta_\rho(e):=T_e(e+\rho,\rho)-C_e$ on the enlarged interval $0\le e\le1$. Since $3\rho\ge e$, the positive part may be removed, and
\[
 \frac{(3\rho-e)^2}{12(e+\rho)}
 =\frac{e-7\rho}{12}+\frac{4\rho^2}{3(e+\rho)}
\]
shows that $\Delta_\rho$ is convex in $e$, because its remaining terms are affine. Its endpoint values are
\[
 \Delta_\rho(0)=\frac{(\rho-1)(4\rho+7)}{8(1+\rho)}\le0,\qquad
 \Delta_\rho(1)=\frac{(\rho-1)(3\rho-1)}{6(1+\rho)}\le0.
\]
Thus $\Delta_\rho(e)\le0$ throughout $[0,1]$, giving $T_e(e+\rho,\rho)\le C_e$. Convexity in $v$ now bounds the middle endpoint throughout its range. Together with Step~3, this proves $\mathcal W+(7e/12)\E\bar g_e(\widetilde X)\le C_e$. Since $u\le\E\bar g_e(\widetilde X)$, this gives $\mathcal W\le C_e-(7e/12)u$ with $C_e=7(1+e)/8$, which is~\eqref{eq:arb-uniform}.
\Halmos\endproof

\subsection{Step 6. The Local Cost Comparison}\label{sec:arb-certificate}
An immediate retailer has zero depletion and reserve. All its shipments in a period form one dispatch, so its setup cost is at most $K_i$ per period, regardless of how its requests are split. Its total contribution is therefore at most $\Lambda_i+K_i$, which the immediate-retailer test bounds by $BG_i$.

For a batch retailer, suppress its index and write $u=k/[GM(\hat w)]\in(0,1]$ for the setup-cost share of its canonical plan at $k=\zeta K$. Its canonical holding and backlog cost is $G-k/M(\hat w)=G(1-u)\ge\Lambda$. Failing the immediate test $\Lambda+K\le BG$ therefore implies
\begin{equation}\label{eq:arb-threshold}
 \frac{k}{G}>\underline\kappa(u):=\zeta(B-1+u)>1.
\end{equation}
Indeed, $\zeta(B-1)=g(2g-3)/(g-1)=2g-1-1/(g-1)\ge52/21$, since $g\ge13/6$. Thus $k>G$, and Lemma~\ref{lem:arb-mean}, with $\nu=uG/k$, gives $h\mu/G\le(1-u)/(1-uG/k)$. This bound decreases in $k/G$, so
\begin{equation}\label{eq:arb-rK}
 \frac{h\mu}{G}\le r(u):=
 \frac{1-u}{1-u/\underline\kappa(u)}.
\end{equation}

For the full-cycle policy, holding and backlog cost $G(1-u)$ plus charged-shipment setup cost $(7/3)K/M(\hat w)$ equals $G[1+(7/(3\zeta)-1)u]$. Adding the width and mirror bound of Lemma~\ref{lem:arb-uniform} and the mean contribution $(3/8)h\mu\le(3/8)r(u)G$ gives the upper bound $GF_e(u)$, where
\begin{equation}\label{eq:arb-full}
 F_e(u):=1+\frac{7(1+e)}8+
 \left(\frac7{3\zeta}-1-\frac{7e}{12}\right)u+\frac38r(u).
\end{equation}

For the half-cycle policy, Lemma~\ref{lem:shift} gives a width $w=\hat w/2$, target $\Gamma=\hat S-\upsilon^*$, and savings
$\Phi=\int_{[0,w)}[G-F(\Gamma-x)]\,U(dx)\ge k/2$.
Its holding and backlog cost plus charged-shipment setups is at most
\[
 G+\frac{(7/3)K-\Phi}{M(w)}
 \le G+\left(\frac7{3\zeta}-\frac12\right)\frac{k}{M(w)}
 \le G\left[1+\left(\frac{14}{3\zeta}-1\right)u\right],
\]
where the last inequality uses $M(\hat w)\le2M(w)$ and $7/(3\zeta)-1/2>0$. The endpoint bound in the proof of Lemma~\ref{lem:arb-uniform} gives $h\hat w/G\le h/H+h/p\le1+e$. Therefore, Lemma~\ref{lem:renewal} bounds its mirror cost by $h\E X_w\le hw/2\le(1+e)G/4$, and its reserve-width cost is $h\delta w\le3(1+e)G/16$. Adding these terms and $(3/8)r(u)G$ gives the half-cycle upper bound $GH_e(u)$, where
\begin{equation}\label{eq:arb-half-cost}
 H_e(u):=1+\frac{7(1+e)}{16}
       +\left(\frac{14}{3\zeta}-1\right)u+\frac38r(u).
\end{equation}

\begin{lemma}[{\mdseries\scshape Local Cost Comparison}]\label{lem:arb-certificate}
With the parameters in~\eqref{eq:arb-parameters}, the half-cycle bound satisfies $H_e(u)\le B$ for $0\le u\le u_*$, and the full-cycle bound satisfies $F_e(u)\le B$ for $u_*\le u\le1$.
\end{lemma}
\proof{Proof.} The two envelopes agree at $u_*$. We show that the half-cycle envelope increases up to this point, the full-cycle envelope decreases after it, and their common value is below $B$.

\emph{Step 1. Check the parameters and monotonicity.}
We verify the inequalities for $e\in[0,1]$, including the limiting endpoint $e=0$. Write $\kappa_0=\zeta(B-1)$. The coarse parameter bounds
\[
 2<g<\frac52,\qquad \frac53<\zeta<2,\qquad \kappa_0>2
\]
follow from~\eqref{eq:arb-parameters}; for the last, use
$\kappa_0=2g-1-1/(g-1)$. Moreover,
\[
 \frac14<\frac{3\zeta}{16}\le u_*<\frac12.
\]
Indeed, $(1+e)/(4+e\zeta)\ge1/4$ because $\zeta<4$, and the upper bound is equivalent to $\zeta(3+e)<8$.

Direct differentiation gives
\[
 r''(u)=-\frac{2\kappa_0(\kappa_0+\zeta-1)}
                   {[\kappa_0+(\zeta-1)u]^3}<0.
\]
Thus $r'(u)\ge r'(1)=-2g/(2g-1)>-g$ on $[0,1]$, and
\[
 H_e'(u)\ge\frac{14}{3\zeta}-1-\frac{3g}{8}
 >\frac73-1-\frac{15}{16}>0.
\]

For the full-cycle envelope, the identity
\[
 r'(u)=-1+
 \frac{\kappa_0(1-2u)-(\zeta-1)u^2}
      {[\kappa_0+(\zeta-1)u]^2}
\]
and $1/4<u_*<1/2$, $\kappa_0>2$ imply
\[
 r'(u_*)\le-1+\frac{1-2u_*}{\kappa_0}<-\frac34.
\]
The affine coefficient $7/(3\zeta)-1-7e/12$ decreases with $e$ and is $10/39$ at $e=0$: its derivative is $7(16/57)/(3g^2)-7/12<0$. Thus
\[
 F_e'(u_*)<\frac{10}{39}-\frac9{32}<0.
\]
Since $r''<0$ makes $F_e$ concave, the full-cycle envelope decreases throughout $[u_*,1]$.

\emph{Step 2. Bound the common endpoint.}
Write
\[
 V=4[(4+e)g-4],\qquad U=3g(1+e),\qquad
 D=g(2g-3)V+U.
\]
These quantities are positive, and $u_*=U/V$. Substitution of $g=(247+32e)/114$ and collection of powers of $e$ give
\[
 VD\,[B-F_e(u_*)]>105+50e-152e^2\ge3,
\]
where the first inequality rounds the constant and linear coefficients downward, discards the positive terms of degrees five to seven, and uses $e^3,e^4\le e^2$ to combine the three negative terms; the second follows because the quadratic function $105+50e-152e^2$ has a negative coefficient on $e^2$, so it is concave, and its endpoint values are $105$ and $3$ at $e=0$ and $e=1$. Hence $F_e(u_*)=H_e(u_*)<B$, and the monotonicity in Step~1 gives $H_e(u)\le B$ for $u\le u_*$ and $F_e(u)\le B$ for $u\ge u_*$.
\Halmos\endproof

Let $\pi_{\rm random}$ be the policy of Section~\ref{sec:arb-policy} from the clean start, with independent uniform offsets. For each batch retailer, let $\Omega_i$ be its selected full- or half-cycle bound, $G_iF_e(u_i)$ or $G_iH_e(u_i)$. It includes holding and backlog cost, charged shipments, mirror cost, reserve-width cost, and the mean contribution $\delta h\mu_i$. Taking long-run averages in~\eqref{eq:arb-credited}, with the reserve bounded by~\eqref{eq:arb-decomposition} and the supplier and trigger costs by Section~\ref{sec:arb-source-budget}, gives
\[
 C(\pi_{\rm random})\le
 \left[G_0+ha_0+\frac34C_0\right]
 +\sum_{i\notin\mathcal P}(\Lambda_i+K_i)
 +\sum_{i\in\mathcal P}\Omega_i,
\]
where~\eqref{eq:arb-source-budget} bounds the bracket by $2G_0$, the immediate-retailer test bounds each term in the first sum by $BG_i$, and Lemma~\ref{lem:arb-certificate} does the same for each $\Omega_i$. Hence,
\begin{equation}\label{eq:arb-result}
 C(\pi_{\rm random})\le2G_0+B\sum_iG_i
 =2\left[G_0+(g-1)\sum_iG_i\right]
 \le2g\OPT=\bar\gamma_m\OPT,
\end{equation}
where $B=2(g-1)$ gives the first equality, while $\zeta=g/(g-1)$, $G_i=J_i(\zeta K_i)$, and $G_0=J_0(gK_0)$ identify the bracket with the shared-cost lower bound~\eqref{eq:lower} at $\gamma=g$. The following deterministic selection and initialization arguments transfer this guarantee to the consolidation policy $\pi_{\rm con}^{(m)}$.

\subsection{Step 7. Deterministic Offset Selection}\label{sec:arb-derandomization}
Index ordinary requests in creation order with deterministic tie-breaking, excluding the initial requests of Section~\ref{sec:arb-initial}. For a request $j$ of size $Q_j$ from a batch retailer $i$, normalize its offset $\xi$ to $\vartheta=\xi/b\in(0,1)$, where $b=3Q_j/4$, as in Section~\ref{sec:arb-policy}. Let $N_j(\vartheta)$ be its number of interior thresholds and $\tau_j$ its creation period. To isolate the cost that depends on the offset of request $j$, its shipment setups and its credited reserve cost over its full lifetime, we use the \emph{credited full-lifetime proxy}
\begin{equation}\label{eq:arb-proxy}
 \mathcal Z_j(\vartheta)=K_i[1+N_j(\vartheta)]
 +(h+p)\sum_{t=\tau_j}^{\tau_j+L_0}
                     (R_{j,t}(\vartheta)-X_t)^+,
\end{equation}
where $X_t$ is total ideal depletion across retailers, and $R_{j,t}$, the stock allocated to request $j$ and retained at the warehouse in period $t$, is zero before any allocation and after completion. The proxy is nonnegative and bounded by $3K_i+(h+p)L_0Q_j$, since $N_j\le2$ and $R_{j,t}\le3Q_j/4$.

At creation, let $\mathcal H_j$ contain the observed ideal state, new request sizes and queue, supplier order, and ordinary request index. Future ideal dynamics do not depend on offsets, by Section~\ref{sec:arb-policy}. Write $\mathcal Z_j(\vartheta)=\Psi(\vartheta,\sigma_j,\mathbf d_j)$, where $\sigma_j$ is that observed ideal state, including $\tau_j$, and $\mathbf d_j$ consists of the fresh demand vectors in periods $\tau_j,\ldots,\tau_j+L_0-1$. Define all conditional values by averaging over these demands: $z_j(\vartheta):=\int\Psi(\vartheta,\sigma_j,\mathbf d)\,\Prob(d\mathbf d)$ and $\bar z_j=\int_0^1z_j(\vartheta)\,d\vartheta$. Thus, $\E[\mathcal Z_j(\vartheta_j)\mid\mathcal H_j]=z_j(\vartheta_j)$ for every $\mathcal H_j$-measurable choice $\vartheta_j$, and the event $\{\tau_j<T\}$ is measurable with respect to the same creation history, as needed for applying the tower property and then summing over requests.

\begin{lemma}[{\mdseries\scshape Deterministic Offset Selection}]\label{lem:arb-selector}
There is a measurable deterministic nonanticipative rule choosing the offset $\vartheta_j$ of each batch retailer's request $j$ with $z_j(\vartheta_j)<\bar z_j+2^{-j}$.
\end{lemma}
\proof{Proof.} We first control how the full-lifetime proxy varies with the offset on each fixed future path. This regularity guarantees a rational offset whose conditional expected proxy is less than the conditional average plus $2^{-j}$, and we select such an offset using only the observed history.

\emph{Step 1. Establish bounded variation.} For a fixed ideal future path, request $j$ is fully allocated by $\tau_j+L_0$, so it reaches at most $L_0+1$ cumulative allocation levels, one per period, none depending on $\vartheta$. Its dispatch periods change only when one of the at most two thresholds crosses one of these levels, or when the second threshold appears at $\vartheta=1/3$. Thus, $\mathcal Z_j(\vartheta)$ is a step function with a uniformly bounded number of discontinuity points. Subtracting the ideal $X_t$ and taking a positive part adds none. Its total variation is at most twice the number of discontinuity points times $3K_i+(h+p)L_0Q_j$,\footnote{As noted after~\eqref{eq:arb-proxy}, the proxy lies between $0$ and $3K_i+(h+p)L_0Q_j$, so each change in its value is at most $3K_i+(h+p)L_0Q_j$. The proxy is constant between consecutive discontinuity points, and at each such point its value can differ from both one-sided limits, so each point contributes at most two such changes.} a deterministic function of the observed $Q_j,K_i,L_0$. Conditional expectation given $\mathcal H_j$ preserves this bound, so $z_j$ has bounded variation and is Riemann integrable.

\emph{Step 2. Select a rational offset measurably.} Since $z_j$ is Riemann integrable, its average over the rational offsets $(2k-1)/(2n)$, $1\le k\le n$, converges to $\bar z_j$, so for large $n$ one of them has $z_j(\vartheta)<\bar z_j+2^{-j}$. Enumerate the rational numbers in $(0,1)$ once and choose the first such offset.\footnote{This rule establishes the existence of a deterministic policy satisfying the guarantee of Theorem~\ref{thm:arbitrary}. Each comparison involves the conditional expectations $z_j(\vartheta)$ and $\bar z_j$ over $L_0$ periods of all retailers' demands, and its tolerance $2^{-j}$ halves with each new request. In practice, drawing each offset uniformly from $(0,1)$, as the randomized policy $\pi_{\rm random}$ does, costs one random number per request and satisfies the same guarantee by~\eqref{eq:arb-result}.}

This choice is measurable, since each comparison is measurable.\footnote{Because lead times are finite, the ideal state of Section~\ref{sec:arb-policy}, in which every allocation ships at once, consists of finitely many quantities and a finite list of queued requests, so it is standard Borel. Future demand vectors are i.i.d.\ over periods and independent of it, and the values of $z_j$ at rational offsets are Borel functions of it.} It also uses only current information, because the creation history $\mathcal H_j$ determines $z_j$, $\bar z_j$, and the index $j$.
\Halmos\endproof

For the deterministic policy, the expected position and setup costs through period $T-1$ are at most the supplier setup costs, the offset-free terms in~\eqref{eq:arb-credited}, the immediate-retailer setup charges, and $\E\sum_{\tau_j<T}\mathcal Z_j(\vartheta_j)$. Each request of a batch retailer is charged for all its shipments, and any reserve held before $T$ belongs to a request created before $T$; extending its reserve charge to its full lifetime only increases the bound. For these requests, Lemma~\ref{lem:arb-selector}, the tower property, and $\sum_j2^{-j}\le1$ give $\E\sum_{\tau_j<T}\mathcal Z_j(\vartheta_j)\le\E\sum_{\tau_j<T}\mathcal Z_j(\Theta_j)+1$, where $\Theta_j$ are independent uniform offsets on $(0,1)$ applied to the same ideal process, so $\E[\mathcal Z_j(\Theta_j)\mid\mathcal H_j]=\bar z_j$. The setup charges in $\mathcal Z_j(\Theta_j)$ have expectation $(7/3)\sum_{i\in\mathcal P}K_i\E N_i(T)$, where $N_i(T)$ counts retailer $i$'s requests created before $T$. All requests created before $T$ finish by $T+L_0$. Only the unfinished request at the head of the ideal queue can retain stock, so the reserve part of the sum is at most the cumulative $(h+p)(R_t-X_t)^+$ of $\pi_{\rm random}$ through period $T+L_0$. Hence, the bounds of Sections~\ref{sec:arb-phase} to~\ref{sec:arb-certificate} for $\pi_{\rm random}$ apply.

To complete the comparison, the remaining terms in \eqref{eq:arb-credited}, $\phi(E)$ and $\sum_i[F_i(\widehat P_i)+hX_i]$, do not depend on the offsets. Dividing by $T$ and taking the limit superior yields \eqref{eq:arb-result} for the deterministic policy, which selects each offset by the nonanticipative rule of Lemma~\ref{lem:arb-selector}.

\subsection{Step 8. Initialization and Finite First Moments}\label{sec:arb-initial}
From a prescribed finite valid initial state, let prescribed pipelines arrive without new orders, send remaining warehouse stock to a fixed retailer, and wait for receipt. Withhold further actions until every retailer net inventory is at most its selected target $\Gamma_i$. Positive inventory during this preparation has deterministic upper bounds, so retailer $i$'s excess is removed once its cumulative demand reaches a constant. The product bound in Step~1 of Appendix~\ref{app:initialization} then gives preparation time finite moments of every order and the demand accumulated during preparation a finite mean. Thus, preparation duration and cost are integrable without a second demand moment.

Initialize positions including outstanding requests at $\Gamma_i$ by inserting initial requests of sizes $\Gamma_i-I_i\ge0$, where $I_i$ is retailer $i$'s net inventory at the end of preparation, and order the supplier to its shifted target. These request sizes and the first supplier order are integrable. All initial requests finish within $L_0$ periods, as shown in Section~\ref{sec:arb-policy}. Give the initial requests of batch retailers fixed offsets, so each has at most three charged shipments; the initial requests of immediate retailers have at most $L_0+1$ dispatch periods. Hence, with integrable request sizes, the transient setup and reserve costs of all initial requests, which arise only at initialization, are integrable.

To transfer the clean-start bound to this initialization, compare with a clean start having retailer inventories $\Gamma_i$, no initial requests, and a first supplier order equal to the shifted target. The actual first supplier order includes exactly $\sum_i(\Gamma_i-I_i)$ extra units. Before period $L_0$, no stock reaches the warehouse in either system; when the first orders arrive, the extra units fill exactly the initial requests, leaving the same stock and ordinary queue in both systems. Hence, their ideal states agree by period $L_0$. Ordinary request histories agree, and so do their indices, because initial requests are excluded from indexing. Deterministic offsets for requests created before $L_0$ can differ because their creation states differ, but all such requests finish by $2L_0$. Later requests see the same ideal states and indices and receive the same offsets. Hence, reserves and positions agree by period $2L_0$, and physical retailer inventories by period $2L_0+\max_iL_i$, the length of the comparison horizon. Stocks during it are bounded by integrable initial stock and orders placed, and backlogs by the integrable initial requests plus demand over the horizon, so its expected cost is finite.

In ordinary periods, $R_t\le(3/4)(\max_i w_i+V_{\max,t})$ pathwise, with $V_{\max,t}=\max_{t-L_0\le\ell<t}D_\ell^{\mathcal P}$ the largest one-period total demand of the batch retailers over the $L_0$ periods before $t$, as in Section~\ref{sec:arb-source-budget}, because any request at the head has size $Q\le w_i+D_{i,\tau-1}$ and retains at most $3Q/4$. This integrable bound applies to randomized and selected offsets. The remaining long-run averages are controlled by the single-location policies: supplier and retailer post-action positions stay in bounded intervals. Request counts, single-location holding and backlog costs, shortage indicators, and bounded functions of the depletion have expected renewal averages. We apply~\eqref{eq:arb-quantile} to $V_{\max,t}$, and the proxy comparison of Section~\ref{sec:arb-derandomization} to the deterministic rule, which uses only the observed history.

With $\mathcal F_t$ the information after dispatches and before demand is realized,~\eqref{eq:retailer-lead} and~\eqref{eq:retailer-value} give $\E[H_iI_{i,t+L_i}^++pI_{i,t+L_i}^-\mid\mathcal F_t]=F_i(P_{i,t})$. Hence, physical costs in periods $L_i,\ldots,T-1$ correspond to positions in periods $0,\ldots,T-L_i-1$. Extend these nonnegative sums to period $T-1$ and add the finite expected physical costs of the first $L_i$ periods, so retailer lead times add no average term.

Finally, let $a_N$ be the expected cost of the first $N$ periods from the clean start, and $B_{\rm net}=2G_0+B\sum_iG_i$. The bounds that give~\eqref{eq:arb-result} prove $\limsup_Na_N/N\le B_{\rm net}$. For each $\epsilon>0$ there is finite $C_\epsilon$ with $a_N\le(B_{\rm net}+\epsilon)N+C_\epsilon$ for every $N\ge0$. Condition on the preparation history, use the i.i.d.\ distribution of the fresh future demands, apply this bound, uniform in $N$, to the remaining horizon, and add the integrable preparation and coupling costs. Divide by the original horizon $T$, let $T\to\infty$, and let $\epsilon\downarrow0$. Since $B_{\rm net}\le \bar\gamma_m\OPT$ by~\eqref{eq:arb-result}, this proves Theorem~\ref{thm:arbitrary}.

\section{Smaller Constants for DFR Demand}\label{app:smooth}\label{sec:smooth}

The two-candidate comparison of Proposition~\ref{prop:scalar} applies to every demand distribution in our model. We obtain a smaller constant by using an additional shape property of demand. When positive demand has a density and a decreasing failure rate (DFR), meaning that its survival function is log-convex, each depletion distribution has a point mass at zero and a nonincreasing density elsewhere. The following Lemma~\ref{lem:quantile} proves a contraction property for distributions with convex quantile functions, including these depletion distributions; Theorem~\ref{thm:dfr} builds on it.

To explain how this property changes the policy, let $Y$ denote the original depletion and let $Y^{[t]}$ have the distribution obtained by retaining and renormalizing its lowest fraction $t$ of probability mass. Then $Y^{[t]}-\E Y^{[t]}\le_{\rm cx}t(Y-\E Y)$: after each variable is centered at its own mean, every convex cost has no larger expectation for the retained distribution than for the scaled original distribution. Here probability mass equals the long-run fraction of periods spent at the corresponding depletion levels. For a fraction above the mass at zero, shortening the cycle at the corresponding cutoff and retuning its retailer target and mirror center has three effects: first, it multiplies the mirror cost by at most $t$; second, it does not increase holding and backlog cost; third, it multiplies order frequency by $1/t$. These comparisons give the tradeoff between inventory and setup costs used in the policy below; fractions at or below the mass at zero use base stock, as specified there.

\begin{definition}[{\mdseries\scshape DFR Policy}]\label{def:dfr}
Fix $\theta\in[1/2,1)$. The deterministic FIFO policy $\pi_{\rm DFR}(\theta)$ is specified by four rules. It uses the initialization in Appendix~\ref{app:initialization} with the selected retailer policies and mirror centers. Choose the smallest minimizing retailer target or mirror center and break all other ties by fixed rules. We suppress the retailer index in the quantities below.

(i)~\emph{Planning quantities.} At each retailer with $K>0$, choose its canonical optimum at planning setup cost $K/\theta$, let $\hat w$ be its width, and compute $z_{\hat w}:=K/M(\hat w)$, $m_{\hat w}:=B_{\hat w}$ from~\eqref{eq:mirror}, and $t_0:=U\{0\}/M(\hat w)$, with $U$ and $M$ defined in~\eqref{eq:renewal}. Set $t:=\min\{1,\sqrt{z_{\hat w}/m_{\hat w}}\}$, with $t=1$ if $m_{\hat w}=0$.

(ii)~\emph{Replenishment rules.} A retailer with $K=0$ uses base stock at a minimizer of $F$. For $K>0$, use the canonical optimum when $t=1$. When $t<1$, use this base stock if $t\le t_0$, and otherwise the width-$w_t$ policy with $U[0,w_t)=tM(\hat w)$ and a target minimizing holding and backlog cost.

(iii)~\emph{Mirror centers.} Each retailer fixes a mirror center minimizing the mirror cost of its selected policy. For base stock, this center is the target and the mirror cost is zero.

(iv)~\emph{Requests, dispatch, and supplier orders.} Follow rules~(ii) and~(iii) of Definition~\ref{def:fifo}, using the retailer policies and mirror centers selected above.\Halmos
\end{definition}

\begin{theorem}[{\mdseries\scshape DFR Demand}]\label{thm:dfr}
Define $\lambda:=2-\sqrt2$ and $\gamma_\lambda:=1+\lambda/2+\sqrt{1+\lambda^2/4}$. Suppose $H_i\ge h$, $p\ge h$, and $K_i\le(\gamma_\lambda-1)K_0$ for every retailer. Assume that each retailer with $K_i>0$ has an absolutely continuous DFR distribution conditional on positive demand. The policy $\pi_{\rm DFR}(\theta)$ of Definition~\ref{def:dfr}, with $\theta=1-1/\gamma_\lambda$, is nonanticipative and feasible and satisfies $C(\pi_{\rm DFR}(\theta))\le \gamma_\lambda\OPT$ from every valid initial state and for every finite supplier lead time $L_0\ge0$.
\end{theorem}

The factor $\gamma_\lambda$ is approximately $2.335$. With $\theta=(\gamma-1)/\gamma$, dividing~\eqref{eq:lower} by $\gamma$ assigns weights $1-\theta$ and $\theta$ to the supplier and retailer values, respectively; the retailer planning setup cost is $K_i/\theta$. The simpler choice $\theta=7/12$ gives the slightly weaker factor $1/(1-\theta)=12/5$.

Zero demand may have point mass, and retailers with $K_i=0$ retain arbitrary demand distributions with finite positive means. At each retailer with $K_i>0$, positive demand must have a density, so integer-valued demand, as in Section~\ref{sec:numerics}, is excluded from this extension.\footnote{The following inequality~\eqref{eq:mirror-contract} can fail for a step quantile function. Also, with integer-valued demand, $U$ is supported on the nonnegative integers, so as the cutoff $w_t$ moves through $[0,\hat w]$, the mass $U[0,w_t)$ jumps at integers and takes only finitely many values. Hence the cutoff equation $U[0,w_t)=tM(\hat w)$ may have no solution for the fraction $t$ selected in Definition~\ref{def:dfr}. Rounding to an attainable fraction would require a separate cost bound.}

For demand conditional on being positive, distributions with positive log-convex densities satisfy the DFR assumption. Finite mixtures of these distributions retain positive log-convex densities and hence are also DFR.\footnote{For a positive density $f$ on $(0,\infty)$, its survival function is $\bar F(x)=\int_0^\infty f(x+u)\,du$. For $x,y>0$ and $0<\eta<1$, log-convexity of $f$ and H\"older's inequality give $\bar F(\eta x+(1-\eta)y)\le\int_0^\infty f(x+u)^\eta f(y+u)^{1-\eta}\,du\le\bar F(x)^\eta\bar F(y)^{1-\eta}$. Thus, $\log\bar F$ is convex, so the failure rate $f/\bar F=-(\log\bar F)'$ is nonincreasing. For a finite mixture $f(x)=\sum_j\omega_j f_j(x)$ of positive log-convex densities $f_j$ on $(0,\infty)$, with fixed weights $\omega_j\ge0$ and $\sum_j\omega_j=1$, log-convexity of each $f_j$ and H\"older's inequality for finite sums give $f(\eta x+(1-\eta)y)\le\sum_j\omega_j f_j(x)^\eta f_j(y)^{1-\eta}\le f(x)^\eta f(y)^{1-\eta}$. Thus, the mixture density is positive and log-convex, and the preceding integral argument makes its survival function log-convex as well.} In particular, finite mixtures of exponential distributions qualify, since each exponential density has an affine log-density. The condition concerns one-period demand; retailer lead-time demand may sum any finite number of these demands. Since the assumption restricts only each retailer's demand distribution and the proof uses Propositions~\ref{prop:upper} and~\ref{prop:lower}, Theorem~\ref{thm:dfr} also allows demand that is dependent across retailers within a period, as Theorem~\ref{thm:main} does.

To compare the retailer and mirror costs of the selected windows, we use the optimized inventory loss from Section~\ref{sec:arb-foundations}: for positive rates $a,b$ and an integrable random variable $X$, define
\begin{equation}\label{eq:nv}
 \NV_{a,b}(X):=\min_z\E\bigl[a(z-X)^++b(X-z)^+\bigr].
\end{equation}
This loss is unchanged by deterministic shifts, which the target $z$ absorbs, and is positively homogeneous: $\NV_{a,b}(cX)=c\NV_{a,b}(X)$ for $c>0$. The next lemma uses convex order to compare a shortened depletion window with the original one, and Lemma~\ref{lem:smooth-mirror} bounds mirror cost by width through a mixture of uniform distributions. Section~\ref{app:smooth-proofs} gives both proofs; the quantile comparison uses the characterization of centered convex order in \citet{ShakedShanthikumar2007}.

\begin{lemma}[{\mdseries\scshape Convex-Quantile Contraction}]\label{lem:quantile}
Let $Y$ be integrable with convex quantile function $Q$. For $0<t\le1$, let $Y^{[t]}$ have quantile function $Q(tu)$, $0<u<1$. For every integrable $Z$ independent of $Y$ and of $Y^{[t]}$, and every $a,b>0$,
\begin{align}
 \NV_{a,b}(Y^{[t]})&\le t\NV_{a,b}(Y),\label{eq:mirror-contract}\\
 \NV_{a,b}(Z+Y^{[t]})&\le(1-t)\NV_{a,b}(Z)+t\NV_{a,b}(Z+Y)
 \le\NV_{a,b}(Z+Y).\label{eq:running-contract}
\end{align}
\end{lemma}

\begin{lemma}[{\mdseries\scshape Nonincreasing-Density Mirror Bound}]\label{lem:smooth-mirror}
Suppose the distribution of $Y$ consists of possible mass at zero and a nonincreasing density on $(0,w)$, with no mass elsewhere. Then
\begin{equation}\label{eq:smooth-mirror}
 \NV_{p,h}(Y)\le p\left(1-\sqrt{\frac p{p+h}}\right)w.
\end{equation}
\end{lemma}

To apply the two lemmas to DFR demand, we now check the shape of the resulting depletion distributions. Let $q=\Prob(D>0)$, let $D^+$ have the conditional distribution of $D$ given $D>0$, and let $U_+$ be its renewal measure including the initial renewal at zero. For absolutely continuous DFR $D^+$, $x\mapsto U_+[0,x]$ is concave on $[0,\infty)$ \citep{Brown1980,Yu2011}. We choose $u_+$ as its right derivative on $(0,\infty)$, which is nonincreasing because the function is concave. Thus, $U_+=\delta_0+u_+(x)\,dx$, with $\delta_0$ the unit mass at zero and $u_+$ a nonincreasing density on $(0,\infty)$. Including zero-demand periods multiplies the expected count at every reached depletion level by $1/q$, giving $U=q^{-1}U_+$.\footnote{The partial sums of positive demands alone count each level they reach once, including the initial zero. The original sequence counts its first occurrence and each repetition after a zero demand, until the next positive demand. This count is geometric with success probability $q$, where success means observing positive demand and moving to a new level. Each zero demand has probability $1-q$ and adds another count, so the mean is $\sum_{n=0}^{\infty}(1-q)^n=1/q$. The counts are independent of the positive demand sizes, giving $U(A)=q^{-1}U_+(A)$ for every measurable set $A$.}

For $w>0$, taking $v(y)=\ind_{\{y\le x\}}$ in~\eqref{eq:occupation} and using $U=q^{-1}U_+$ gives the distribution function $\mathsf F_w(x):=\Prob(X_w\le x)=[1+\int_0^x u_+(y)\,dy]/[qM(w)]$ for $0\le x\le w$. Thus, $X_w$ has point mass $1/[qM(w)]$ at zero and density $u_+(x)/[qM(w)]$ on $(0,w)$. This density is nonincreasing because $u_+$ is nonincreasing and $qM(w)>0$. The distribution function $\mathsf F_w$ is therefore continuous and concave on $[0,w]$, and its quantile function is convex, including the initial interval on which it is zero.\footnote{The generalized inverse $Q_w(u):=\inf\{x\in[0,w]:\mathsf F_w(x)\ge u\}$, $0<u<1$, gives the smallest depletion at which cumulative probability reaches $u$. This definition allows jumps and flat portions of a distribution function; here $Q_w(u)=0$ for $u\le\mathsf F_w(0)=1/[qM(w)]$. To check convexity, set $x_j=Q_w(u_j)$ for $j=1,2$. For $0\le\eta\le1$, concavity gives $\mathsf F_w(\eta x_1+(1-\eta)x_2)\ge\eta u_1+(1-\eta)u_2$, hence $Q_w(\eta u_1+(1-\eta)u_2)\le\eta Q_w(u_1)+(1-\eta)Q_w(u_2)$. The point mass at zero creates a jump from negative arguments to zero, outside the interval $[0,w]$ on which concavity is used.} Thus, this depletion distribution meets the assumptions of Lemmas~\ref{lem:quantile} and~\ref{lem:smooth-mirror}.

\proof{Proof of Theorem~\ref{thm:dfr}.} We bound the full window's mirror cost and use the contraction inequalities to select a window or base stock. We then combine the retailer and supplier bounds and choose the weight $\theta\in[1/2,1)$ that balances the retailer and supplier approximation factors.

\emph{Step 1. Bound the full-window mirror cost.} Fix a retailer with $K>0$ and a weight $\theta\in[1/2,1)$. Take its canonical optimum at planning setup cost $K/\theta$, with width $\hat w$, long-run depletion $Y=X_{\hat w}$, and order frequency $\nu=1/M(\hat w)$. Let $Z$ be its lead-time demand, independent of $Y$, and write
\begin{equation}\label{eq:smooth-data}
 r_{\hat w}:=\NV_{H,p}(Z+Y),\qquad m_{\hat w}:=\NV_{p,h}(Y),\qquad
 z_{\hat w}:=K\nu,\qquad j:=J(K/\theta)=r_{\hat w}+z_{\hat w}/\theta.
\end{equation}
The canonical target minimizes holding and backlog cost at its fixed width: otherwise shifting both endpoints would reduce total cost without changing the order frequency. Together with the canonical average-cost identity~\eqref{eq:canonical-average-cost}, this proves the last equality in~\eqref{eq:smooth-data}. By~\eqref{eq:mirror}, $m_{\hat w}$ is the mirror cost $B_{\hat w}$. We bound this cost through the canonical width: the endpoint argument in the proof of Lemma~\ref{lem:reward}, with $\beta_H=Hp/(H+p)$, gives $\hat w\le j/\beta_H$. Consequently, Lemma~\ref{lem:smooth-mirror} yields
\begin{equation}\label{eq:smooth-mirror-ratio}
 \frac{m_{\hat w}}{j}\le\frac{p(1-\sqrt{p/(p+h)})\hat w}{j}\le\frac{p(1-\sqrt{p/(p+h)})}{\beta_H}
 \le\frac1{1+\sqrt{p/(p+h)}}\le\lambda,
\end{equation}
where the second inequality uses $\hat w\le j/\beta_H$, and $H\ge h$ gives $\beta_H\ge hp/(h+p)$ for the third. For the last, $p\ge h$ gives $p/(p+h)\ge1/2$, and $1/(1+1/\sqrt2)=2-\sqrt2=\lambda$; neither step compares $p$ and $H$.

\emph{Step 2. Choose the window length.} Let $t_0=U\{0\}/M(\hat w)$. Continuity of the depletion distribution away from zero makes every time fraction $t\in(t_0,1]$ attainable by a cutoff. For such a fraction, the window $[0,w_t)$ with $U[0,w_t)/M(\hat w)=t$ has order frequency $\nu/t$. By~\eqref{eq:occupation}, its depletion has distribution function $\Prob(Y\le x)/t$ on $[0,w_t)$, so it is $Y^{[t]}$ in Lemma~\ref{lem:quantile}, with $Q$ the quantile function of $Y$. Retune its retailer and mirror targets to minimize the two losses. The mirror contraction~\eqref{eq:mirror-contract} and the outer inequality of~\eqref{eq:running-contract} bound its actual augmented cost by $r_{\hat w}+tm_{\hat w}+z_{\hat w}/t$.\footnote{The three terms bound holding and backlog cost, mirror cost, and actual setup cost, respectively. After retuning the target, holding and backlog cost is $\NV_{H,p}(Z+Y^{[t]})\le\NV_{H,p}(Z+Y)=r_{\hat w}$, by comparing the first and last expressions in~\eqref{eq:running-contract}; this is its outer inequality. After retuning the mirror center, mirror cost is $\NV_{p,h}(Y^{[t]})\le t\NV_{p,h}(Y)=tm_{\hat w}$ by~\eqref{eq:mirror-contract}. Finally, $M(w_t)=tM(\hat w)$, so the order frequency is $1/M(w_t)=\nu/t$. Charging the actual setup cost $K$ per order gives $K\nu/t=z_{\hat w}/t$, with $z_{\hat w}=K\nu$ from~\eqref{eq:smooth-data}, which is the third term.}

The preceding cost bound holds for every fraction $t\in(t_0,1]$. We now take the fraction prescribed in Definition~\ref{def:dfr}; Step 3 treats a selected fraction $t\le t_0$ by base stock. When $0<z_{\hat w}<m_{\hat w}$, this choice is $t=\sqrt{z_{\hat w}/m_{\hat w}}$. Substituting this value and using $(\sqrt{z_{\hat w}/\theta}-\sqrt{\theta m_{\hat w}})^2\ge0$ gives $r_{\hat w}+tm_{\hat w}+z_{\hat w}/t=r_{\hat w}+2\sqrt{z_{\hat w}m_{\hat w}}\le j+\theta m_{\hat w}\le(1+\theta\lambda)j$, where the first inequality also uses $j=r_{\hat w}+z_{\hat w}/\theta$ from~\eqref{eq:smooth-data}, and the second uses $m_{\hat w}\le\lambda j$ from~\eqref{eq:smooth-mirror-ratio}. When $z_{\hat w}\ge m_{\hat w}$, including $m_{\hat w}=0$, we retain the full cycle: subtracting its cost $r_{\hat w}+m_{\hat w}+z_{\hat w}$ from $j+\theta m_{\hat w}$ gives $(1-\theta)(z_{\hat w}/\theta-m_{\hat w})\ge0$, because $z_{\hat w}/\theta\ge z_{\hat w}\ge m_{\hat w}$, so the full cycle's actual cost is then also at most $(1+\theta\lambda)j$.

\emph{Step 3. Handle the atom at zero and zero setup costs.} Since $U\{0\}=1/q$, the time fraction at depletion zero is $t_0=\nu/q$. If the selected $t$ is at most $t_0$, use base stock at a minimizer of $F$, the expected holding and backlog cost against the lead-time demand $Z$. Its mirror loss is zero and its cost is $\NV_{H,p}(Z)+Kq=\NV_{H,p}(Z)+z_{\hat w}/t_0\le r_{\hat w}+tm_{\hat w}+z_{\hat w}/t$, because $t\le t_0$ and conditional Jensen's inequality gives $\NV_{H,p}(Z)\le r_{\hat w}$.\footnote{\label{fn:nv-conditional-jensen}For every target $\Gamma$, conditioning on $Z$ and using the independence of $Y$ and $Z$ gives $\E[H(\Gamma-Z-Y)^++p(Z+Y-\Gamma)^+\mid Z]\ge H(\Gamma-Z-\E Y)^++p(Z+\E Y-\Gamma)^+$, since the loss is convex in $Y$. Taking expectations and minimizing over $\Gamma$ yields $r_{\hat w}=\NV_{H,p}(Z+Y)\ge\NV_{H,p}(Z+\E Y)=\NV_{H,p}(Z)$, where the first equality is~\eqref{eq:smooth-data} and the last follows by replacing $\Gamma$ with $\Gamma-\E Y$ in the minimization, which shifts the target without changing the minimum.} Hence base stock also meets the bound of Step~2.

The policy choices used above are well defined: the retailer order-up-to target and mirror center can each be chosen at a finite value that minimizes the corresponding expected loss, by the continuity and coercivity argument for the newsvendor loss in Section~\ref{sec:arb-foundations}. Finally, a retailer with $K_i=0$ uses base stock, with zero mirror cost, and meets the same factor $1+\theta\lambda$.

\emph{Step 4. Combine the network bounds.} We have, therefore, constructed deterministic retailer policies with bounded post-action positions with $\cL_i\le(1+\theta\lambda)J_i(K_i/\theta)$. The policy $\pi_{\rm DFR}(\theta)$ of Definition~\ref{def:dfr} implements these retailer policies by the FIFO allocation, supplier rule, and initialization of Definition~\ref{def:fifo}, using their optimized mirror centers and fixed rules for all ties. A retailer with $K_i>0$ that selects base stock has bounded per-period setup charges, at most $K_i$, and bounded position costs. Thus, the proof of Proposition~\ref{prop:upper}, including the initialization and expected-average argument of Appendix~\ref{app:initialization}, applies to these policies and centers, so $C(\pi_{\rm DFR}(\theta))\le J_0(K_0+K_{\max})+\sum_i\cL_i$. For the matching lower bound, substituting $\gamma=1/(1-\theta)$ into~\eqref{eq:lower} and dividing by $\gamma$ gives
\begin{equation}\label{eq:smooth-weighted-lower}
 (1-\theta)J_0\bigl(K_0/(1-\theta)\bigr)
       +\theta\sum_iJ_i(K_i/\theta)\le\OPT.
\end{equation}
If $K_{\max}\le\theta K_0/(1-\theta)$, then $K_0+K_{\max}\le K_0/(1-\theta)$. Since the optimal supplier cost $J_0(k)$ is nondecreasing in its setup cost $k$, this gives
\begin{equation}\label{eq:smooth-supplier-comparison}
 J_0(K_0+K_{\max})\le J_0\bigl(K_0/(1-\theta)\bigr),
\end{equation}
where the left side is the supplier term in the FIFO upper bound and the right side is the supplier value used in~\eqref{eq:smooth-weighted-lower}. Combining this comparison with $\cL_i\le(1+\theta\lambda)J_i(K_i/\theta)=(1/\theta+\lambda)\,\theta J_i(K_i/\theta)$ gives
$C(\pi_{\rm DFR}(\theta))\le J_0(K_0/(1-\theta))+(1/\theta+\lambda)\,\theta\sum_iJ_i(K_i/\theta)$. The supplier term in this upper bound is $1/(1-\theta)$ times its weighted contribution in~\eqref{eq:smooth-weighted-lower}, and the retailer terms are bounded by $1/\theta+\lambda$ times their weighted contribution. Both contributions are nonnegative, so the larger of these two factors bounds their sum relative to~\eqref{eq:smooth-weighted-lower}. Thus, under the stated setup-cost condition,
\begin{equation}\label{eq:smooth-network-factor}
 C(\pi_{\rm DFR}(\theta))\le
 \max\left\{\frac1{1-\theta},\frac1\theta+\lambda\right\}\OPT.
\end{equation}

\emph{Step 5. Balance the approximation factors.} For the simple choice $\theta=7/12$, the two factors in~\eqref{eq:smooth-network-factor} satisfy $1/(1-\theta)=12/5$ and $1/\theta+\lambda=12/7+\lambda<81/35<12/5$, where $\lambda=2-\sqrt2<3/5$ because $(7/5)^2<2$. The larger factor is therefore $12/5$. At this choice of $\theta$, the theorem's setup-cost assumption gives $K_{\max}\le(\gamma_\lambda-1)K_0\le(7/5)K_0=\theta K_0/(1-\theta)$, so the supplier cost comparison~\eqref{eq:smooth-supplier-comparison} applies. Thus,~\eqref{eq:smooth-network-factor} gives $C(\pi_{\rm DFR}(7/12))\le(12/5)\OPT$, as stated after Theorem~\ref{thm:dfr}.

To improve this factor, we now vary $\theta$. The supplier factor $1/(1-\theta)$ increases with $\theta$ and the retailer factor $1/\theta+\lambda$ decreases, so their intersection minimizes the maximum in~\eqref{eq:smooth-network-factor} on $[1/2,1)$. Setting both terms equal to $\gamma$, so that $\theta=1-1/\gamma$ and $1/\theta=\gamma/(\gamma-1)$, gives $\gamma^2-(2+\lambda)\gamma+\lambda=0$; its larger root is the factor $\gamma_\lambda=1+\lambda/2+\sqrt{1+\lambda^2/4}$. Since $\gamma_\lambda\ge2$, this $\theta=1-1/\gamma_\lambda$ lies in $[1/2,1)$. Here $\theta/(1-\theta)=\gamma_\lambda-1$, so~\eqref{eq:smooth-network-factor} proves the claimed bound for $\pi_{\rm DFR}(\theta)$.\Halmos\endproof

\subsection{Proofs of the Contraction and Mirror Bounds}\label{app:smooth-proofs}

\proof{Proof of Lemma~\ref{lem:quantile}.} We first show that every gap between two quantiles of $Y^{[t]}$ is at most the corresponding gap for $tY$. We then subtract each variable's own mean and compare the resulting variables in convex order: every convex loss has no larger expectation for $Y^{[t]}-\E Y^{[t]}$ than for $t(Y-\E Y)$, whenever these expectations are defined. Applying it to inventory losses and using convexity gives both contraction inequalities, for the mirror loss and the holding and backlog loss.

\emph{Step 1. Compare gaps between quantiles.} For $0<u<v<1$, the gap between the $u$ and $v$ quantiles is $Q(tv)-Q(tu)$ for $Y^{[t]}$ and $t\{Q(v)-Q(u)\}$ for $tY$. Since $Q$ is nondecreasing and convex, its almost-everywhere derivative $Q'$ is nonnegative and nondecreasing, so $Q'(tx)\le Q'(x)$. Integrating on compact subintervals, then taking endpoint limits, gives $Q(tv)-Q(tu)=t\int_u^vQ'(tx)\,dx\le t\int_u^vQ'(x)\,dx=t\{Q(v)-Q(u)\}$, so each gap of $Y^{[t]}$ is at most the corresponding gap of $tY$.

\emph{Step 2. Establish convex order after centering.} Hence, $tQ(u)-Q(tu)$ is nondecreasing, including across any initial zero segment. This difference governs the variability of the centered depletion. If $Q_1,Q_2$ are nondecreasing integrable functions and $Q_2-Q_1$ is nondecreasing, then its average over $[0,a]$ is at most its average over $[0,1]$: $\int_0^a\{Q_2(u)-Q_1(u)\}\,du\le a\int_0^1\{Q_2(u)-Q_1(u)\}\,du$ for $0\le a\le1$. After subtracting each quantile function's mean, the preceding inequality says that the lower quantile integrals of $Q_1$ dominate those of $Q_2$. For a uniform random variable $V$ on $(0,1)$ and $j=1,2$, applying the identity $\E(s-Q_j(V))^+=\max_{0\le a\le1}\{as-\int_0^a Q_j(u)\,du\}$ to the centered quantile functions orders all lower stop-loss expectations\footnote{For an integrable random variable $X$, the lower stop-loss expectation at threshold $s$ is $\E(s-X)^+$, the expected amount by which $X$ falls below $s$, with zero contribution when $X\ge s$. See \citet[Section~3.A.1, equation~(3.A.6)]{ShakedShanthikumar2007} for these expected losses and their role in convex order. Here the centered variables are $X_j=Q_j(V)-\E Q_j(V)$, and ``orders all'' means $\E(s-X_1)^+\le\E(s-X_2)^+$ for every real $s$.} of the centered variables. Write $X_1\le_{\rm cx}X_2$ when the expectation of every convex function of $X_1$ is at most that of the same function of $X_2$, whenever both are defined. Equal means and ordered lower stop-loss transforms imply this convex order \citep[Theorem~3.A.1]{ShakedShanthikumar2007}. Applying this characterization to $Q_1(u)=Q(tu)$ and $Q_2(u)=tQ(u)$, whose difference is nondecreasing, gives $Y^{[t]}-\E Y^{[t]}\le_{\rm cx}t(Y-\E Y)$.

\emph{Step 3. Compare the optimized losses.} Adding an independent $Z$ preserves this order, since averaging any convex test function over $Z$ leaves a convex function of its remaining argument. For each fixed target, the inventory loss in~\eqref{eq:nv} is convex, so minimizing over the target and using shift invariance give $\NV_{a,b}(Z+Y^{[t]})\le\NV_{a,b}(Z+tY)$. Taking $Z=0$ gives $\NV_{a,b}(Y^{[t]})\le\NV_{a,b}(tY)$, and positive homogeneity, $\NV_{a,b}(tY)=t\NV_{a,b}(Y)$, proves the mirror contraction~\eqref{eq:mirror-contract}.

For the holding and backlog comparison, we also use convexity of $\NV_{a,b}$ on a common probability space: combine two feasible targets in the same proportions as their random variables and use convexity of the positive part. The comparison just proved gives the first inequality, and this convexity applied to $Z+tY=(1-t)Z+t(Z+Y)$ gives the second: $\NV_{a,b}(Z+Y^{[t]})\le\NV_{a,b}(Z+tY)\le(1-t)\NV_{a,b}(Z)+t\NV_{a,b}(Z+Y)$. Conditional Jensen's inequality (footnote~\ref{fn:nv-conditional-jensen}) gives $\NV_{a,b}(Z)\le\NV_{a,b}(Z+Y)$, proving~\eqref{eq:running-contract}. Theorem~\ref{thm:dfr} uses it to bound the holding and backlog cost of a shorter window by $r_{\hat w}$; finite first moments make each optimized loss finite.\Halmos\endproof

\proof{Proof of Lemma~\ref{lem:smooth-mirror}.} We represent the depletion distribution as a mixture of uniform distributions. Evaluating every uniform component at one common mirror target gives a cost bound that also holds after averaging over the mixture, which proves the mirror bound~\eqref{eq:smooth-mirror}.

\emph{Step 1. Represent the distribution as a uniform mixture.} The case $w=0$ is immediate, so assume $w>0$. A distribution with mass $t_0$ at zero and a nonincreasing density $\varrho$ on $(0,w)$ is a mixture of uniform distributions on $[0,a]$, $0\le a\le w$, and the mass at zero is the component $a=0$.\footnote{Take a right-continuous version of $\varrho$, place mass $t_0$ at $a=0$, use mixing measure $-a\,d\varrho(a)$ on $(0,w)$, and place mass $w\varrho(w-)$ at $w$. As $a\downarrow0$, $a\varrho(a)\le\int_0^a\varrho(x)\,dx\to0$ because $\varrho$ is nonincreasing and integrable. Integration by parts gives total mass one, and integrating $1/a$ over components with $a>x$ gives $\varrho(x)$ almost everywhere. The argument allows jumps in $\varrho$ and a density that is unbounded near zero, so no further regularity is needed.}

\emph{Step 2. Bound each component at one target.} Choose the common mirror target $m=w(1-\sqrt{p/(p+h)})$. For a uniform component with $a\le m$, the expected loss is $p(m-a/2)\le pm$. For $a\ge m$, direct integration gives $\mathcal L(a)=ha/2-hm+(p+h)m^2/(2a)$, $\mathcal L(m)=pm/2$, and $\mathcal L(w)=pm$. The function $\mathcal L$ is convex on $[m,w]$, as its last term is convex in $a$ and the others are affine, so its maximum occurs at an endpoint and is at most $pm$. Averaging over the mixture, the expected loss at the fixed target $m$ is at most $pm$, and minimizing over targets gives~\eqref{eq:smooth-mirror}.\Halmos\endproof

\section{Details of the Numerical Study}\label{app:numerics}

\subsection{Items and an Example Instance}

We select an item for a state if, over the last 730 days, every store in the state records a sale within the first 28 days and has mean daily sales between 0.5 and 40. Among them, network-average store sales define high ($\ge10$), medium ($3$ to below $10$), and low ($0.5$ to below $3$) tiers. We keep the two best-selling eligible items in each state and tier, giving the 18 state-item pairs in Table~\ref{tab:numerics-data}, with 16 distinct items, as \texttt{FOODS\_3\_587} is kept in all three states. The empirical distributions of Section~\ref{sec:numerics} ignore weekday and seasonal patterns, which the replay in Appendix~\ref{app:numerics-checks} keeps.

\begin{table}[htbp]
\SingleSpaced
\centering
\caption{Selected items and store-level sales, zero-sales frequencies, and average unit prices.}
\label{tab:numerics-data}
\footnotesize\setlength{\tabcolsep}{5pt}
\begin{tabular}{@{}lllr rr@{}}
\toprule
State (stores) & Tier & Item & Mean sales & Zero days (\%) & Price (\$)\\
\midrule
CA (4) & High & \texttt{FOODS\_3\_282} & 9.1--38.2 & 4--8 & 2.51\\
 & High & \texttt{FOODS\_3\_587} & 12.7--33.1 & 3--11 & 2.47\\
 & Medium & \texttt{FOODS\_3\_389} & 2.4--15.6 & 0--38 & 1.67\\
 & Medium & \texttt{FOODS\_3\_406} & 4.3--14.5 & 9--19 & 0.80\\
 & Low & \texttt{FOODS\_3\_121} & 1.1--5.1 & 4--43 & 2.98\\
 & Low & \texttt{FOODS\_3\_769} & 1.9--4.0 & 32--44 & 1.29\\
\addlinespace
TX (3) & High & \texttt{FOODS\_3\_587} & 29.8--34.3 & 1--2 & 2.27\\
 & High & \texttt{FOODS\_3\_714} & 18.2--24.5 & 0--0 & 1.58\\
 & Medium & \texttt{HOUSEHOLD\_1\_459} & 8.0--10.8 & 5--13 & 0.96\\
 & Medium & \texttt{FOODS\_3\_234} & 6.4--12.6 & 27--50 & 0.20\\
 & Low & \texttt{HOUSEHOLD\_1\_457} & 2.4--3.5 & 9--14 & 10.10\\
 & Low & \texttt{HOBBIES\_1\_404} & 2.0--3.9 & 34--51 & 0.47\\
\addlinespace
WI (3) & High & \texttt{FOODS\_3\_586} & 14.7--33.3 & 0--0 & 1.67\\
 & High & \texttt{FOODS\_3\_587} & 11.5--25.1 & 12--16 & 2.48\\
 & Medium & \texttt{FOODS\_3\_458} & 5.9--13.2 & 0--2 & 1.94\\
 & Medium & \texttt{FOODS\_3\_329} & 6.4--14.1 & 39--40 & 2.98\\
 & Low & \texttt{FOODS\_2\_099} & 1.5--4.1 & 14--40 & 1.99\\
 & Low & \texttt{FOODS\_1\_162} & 2.0--4.5 & 21--28 & 1.93\\
\bottomrule
\end{tabular}\par\vspace{6pt}
\begin{minipage}{.98\linewidth}
\footnotesize\emph{Notes.} CA, TX, and WI denote California, Texas, and Wisconsin. Mean sales and Zero days give ranges across the state's stores of mean daily unit sales and of the percentage of days without sales. Price is the average store unit price in dollars. All summaries use the final 730 days; tiers use the network-average store mean.
\end{minipage}
\end{table}

For California's \texttt{FOODS\_3\_587} at $L_0=3$, service target $0.90$, and $(T_s,T_w)=(3,7)$, daily unit rates $H_i,h,p$ are about $0.17$, $0.13$, and $1.52$ cents, retailer setup costs are $10$ to $25$ cents, and $K_0=\$2.64$. Its simulated daily cost of $\$1.28$ is $1.50$ times the best evaluated lower bound of $\$0.86$ and $0.95$ times the upper bound of $\$1.36$ of Section~\ref{sec:numerics}, which exceeds it by about 6 percent.

\subsection{Single-Location Computation and Benchmark Searches}\label{app:numerics-methods}

\emph{Single-location values and search ranges.} The empirical demands are integer-valued, so the renewal measure is supported on nonnegative integers. Replacing a positive width $w$ by $\lceil w\rceil$ leaves its depletion distribution and order frequency unchanged. For a fixed integer width $w$, write $u_j=U\{j\}$. Its single-location cost at target $S$ is
\[
 \frac{k+\sum_{j=0}^{w-1}u_j f(S-j)}{\sum_{j=0}^{w-1}u_j},
\]
where $f$ is the retailer or supplier expected loss. Lead-time demand is integer-valued, so $f$, and hence this objective as a function of $S$, is convex and piecewise linear with integer breakpoints. It therefore has an integer minimizer. If $x_*$ is an integer minimizer of $f$, some minimizing target lies in $[x_*,x_*+w-1]$: below this interval all arguments of $f$ are at or below $x_*$, and above it all are at or above $x_*$. Moving a target toward the interval cannot increase cost. Thus integer widths and targets recover the unrestricted single-location optimum; they do not impose an integer-action restriction on the lower bound, which therefore remains a bound for the divisible-quantity model.

The width search range also provably contains the canonical width. Let $a,b$ be the holding and penalty rates of the single-location problem, and let $\widetilde J(k)$ be a feasible value from the current search. The canonical endpoint bound used in Lemma~\ref{lem:reward} gives
\[
 \hat w\le J(k)(1/a+1/b)\le\widetilde J(k)(1/a+1/b).
\]
The implementation enlarges its maximum width $w_{\max}$ until the last expression is at most $w_{\max}-2$, leaving room for rounding the canonical width upward. Targets are enumerated from $x_*$ through $x_*+w_{\max}$. Tables shared across setup costs are sized at the largest required setup cost, and the width bound is checked for every positive setup cost queried. At zero setup cost, the value is $\min f$. Mirror losses likewise have an integer minimizing center. These floating-point calculations have no grid-discretization or width-truncation error. The separate rounding of the half-window candidate in Section~\ref{sec:numerics-instances} changes the implemented policy, whose own augmented cost enters $\mathrm{UB}$.

\emph{Simulation searches.} Each benchmark uses two starting policies and a separate 50,000-day search path, discarding its first 5,000 days. The two starting $(R,nQ)$ policies come from single-location optima at actual setup costs and from translating our implemented FIFO parameters into reorder points and batch sizes. Coordinate search adjusts reorder points and positive integer batch sizes in steps of $1,3,10,30,$ and $100$, testing both signs and retaining improvements. It dispatches only complete queued batches, in FIFO order, whereas our policy may fill a request in part.

For the integer-ratio benchmark, review intervals lie between 1 and 56 days, and each retailer interval is an integer multiple or divisor of the supplier interval. Its starting policies round the calibrating EOQ intervals and our implemented policy's width-to-mean-demand ratios to integer ratios. Initial order-up-to levels are newsvendor quantiles over each review interval plus its lead time. Search moves to neighboring feasible retailer intervals or shifts the supplier interval by one day, mapping retailer intervals to the nearest feasible values and adjusting targets by the associated newsvendor-level changes. It also moves targets by $1,3,10,30,$ and $100$ units in either direction. Retailer review phases match supplier receipts, and dispatches use FIFO with partial fills.

Each search stops after a pass without improvement or at a nominal budget of 2,500 trial evaluations per search; no reported run reached that budget. Summed over the two searches, the mean counts are 1,011 trials per instance for $(R,nQ)$ and 1,067 for integer-ratio policies. We keep the resulting policy with the lower search-path cost and evaluate it on the common 100,000-day evaluation path after 10,000 warm-up days. These are the best policies found by the specified searches, rather than proven optima within their policy classes. In contrast, our FIFO implementation sets its parameters by single-location calculations, without any simulation search.

\subsection{Results by Group}\label{app:numerics-groups}

Table~\ref{tab:numerics-main} reports the common-penalty results of Section~\ref{sec:numerics} by supplier lead time and setup calibration. Means and maxima are within groups, and the last row pools all 864 instances.

\begin{table}[htbp]
\SingleSpaced
\centering
\caption{Common-penalty performance by supplier lead time and setup-cost calibration.}
\label{tab:numerics-main}
\footnotesize\setlength{\tabcolsep}{4pt}
\begin{tabular}{@{}lrrrrrrrr@{}}
\toprule
 & & \multicolumn{2}{c}{$C_{\rm sim}/\mathrm{LB}^*$} & \multicolumn{2}{c}{$\mathrm{UB}/\mathrm{LB}^*$} & $\mathrm{UB}/C_{\rm sim}$ & Uncoord. & Full only\\
Group & Count & Mean & Max. & Mean & Max. & Mean & Mean & Mean\\
\midrule
\multicolumn{9}{l}{Panel A: supplier lead time $L_0$}\\
1 & 216 & 1.46 & 1.65 & 1.53 & 1.82 & 1.047 & 1.48 & 1.54\\
3 & 216 & 1.45 & 1.65 & 1.54 & 1.84 & 1.053 & 1.48 & 1.53\\
7 & 216 & 1.45 & 1.70 & 1.54 & 1.87 & 1.059 & 1.47 & 1.53\\
14 & 216 & 1.44 & 1.75 & 1.54 & 1.93 & 1.064 & 1.46 & 1.51\\
\addlinespace
\multicolumn{9}{l}{Panel B: calibrating EOQ cycles $(T_s,T_w)$}\\
$(3,7)$ & 144 & 1.55 & 1.72 & 1.65 & 1.87 & 1.064 & 1.58 & 1.63\\
$(3,14)$ & 144 & 1.38 & 1.64 & 1.44 & 1.72 & 1.042 & 1.39 & 1.43\\
$(3,28)$ & 144 & 1.21 & 1.45 & 1.24 & 1.49 & 1.030 & 1.22 & 1.23\\
$(7,7)$ & 144 & 1.64 & 1.75 & 1.80 & 1.93 & 1.101 & 1.67 & 1.78\\
$(7,14)$ & 144 & 1.56 & 1.66 & 1.65 & 1.80 & 1.060 & 1.59 & 1.66\\
$(7,28)$ & 144 & 1.37 & 1.54 & 1.42 & 1.59 & 1.036 & 1.39 & 1.43\\
\midrule
All & 864 & 1.45 & 1.75 & 1.53 & 1.93 & 1.056 & 1.47 & 1.53\\
\bottomrule
\end{tabular}\par\vspace{6pt}
\begin{minipage}{.98\linewidth}
\footnotesize\emph{Notes.} Lead times and calibrating EOQ cycles are in days. $C_{\rm sim}$ is our policy's simulated cost, $\mathrm{UB}$ is its analytical upper bound, and $\mathrm{LB}^*$ is the largest evaluated shared-cost lower bound. Uncoord. and Full only give the mean ratio of simulated cost to $\mathrm{LB}^*$ for the uncoordinated $(s,S)$ policy and for the variant that uses only full windows. Ratios use two decimals, except the three-decimal $\mathrm{UB}/C_{\rm sim}$ column. Count is the number of instances in each group.
\end{minipage}
\end{table}

\subsection{Further Checks}\label{app:numerics-checks}

Three further checks vary the retailer lead time, network size, and demand process. First, on the subset with service target $0.90$ and $(T_s,T_w)=(3,7)$, retailer lead times $0$, $1$, and $2$ give mean ratios $C_{\rm sim}/\mathrm{LB}^*$ of $1.48$, $1.53$, and $1.55$. Second, a national instance with all ten stores and item \texttt{FOODS\_3\_714} uses a common assumed unit price of $\$1.58$; over four supplier lead times and both service targets at $(T_s,T_w)=(3,7)$, simulated ratios $C_{\rm sim}/\mathrm{LB}^*$ range from $1.49$ to $1.56$ and upper-bound ratios are below $1.60$. Third, we replay actual sales, which keep their weekday and seasonal patterns and shift from one year to the next. The replay designs our policy from the sales of days 1,212 to 1,576 and operates it on actual sales from days 1,577 to 1,941, using the 18 state-item pairs at $L_0=3$, service target $0.90$, and $(T_s,T_w)=(3,7)$. We do not use the lower bound in this replay comparison, because its validity relies on demand being i.i.d.\ over periods and it does not bound optimal cost on a fixed sales sequence. The design-period $\mathrm{UB}$ bounds expected cost under i.i.d.\ demand and gives no pathwise guarantee, so we use it as a benchmark; in i.i.d.\ simulation of these instances, cost is $0.92$ to $0.95$ times $\mathrm{UB}$. After 30 warm-up days, the 335-day average cost is $1.03$ times the design-period $\mathrm{UB}$ on average (median $1.00$, range $0.79$ to $1.38$). The spread across items mainly reflects changes in demand between the two years: recomputed from the replay year's empirical distributions, with all cost parameters held at their original values, $\mathrm{UB}$ ranges from $0.77$ to $1.25$ times its design-period value, and this factor has correlation $0.90$ with the cost ratio across items. On the same sales, the policy designed with hindsight from the replay year's distributions has realized cost close to its recomputed i.i.d.\ upper bound: it costs $0.98$ times that $\mathrm{UB}$ on average and exceeds it for 3 of the 18 items. This comparison mixes the weekday and seasonal patterns of actual sales with single-path variation and slack in the analytical bound, so it does not measure the effect of these patterns alone. Our policy's realized cost is $1.07$ times the hindsight policy's on average, ranging from $1.02$ to $1.18$. The latter is the observed premium relative to a hindsight-refitted policy in the same class, on these items and paths with all cost parameters held fixed.

In the integer-ratio comparison of Section~\ref{sec:numerics}, our policy's estimated cost is lower in all 20 instances of the five intermittent items at $(3,28)$, by $2.9$ percent on average and up to $7.2$ percent, and in 8 and 3 of their 20 instances at $(3,7)$ and $(7,7)$. The batch-means interval of Section~\ref{sec:numerics} for our policy's cost minus the benchmark's, computed from the ten batch differences on the common demand path, lies below zero in 27 of these 31 instances. It includes zero in the other four, where our policy is cheaper by $0.3$ to $0.7$ percent, and in two instances of Wisconsin's \texttt{FOODS\_2\_099} at $(3,28)$, where it costs at most $0.15$ percent more, a difference within simulation error.

\subsection{Dynamic Programming for Two Stores}\label{app:numerics-dp}

A dynamic program for an $n$-store network tracks warehouse stock, each store's inventory position, and the $L_0-1$ supplier orders in transit, giving $n+L_0$ state variables when $L_0\ge1$. A four-store network at $L_0=14$ has 18 state variables, so even 20 values per variable give more than $10^{23}$ states. Therefore, we solve bounded-state dynamic programs only for two-store networks. For each of the six low-tier items, we keep the network's two stores with the smallest mean daily sales and apply the calibration of Section~\ref{sec:numerics-instances} to this two-store network, with $L_0\in\{0,1\}$, $L_i=1$, both service targets, and all six pairs $(T_s,T_w)$. These 144 instances satisfy the conditions for the $5/2$ specialization of Theorem~\ref{thm:main}, with $K_{\max}/K_0\le0.93<3/2$; our largest program has 5.2 million states.

The computation first restricts inventories, orders, and dispatches to integers. We next express running cost in terms of current inventory positions. Retailer $i$'s conditional expected holding and backlog cost in period $t+L_i$ is $F_i(P_{i,t})$ by~\eqref{eq:retailer-lead} and~\eqref{eq:retailer-value}, so charging it in period $t$ leaves the expected average cost unchanged. With $L_0\le1$, no supplier order is in transit at the start of a period, so warehouse stock and the two store inventory positions form a sufficient state. We then restrict these states to a finite box with bounds chosen from single-location $(s,S)$ levels; a store position that falls below the box after demand is reset to the box's lower boundary.

For this finite model, relative value iteration stops when its lower and upper bounds on the optimal average cost of this truncated integer model differ by at most $0.002$ percent of that cost; write $C_{\rm DP}$ for their midpoint. This tolerance controls only the iteration error; the calculation gives no error bound for the state truncation or the integer restriction. Therefore, we use $C_{\rm DP}$ as a numerical benchmark for $\OPT$ of the divisible-quantity model, not as a bound on it.

To assess the effect of the finite box, we simulate each computed DP policy on the untruncated system for 400,000 days after 10,000 warm-up days: it leaves the box on at most two days, and its estimated cost is within $0.5$ percent of $C_{\rm DP}$. On the same demand paths, the ratio of the FIFO policy's simulated cost to $C_{\rm DP}$ averages $1.22$ and ranges from $1.08$ to $1.45$. To assess the strength of the lower bound in these instances, we also compare it with the DP benchmark: $\mathrm{LB}^*$ averages $0.82$ times $C_{\rm DP}$ and is at least $0.71$ times it, and $C_{\rm sim}/\mathrm{LB}^*$ averages $1.49$. The policy's ratio to $C_{\rm DP}$ averages $1.23$ at $L_0=0$ and $1.22$ at $L_0=1$ and rises with store setup costs relative to the supplier's, which scale with $(T_s/T_w)^2$, from a mean of $1.10$ at $(T_s,T_w)=(3,28)$ to $1.38$ at $(7,7)$.

At $(7,7)$, where this ratio is largest, the computed DP policies leave the warehouse empty after dispatch in 86 to 99 percent of simulated periods and send 74 to 97 percent of store shipments in periods with supplier receipts. Our policy keeps warehouse stock to serve half-window store cycles: its store setup costs are $1.75$ to $1.93$ times those of the computed DP policies, and warehouse holding is 24 to 32 percent of its simulated cost, against at most 3 percent of theirs.

\end{document}